\documentclass[11pt,a4paper]{article}
\usepackage[]{times}
\usepackage{fullpage}
\usepackage{graphicx}
\usepackage{subfigure}
\usepackage{amsmath}
\usepackage{fancyhdr}
\usepackage[colorlinks]{hyperref}
\usepackage{xcolor}
\date{}

\newcommand{\prov}{{\sc Proof}.\hspace*{3mm} }

\newcommand{\QED}{$\rule{2mm}{2mm}$}

\newcommand{\natu}{{\sf I \! N}}
\newtheorem{theorem}{Theorem}[section]
\newtheorem{lemma}[theorem]{Lemma}
\newtheorem{e-proposition}[theorem]{Proposition}
\newtheorem{corollary}[theorem]{Corollary}
\newtheorem{e-definition}[theorem]{Definition\rm}
\newtheorem{remark}{\it Remark\/}

\title{Methods of arbitrary optimal order with tetrahedral finite-element meshes forming polyhedral approximations of curved domains}
\author{
    Vitoriano Ruas$^{1}$\thanks{This work was partially supported by CNPq, the National Research Council of Brazil}
		\\[1mm]
  {\small $^{1}$ Institut Jean Le Rond d'Alembert, CNRS UMR 7190, 
  Sorbonne Universit\'e, F-75005 Paris, France.}\\[1mm]
  {\small e-mail: {\it vitoriano.ruas@upmc.fr}}}

\begin{document}
\maketitle

\begin{abstract}

In recent papers (see e.g. \cite{arXiv}, \cite{Maugin}, \cite{CAMWA} and \cite{PAMM}) the author introduced a simple alternative of the $n$-simplex type, to enhance the accuracy of approximations of second-order boundary value problems with Dirichlet conditions, posed in 
smooth curved domains. This technique is based upon trial-functions consisting of piecewise polynomials defined on straight-edged triangular or tetrahedral meshes, interpolating the Dirichlet boundary conditions at points of the true boundary. In contrast the test-functions are defined upon the standard degrees of freedom associated with the underlying method for polytopic domains. While method's mathematical analysis for two-dimensional domains was carried out in \cite{arXiv} and \cite{PAMM}, this paper is devoted to the study of the three-dimensional case. Well-posedness, uniform stability and optimal a priori error estimates in the energy norm are demonstrated for a tetrahedron-based Lagrange family of finite elements. Unprecedented $L^2$-error estimates for the class of problems considered in this work are also proved. A series of numerical examples illustrates the potential of the new technique. In particular its better accuracy at equivalent cost as compared to the isoparametric technique is highlighted. Moreover the great generality of the new approach is exemplified through a method with degrees of freedom other than nodal values.
\end{abstract}


\section{Introduction}
Petrov-Galerkin formulations of boundary value problems showed in the past decades to be a powerful tool to overcome difficulties brought about by the space discretization of certain types of partial differential equations. A significant illustration is provided by the SUPG method introduced by Hughes \& Brooks \cite{HughesBrooks} in 1982, in order to stably handle convection-diffusion equations. Other examples are the families of methods proposed by Hughes \& Franca and collaborators in the late eighties for the finite-element modeling of various  problems in Continuum Mechanics, in particular as a popular alternative to Galerkin methods for viscous incompressible flow (see e. g. \cite{HughesFranca}). The outstanding contributions about ten years earlier of Babu\v{s}ka (see e.g. \cite{Babuska}) and Brezzi \cite{Brezzi}, among other authors, were decisive to provide the theoretical background that allowed to formally justify the reliability of Petrov-Galerkin formulations, namely, the so-called \textit{inf-sup} condition. In this paper we endeavor to show another application of this approach in a rather different framework, though not less important.\\    
More precisely this work deals with finite element methods of optimal order greater than one 
to solve boundary value problems with Dirichlet conditions, posed in domains with a smooth curved boundary of arbitrary shape. The method 
is similar to the technique known as \textit{interpolated boundary conditions}, or simply IBC,  
studied in \cite{BrennerScott}. However, in spite of being very intuitive and known since the seventies (cf. \cite{Nitsche} and \cite{Scott}) IBC has not been much used so far. This is certainly due to its difficult implementation, the lack of an extension to three-dimensional problems and, most of all, restrictions on the choice of boundary nodal points to reach optimal convergence rates. In contrast the  implementation of our method is straightforward in both two- and three-dimensional geometries. This is due to the fact that only polynomial algebra is necessary, while the domain is simply approximated by the polytope formed by the union of standard $n$-simplexes of a finite-element mesh. Furthermore approximations of optimal order can be obtained for non-restrictive choices of boundary nodal points. \\
\indent Generally speaking, our methodology is designed to handle Dirichlet conditions to be prescribed at boundary points different from mesh vertices, or yet over entire boundary edges or faces, in connection with methods of order greater than one in problem's natural norm, for a wide spectrum of boundary value problems. For example, the application of its principle should avoid order erosion of the $RT_1$ mixed method (cf. \cite{RaviartThomas}) or yet the second order modification of the $BDM_1$ mixed method considered in \cite{Brunner} in the case where fluxes are prescribed all over disjoint smooth curved portions of the boundary. \\
\indent In order to avoid non essential difficulties we confine the study of our technique taking as a model the Poisson equation solved by the classical Lagrange tetrahedron-based methods of degree greater than one. For instance, if quadratic finite elements are employed and we shift prescribed solution boundary values from the true boundary to the mid-points of the boundary edges of the approximating polyhedron, the error of the numerical solution will be of order not greater than $1.5$ in the energy norm (cf. \cite{Ciarlet}),  
instead of the best possible second order. Unfortunately this only happens if the true domain itself is a polyhedron, assuming of course that the solution is sufficiently smooth.\\ 
\indent Since early days finite element users considered method's isoparametric version, with meshes consisting of curved triangles or tetrahedra, as the ideal way to recover optimality in the case of a curved domain (cf. \cite{Zienkiewicz}). However, besides an elaborated description of the mesh, the isoparametric technique inevitably leads to the integration of rational functions to compute the system matrix. In the case of complex non linear problems, this raises the delicate question on what numerical quadrature formula should be used 
to compute element matrices, in order to avoid qualitative losses in the error estimates or ill-posedness of approximate problems. In contrast, in the technique described in \cite{Maugin} and analyzed in \cite{arXiv} for two-dimensional problems, exact numerical integration can be used for the most common non linearities, since we only have to deal with polynomial integrands. Furthermore the element geometry remains the same as in the case of polytopic domains. It is noteworthy that both advantages do not bring about any order erosion in the error estimates that hold for our method, as compared to the equivalent isoparametric version. As a matter of fact the former can be viewed as a small perturbation of the usual Galerkin formulation with conforming Lagrange finite elements based on meshes consisting of triangles or tetrahedra with straight edges. The two-dimensional case was addressed in detail in \cite{arXiv}, in \cite{PAMM} and references therein, in connection with both Lagrange elements and Hermite finite elements with normal-derivative degrees of freedom.  
This work focuses on three-dimensional Lagrange finite element methods, in which the trial functions are discontinuous, in contrast to their two-dimensional counterparts. Likewise the classical conforming Lagrange family is thoroughly studied. Furthermore we consider a situation among many others, in which an isoparametric construction in the strict sense of the term (cf. \cite{Zienkiewicz}) is helpless. More 
precisely we study a second-order method which is nonconforming even for polyhedral domains. \\
\indent An outline of the paper is as follows. Section 2 is devoted to the model problem in a smooth three-dimensional domain 
selected for the presentation of our method. Some pertaining notations are also given therein, followed by preliminary considerations concerning the boundary of this domain as related to the family of meshes to be used in the sequel. In Section 3 we describe our method's main ingredients to solve the model problem, by means of the standard Lagrange family of finite element methods. The underlying approximate problem is posed in Section 4; corresponding stability and well-posedness results are given therein. In Section 5 error estimates are first proved in the energy norm. In the same section error estimates in the $L^2$-norm are also provided, which to the best of author's knowledge are unprecedented in the framework of the class of problems addressed in this article. In Section 6 we illustrate the approximation properties of our method studied in the previous sections, by solving some test-problems with the standard quadratic Lagrange finite element. Application of our technique to a nonconforming Lagrange second-order method having no effective isoparametric counterpart is considered in Section 7. We conclude in Section 8 with some comments on the whole work.

\section{Preliminaries}

Before introducing and studying our method we specify the particular framework in which its application is considered in this work.

\subsection{The model problem}  

Although the method studied in this work extends in a straightforward manner to more complex second-order boundary value problems, symmetric or non symmetric, linear or non linear, in order to simplify the presentation we consider as a model the Poisson equation with Dirichlet boundary conditions in a  
three-dimensional domain $\Omega$ with boundary $\Gamma$ having suitable regularity properties, that is,

\begin{equation}
\label{Poisson}
\left\{
\begin{array}{l}
 - \Delta u = f \mbox{ in } \Omega \\
u = g \mbox{ on } \Gamma,
\end{array}
\right.
\end{equation} 
\noindent where 
$f$ and $g$ are given functions defined in $\Omega$ and on $\Gamma$. \\
Our technique is most effective in connection with methods of order $k > 1$ in the energy norm $\parallel {\bf grad} (\cdot) \parallel_{0}$, in case $u \in H^{k+1}(\Omega)$, where $\parallel \cdot \parallel_{0}$ equals $[ \int_{\Omega} (\cdot)^2 ]^{1/2}$ (that is, the standard norm of $L^2(\Omega)$).  
Accordingly, in order to make sure that $u$ possesses the $H^{k+1}$-regularity property we shall assume that $f \in H^{k-1}(\Omega)$ and $g \in H^{k+1/2}(\Gamma)$ (cf. \cite{Adams}). At this point we observe that, owing to the Sobolev Embedding Theorem \cite{Adams} $g$ is necessarily continuous since $k$ is not less than one. We must further assume that $\Omega$ is sufficiently smooth. For instance, if $k=2$ we assume that $\Gamma$ is at least of the $C^1$-class. Actually, more than this, we make the assumption that, whatever $k$, the principal curvatures of $\Gamma$ (cf. \cite{Cartan}) are uniquely defined almost everywhere. Notice that in doing so we are not necessarily requiring that $\Gamma$ be of the $C^2$-class.
We also note that our regularity assumptions rule out the case where $\Gamma$ is the union of smooth curved portions which do not form a manifold of the $C^1$-class.

\subsection{Meshes and related sets} 

Let us be given a mesh ${\mathcal T}_h$ consisting of straight-edged tetrahedra satisfying the usual compatibility conditions (see e.g. \cite{Ciarlet}). Every element of ${\mathcal T}_h$ is to be viewed as a closed set. Moreover this mesh is assumed to fit $\Omega$ in such a way that all the vertices of the polyhedron $\cup_{T \in {\mathcal T}_h} T$ lie on $\Gamma$. We denote the interior of this union set by $\Omega_h$ and define $\tilde{\Omega}_h := \Omega \cap \Omega_h$, $\Omega^{'}_h := \Omega \cup \Omega_h$. The boundaries of $\Omega_h$ and $\tilde{\Omega}_h$ are respectively denoted by $\Gamma_h$ and $\tilde{\Gamma}_h$ and moreover  
$\Gamma^{'}_h:= \bar{\Omega}_h \cap \Gamma$. ${\mathcal T}_h$ is assumed to belong to a regular family of partitions 
in the sense of \cite{Ciarlet} (cf. Section 3.1), though not necessarily quasi-uniform. 
The boundary of every $\forall T \in {\mathcal T}_h$ is represented by $\partial T$ and its diameter by $h_T$,  while $h :=  \max_{T \in {\mathcal T}_h} h_T$. We make the non essential and yet reasonable assumption that any element in ${\mathcal T}_h$ have at most either one edge or one face contained in $\Gamma_h$. Actually such a condition is commonly fulfilled in practice, for thereby excessively flat tetrahedra are avoided. \\
Let ${\mathcal S}_h$ be the subset of ${\mathcal T}_h$ consisting of tetrahedra having 
one face on $\Gamma_h$ and ${\mathcal R}_h$ be the subset of ${\mathcal T}_h \setminus {\mathcal S}_h$ of tetrahedra having exactly one edge on $\Gamma_h$. We also set ${\mathcal O}_h := {\mathcal S}_h \cup {\mathcal R}_h$. 
Notice that, owing to our initial assumption, no tetrahedron in ${\mathcal T}_h \setminus {\mathcal O}_h$ has a nonempty intersection with $\Gamma_h$. \\
For every $T \in {\mathcal S}_h$ we denote by $O_T$ the vertex of $T$ not belonging to $\Gamma$.

\subsection{Notations} 

Hereafter $\parallel \cdot \parallel_{r,D}$ and $| \cdot |_{r,D}$ represent, respectively, the standard norm and semi-norm of Sobolev space $H^{r}(D)$ (cf. \cite{Adams}), for  
$r \in \Re^{+}$ with $H^0(D)=L^2(D)$, $D$ being a subset of $\overline{\Omega^{'}_h}$. We also denote by $\parallel \cdot \parallel_{m,p,D}$ the usual norm of $W^{m,p}(D)$ for $m \in \natu^{*}$ and 
$p \in [1,\infty] \setminus \{2\}$ with $W^{0,p}(D)=L^p(D)$. Whenever $D$ is $\Omega$ the subscript $,D$ is dropped.\\
 
Finally we introduce the notations $\parallel \cdot \parallel_{0,h}$ and $\parallel \cdot \parallel_{\widetilde{0,h}}$ 
for the standard norms of $L^2(\Omega_h)$ and $L^2(\tilde{\Omega}_h)$, respectively, 
 which will play a key role in the reliability analysis of our method. This is because all our error estimates will be given in the former norm if $\Omega$ is convex and in the latter otherwise. \\
In this respect it is noticeable that for a given mesh and a function $v \in L^2(\Omega)$, $\parallel v \parallel_{0,h}$ (resp. $\parallel v \parallel_{\widetilde{0,h}}$) may equal zero, even if $v$ does not vanish in 
$\Omega \setminus \Omega_h$. However in asymptotic terms this situation is ruled out as far as $u$ is concerned. Indeed the estimates are supposed to hold as $h$ goes to zero, since the family of meshes under consideration is regular (cf. \cite{Ciarlet}, Sect. 3.1). Thus the meshes asymptotically cover the whole $\Omega$. Incidentally this apparently 
indefinite error measure in the case of curved domains is the one used in classical textbooks on the mathematical analysis of the finite element method, such as \cite{Ciarlet} (cf. Section 4.4. p.266 and on) and \cite{StrangFix} (cf. Section 4.4, p.192 and on).

\subsection{Basic assumptions for the formal analysis}
 
\begin{figure}[h]
\label{fig0}
\begin{center}
\includegraphics[scale=0.5]{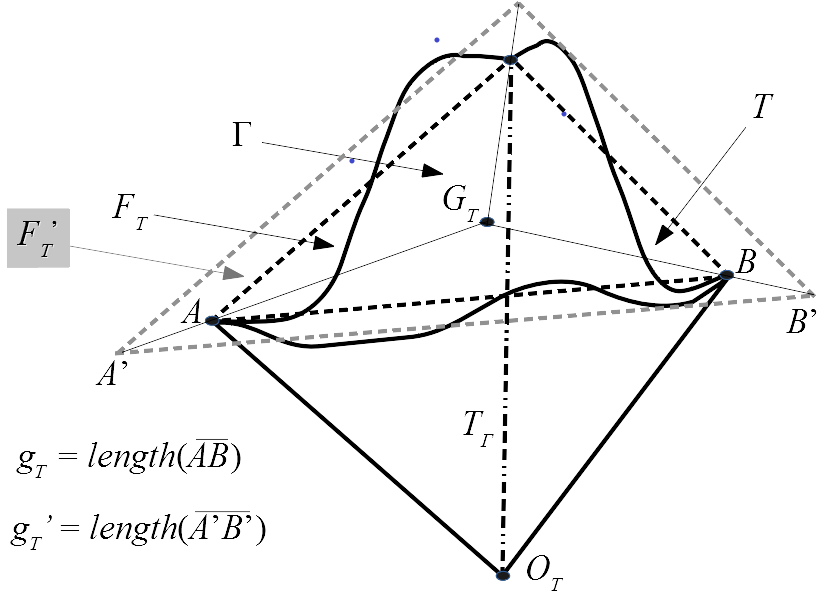}
\end{center}
\par
\caption{Set $T_{\Gamma}$ and triangles $F_T$ and $F^{'}_T$ associated with a tetrahedron $T \in {\mathcal T}_h$ with three vertices on $\Gamma$}  
\end{figure}

Although this is by no means necessary, neither to define our method, nor to implement it, 
henceforth we assume that the meshes in use are fine enough to satisfy some geometric criteria. This assumption is a key sufficient condition for the subsequent reliability results to hold. It also enables the capture of all the nuances of $\Gamma$ by its discrete counterpart $\Gamma_h$, taking advantage of the great flexibility of tetrahedral meshes to fit curvilinear boundaries, even those with sharp variations of shape. \\
Referring to Figure \ref{fig0}, we first associate with every $T \in {\mathcal S}_h$ a closed set $T_{\Gamma}$ delimited by $\Gamma$ and the planes of the faces of $T$ intersecting at $O_T$. Further referring to Figure \ref{fig0}, let $G_T$ be the centroid of $F_T$, $g_T$ be the largest edge of $F_T$ and $F^{'}_T$ be a homothetic transformation of $F_T$ in its plane with center $G_T$ and ratio $\gamma_T = g^{'}_T/g_T$, where $g^{'}_T$ is the maximum edge length of $F^{'}_T$. We take $\gamma_T = 1+ C_H$ where $C_H$ is a small non negative constant independent of $T$ and $h$, though sufficiently large for any point of $F_T^{'}$ to be the orthogonal projection onto the plane of $F_T$ of at most one point $P$ in a simply connected portion of $\Gamma$ 
\footnote{It is not difficult to figure out that $C_H$ can even be proportional to $h_T$.}. \\
We first require the following condition:\\

\noindent \underline{\textit{Assumption}$^{+}$ :} $h$ is small enough for the intersection $P$ with $\Gamma$ belonging to $T_{\Gamma}$ of any segment joining $O_T$ to a point $M \in F_T$ to be uniquely defined $\forall T \in {\mathcal S}_h$.  \rule{2mm}{2mm} \\

In addition to \textit{Assumption}$^{+}$ the following condition is also supposed to be satisfied by the meshes: \\
Let $H_T$  be the closest intersection with $\Gamma$ of the perpendicular to $F_T$ passing through $G_T$. We know that there exists a ball $B(H_T,r_H)$ and a plane $\Pi_H$ swept by the coordinates $x_T,y_T$ of an orthogonal coordinate system $(x_T,y_T,z_T)$ with origin $O_H$, such that a function $f_T(x_T,y_T)$ of the piecewise $C^2$-class uniquely expresses the coordinate $z_T$ of points located on $\Gamma$, as long as they lie in $B(H_T,r_H)$ (cf. \cite{Evans}).\\ 

\noindent \underline{\textit{Assumption}$^{*}$ :} $h$ is small enough for $\Pi_H$ to be taken parallel to $F_T$ and the ball $B(H_T,r_H)$ to contain $F^{'}_T$ $\forall T \in {\mathcal S}_h$  \rule{2mm}{2mm} \\

Some important consequences of both assumptions above are as follows:

\begin{e-proposition}
\label{prop01}
If \textit{Assumption}$^{+}$ and \textit{Assumption}$^{*}$ hold there exists a constant $C^{+}_{\Gamma}$ depending only on $\Gamma$ such that $\forall M \in F_T$ the length of the segment joining $M$ and $P \in T_{\Gamma} \cap \Gamma$ aligned with $O_T$ and $M$ is bounded above by $C^{+}_{\Gamma} h_T^2$ . 
\end{e-proposition}

\prov The proof is based on the fact that, provided $h$ is sufficiently small, the maximum ${\mathcal H}_{max}$ of the euclidean norm in $F^{'}_T$ of the Hessian ${\mathcal H}(f_T)$ of the function $f_T$, is bounded above by an expression depending only on the Gausssian curvature and the mean curvature of $\Gamma$ multiplied by a constant independent of $T$. In the Appendix we give a rigorous justification of this assertion. Taking it for granted, since $f_T$ vanishes at the end-points of every edge $e$ of $F_T$, the first order derivative of $f_T$ in the direction of $e$, say $\partial f_T/\partial e$, must vanish at at least a point $N_e \in e$. Therefore at any point $N \in F^{'}_T$ we have $[\partial f_T/\partial e](N) \leq {\mathcal H}_{max} length(\overline{N_eN}) \leq {\mathcal H}_{max} g^{'}_T$. Since this bound holds for all the three edges of $F_T$ the maximum of the euclidean norm of the gradient of $f_T$ in $F^{'}_T$ denoted by ${\mathcal G}_{max}$ is uniformly bounded above by $2{\mathcal H}_{max} g^{'}_T$, or yet by $2(1+C_H){\mathcal H}_{max} h_T$.\\
Next, since $f_T(A)=0$ where $A$ is a vertex of $F_T$, we note that $\forall N \in F^{'}_T$, $|f_T(N)| \leq length(\overline{AN}) {\mathcal G}_{max} $. Therefore $\forall N \in F^{'}_T$, $|f_T(N)| \leq 
(1+C_H) {\mathcal G}_{max} h_T$. Finally let $\theta_0$ denote the smallest angle between $F_T$ and $\overline{MP}$, which is bounded below independently of $T$ and $h$ for a regular family of meshes. Letting $N$ be the orthogonal projection of $P$ onto the plane of $F_T$ (supposedly a point of $F^{'}_T$), the result follows with $C^{+}_{\Gamma} = 2(1+C_H)^2 {\mathcal H}_{max}[sin(\theta_0)]^{-1}$.  \rule{2mm}{2mm}

\begin{e-proposition}
\label{prop02}
Assume that $\Gamma$ is of the piecewise $C^{k+1}-class$ for $k>1$. Let $D_{x,y}^j v$ be the $j$-th order tensor, whose components are the partial derivatives of order $j$ with respect to $x$ and $y$ of a sufficiently differentiable function $v(x,y)$.   
If \textit{Assumption}$^{*}$ holds, there exists constants $C_{\Gamma}^{j}$ depending only of $\Gamma$ such that $|[D_{x,y}^j f_T](M)| \leq C^j_{\Gamma} h_T^{\max[2-j,0]}$ $\forall M \in F^{'}_T$ for $j=1,2\ldots,k+1$.
\end{e-proposition} 

\prov From the proof of Proposition \ref{prop01} we infer that the result holds true for $j=0$ with $C^0_{\Gamma} = C^{+}_{\Gamma}$ and for $j=1$ with $C^1_{\Gamma}=2(1+C_H){\mathcal H}_{max}$. As for $j=2$ we also saw that the result holds true with 
$C^2_{\Gamma} = {\mathcal H}_{max}$. Finally for $2 < j \leq k+1$ the bound is a simple consequence of the regularity assumptions on 
$\Gamma$. \rule{2mm}{2mm}

\section{Method description} 

\indent First of all we need some additional definitions regarding the set $(\Omega \setminus \Omega_h) \cup (\Omega_h \setminus \Omega)$.\\
With every edge $e$ of $\Gamma_h$ we associate a plane set $\delta_e$ containing $e$, delimited by $\Gamma$ and $e$ itself and set $\delta^{'}_e := \delta_e \cap \bar{\Omega}$. The plane of $\delta_e$ can be arbitrarily chosen about $e$. However for better results it should be close to the bisector of the faces of the pair of elements in ${\mathcal S}_h$ intersecting at $e$, which can eventually be a face shared by both. Such a choice will be assumed throughout this work. 
Although the contrary is perfectly possible, in order to avoid more cumbersome descriptions, $\delta_e$ is supposed not to lie in the plane of a face common to a tetrahedron in ${\mathcal S}_h$ and a tetrahedron in ${\mathcal R}_h$. In Figure 1 we illustrate one out of three such plane sets corresponding to the edges of the faces $F_T$ and $F_{T^{'}}$ contained in $\Gamma_h$ of tetrahedra $T$ and $T^{'}$ belonging to ${\mathcal S}_h$. More precisely we show $\delta_e$ for an edge $e$ common to $F_T$ and $F_{T^{'}}$.\\

For theoretical purposes ${\mathcal T}_h$ is supposed to fulfill a condition analogous to \textit{\textit{Assumption}}$^{+}$, namely, \\

\noindent \underline{\textit{Assumption}$^{++}$ :} $h$ is small enough for the intersection $Q$ with $\Gamma$ belonging to any plane set 
$\delta_e$ of every perpendicular to $e$ through a point of $M \in e$ to be uniquely defined.  \rule{2mm}{2mm} \\

In view of this assumption, akin to Proposition \ref{prop01}, the following result can be established:
\begin{e-proposition}
\label{prop03}
Let $e \subset \Gamma_h$ be an edge of $T \in {\mathcal O}_h$. If \textit{Assumption}$^{++}$ and \textit{Assumption}$^{*}$ hold there exists a constant $C^{++}_{\Gamma}$ depending only on $\Gamma$ such that $\forall M \in e$ the length of the segment joining $M$ and the point $Q$ defined in the former is bounded above by $C^{++}_{\Gamma} h_T^2$. \rule{2mm}{2mm}
\end{e-proposition}

Henceforth we refer to $C_{\Gamma}$ as the maximum between $C^{+}_{\Gamma}$ and $C^{++}_{\Gamma}$. \\

Further, for every $T \in {\mathcal S}_h$, we define a closed set $\Delta_T$ delimited by $\Gamma$, $\partial T$ and the nonempty sets $\delta^{'}_e$ associated with the edges of $F_T$, as illustrated in Figure 1. 
In this manner we can assert that, if $\Omega$ is convex, $\Omega_h$ is a proper subset of $\Omega$ and $\bar{\Omega}$ is the union of the disjoint sets $\Omega_h$ and $\displaystyle \cup_{T \in {\mathcal S}_h} \Delta_T$. Otherwise $\Omega_h \setminus \Omega$ is a nonempty set containing subsets of $T \in {\mathcal S}_h$ whose volume is an $O(h_T^4)$ and subsets of $T \in {\mathcal R}_h$ whose volume is an $O(h_T^5)$, both types of subsets corresponding to non-convex portions of $\Gamma$. Whatever the case, the above configurations are of  merely academic interest and carry no practical meaning, as much as the sets 
$T_{\Delta}:=T \cup \Delta_T$ $\forall T \in {\mathcal S}_h$ or $T_{\Delta}:=T \cup \delta_e$ $\forall T \in {\mathcal R}_h$ and $\tilde{T}:= T \cap \Omega$ $\forall T \in {\mathcal O}_h$.
\begin{figure}[h]
\label{fig1}
\begin{center}
\includegraphics[scale=0.5]{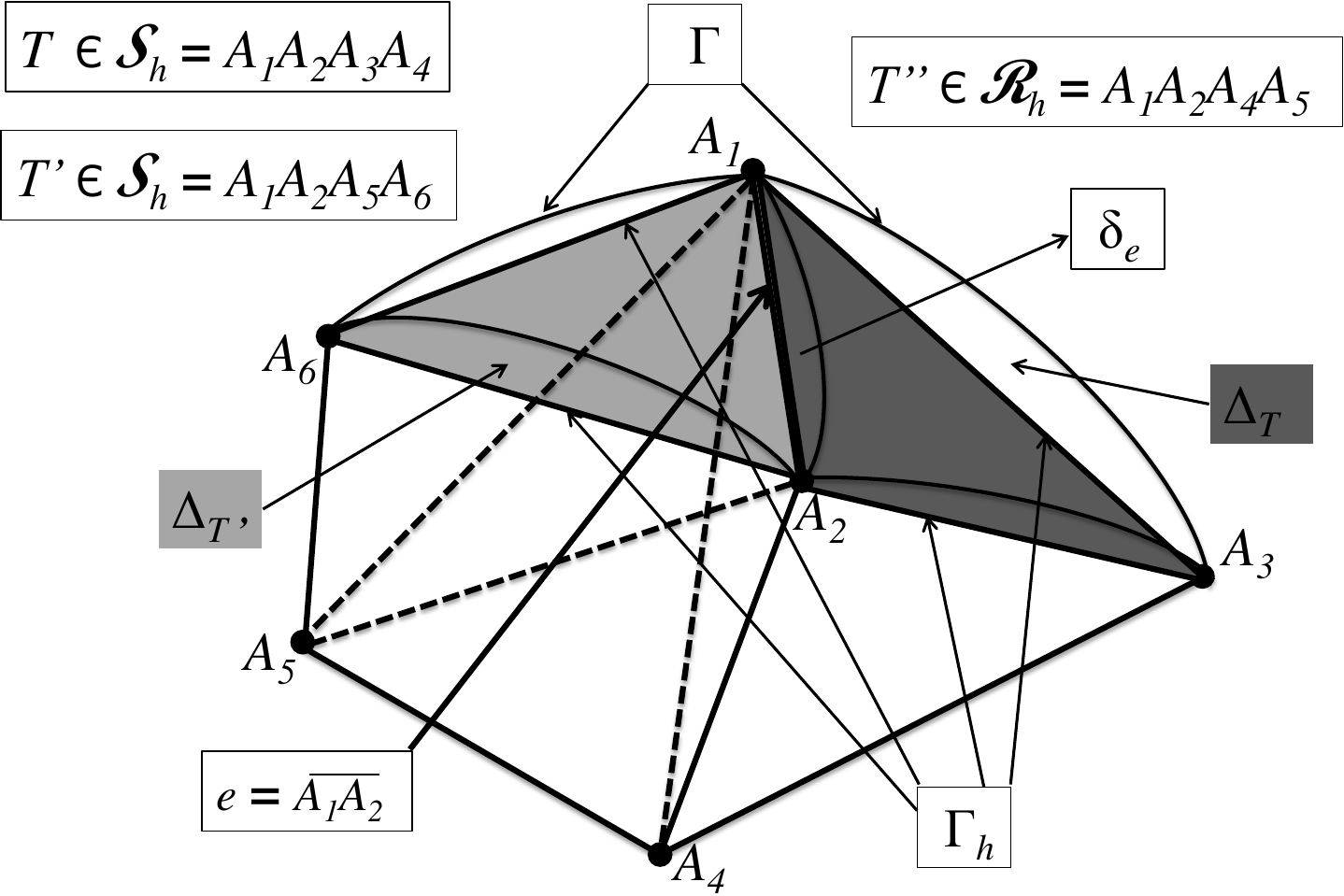}
\end{center}
\par
\caption{Sets $\Delta_T$, $\Delta_{T^{'}}$, $\delta_e$ for tetrahedra $T,T^{'} \! \in \! {\mathcal S}_h$ with a common edge $e$ and a 
tetrahedron $T^{''} \! \in \! {\mathcal R}_h$}
\end{figure} 

\indent Next we introduce a space $V_h$ and a linear manifold $W_h^g$ associated with ${\mathcal T}_h$. With this aim we denote by ${\mathcal P}_m(D)$ the space of polynomials of degree less than or equal to $m$ in a bounded subset $D$ of $\Re^n$.  \\
$V_h$ is the standard Lagrange finite element space consisting of continuous functions $v$ defined in $\Omega_h$ that vanish on $\Gamma_h$, whose restriction to every $T \in {\mathcal T}_h$ belongs to ${\mathcal P}_k(T)$ for $k \geq 2$. For convenience we extend by zero every function $v \in V_h$ to $\Omega \setminus \Omega_h$. We recall that a function in $V_h$ is uniquely defined by its values at the points which are vertices of the partition of each tetrahedron in ${\mathcal T}_h$ into $k^3$ equal tetrahedra (cf. \cite{Zienkiewicz}). Henceforth such points will be referred to as the \textit{Lagrangian nodes} (of order $k$ if necessary). \\
$W_h^g$ in turn is the set of functions defined in $\bar{\Omega}_h$ 
having the properties listed below. \\

\begin{enumerate} 
\item The restriction of $w \in W_h^g$ to every $T \in {\mathcal T}_h$ belongs to ${\mathcal P}_k(T)$;
\item Every $w \in W_h^g$ is single-valued at the vertices of $\Omega_h$ and the inner Lagrangian nodes of the mesh, i.e., at all its Lagrangian nodes of order $k$, but those located on $\Gamma_h$ which are not vertices of $\Omega_h$;
\item A function $w \in W_h^g$ takes the value $g(S)$ at any vertex $S$ of $\Gamma_h$; 
\item $\forall T \in {\mathcal S}_h$, $w(P) = g(P)$ at every $P$ among the $(k-1)(k-2)/2$ nearest intersections with $\Gamma$ 
of the line passing through $O_T$ and the $(k-1)(k-2)/2$ points $M$ not belonging to any edge of $F_T$ among the  
$(k+2)(k+1)/2$ points of $F_T$ that subdivide this face (opposite to $O_T$) into $k^2$ equal triangles (see illustration in Figure 2 for $k=3$); 
\item $\forall T \in {\mathcal O}_h$, $w(Q) = g(Q)$ at every $Q$ among the $k-1$ nearest intersections with $\Gamma$ 
of the line orthogonal to $e$ in the plane set $\delta_e$, passing through the points $M \in e$ different from vertices of $T$, subdividing $e$ into $k$ equal segments, where $e$ generically represents the edge of $T$ contained in $\Gamma_h$ (see illustration in Figure 3 for $k=3$).  
\end{enumerate}
For the subsequent reliability analysis it is convenient to extend to $\bar{\Omega} \setminus \bar{\Omega}_h$ any function $w \in W_h^g$  in such a way that its polynomial expression in $T \in {\mathcal O}_h$ also applies to points in $T_{\Delta} \setminus T$. In doing so, except for the nodes located on $\Gamma$, a function $w \in W_h^g$ is multi-valued in $\delta_e \setminus \bar{\Omega}_h$ if this set is nonempty. In this case the distinct expressions of $w$ therein are those in the tetrahedra belonging to ${\mathcal O}_h$ to which $\delta_e$ is attached.
\begin{figure}[h]
\label{fig2}
\begin{center}
\includegraphics[scale=0.5]{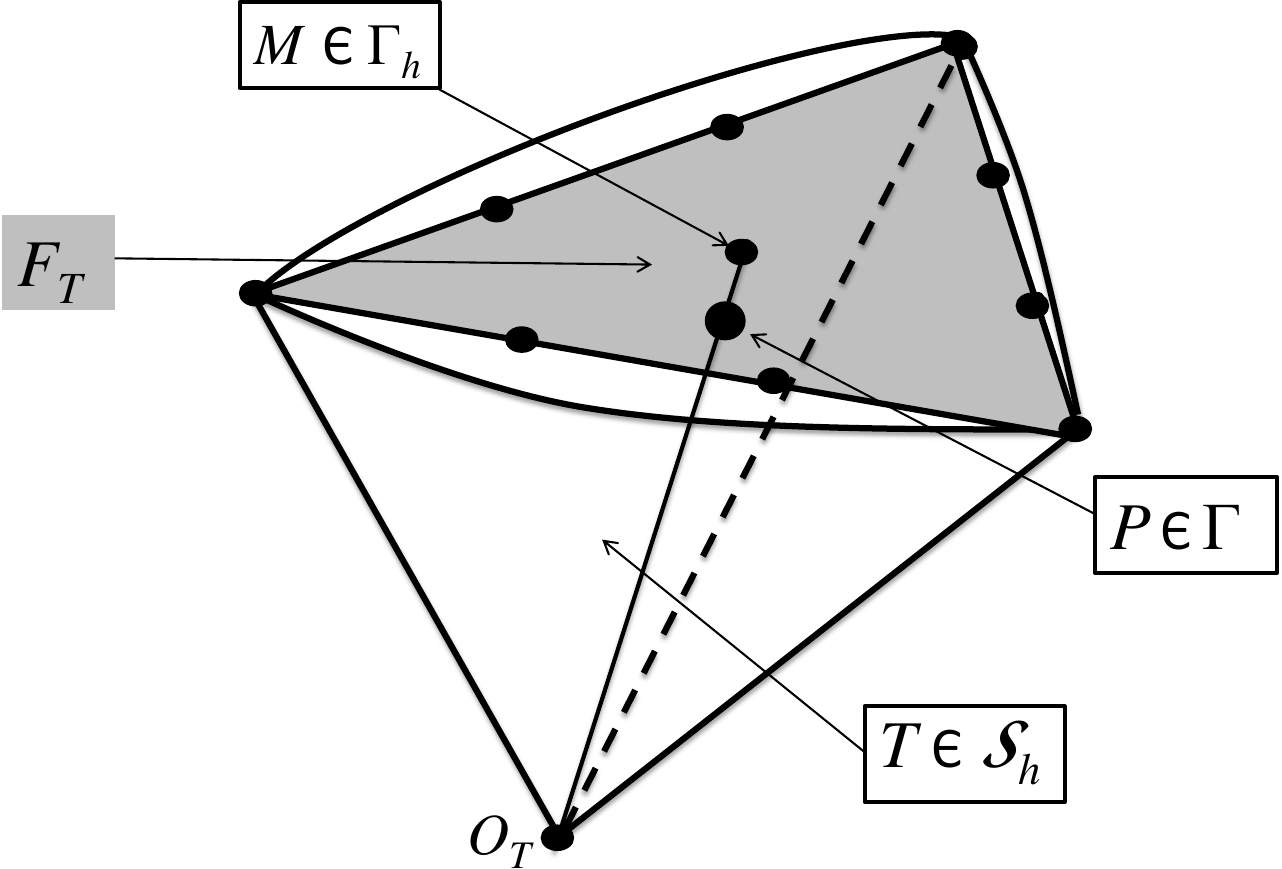}
\end{center}
\par
\caption{Node $P \in \Gamma$ of $W_h^g$ corresponding to the Lagrangian node $M$ in the interior of $F_T \subset \Gamma_h$}
\end{figure}  
\begin{figure}[h]
\label{fig3}
\begin{center}
\includegraphics[scale=0.5]{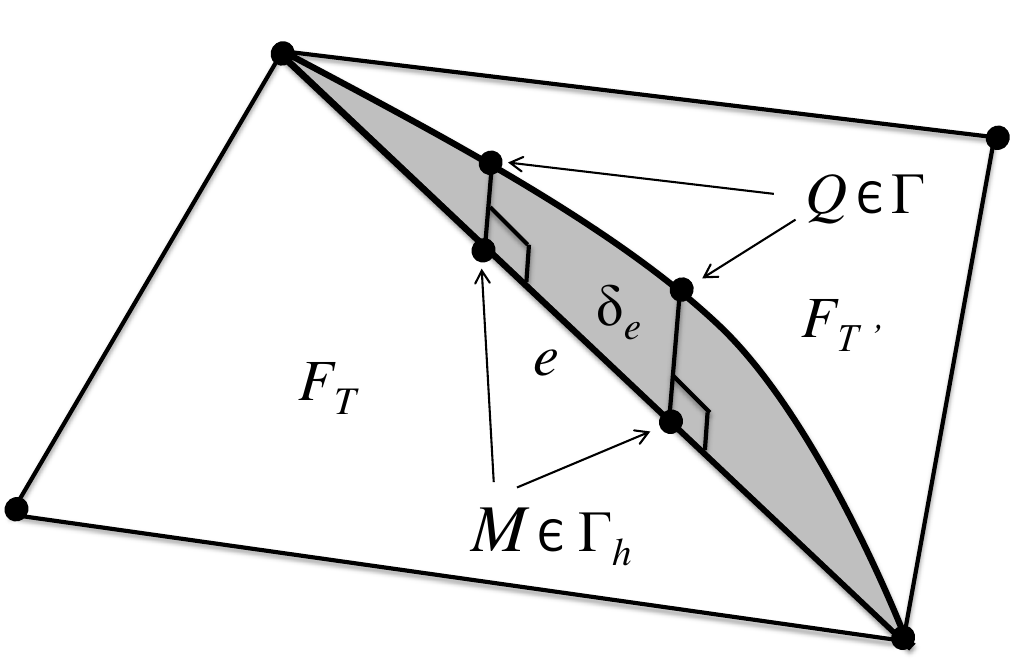}
\end{center}
\par
\caption{Nodes $Q \in \Gamma \cap \overline{\delta_e}$ of $W_h^g$ related to the Lagrangian nodes $M \in e \subset \Gamma_h$}
\end{figure} 

\begin{remark}	
It is important to stress that the sets $\Delta_T$, $T_{\Delta}$ enable the extension of $w \in W^g_h$ to $\bar{\Omega} \setminus \bar{\Omega}_h$, but play no role in the practical implementation of our method (cf. Remark 3 hereafter). \QED
\end{remark}

\begin{remark}
Unless two elements in ${\mathcal S}_h$ sharing an edge $e \subset \Gamma_h$ also have a common face (in which $\delta_e$ is necessarily contained by construction), a function in $w \in W^g_h$ will not be continuous across their faces intersecting at $e$. This is because otherwise the traces of $w$ from both sides of a face $F$ common to two tetrahedra in ${\mathcal S}_h$ and ${\mathcal R}_h$ necessarily coincide only at a number of nodes on $F$ less by $k-1$ the amount of $(k+2)(k+1)/2$ nodes necessary to uniquely define a polynomial of $P_k$ in two variables. Indeed, if $\delta_e$ is not in the plane of $F$, the $k-1$ nodes $Q$ lying in $\delta_e \cap \Gamma$ which are not vertices of $F$ will not belong to $F$. 
Notice that this situation is in contrast to the two-dimensional counterpart of $W_h^g$, which is a subspace of $H^1(\Omega)$. Notice however that in three-dimensional space $w \in W^g_h$ is necessarily continuous across all faces common to two tetrahedra in the mesh having no edge on $\Gamma_h$  \QED 
\end{remark}

\begin{remark}
The construction of the nodes associated with $W_h^g$ located on $\Gamma$ advocated in items 4. and 5. is not mandatory. Notice that it differs 
from the intuitive construction of such nodes lying on normals to faces of $\Gamma_h$ commonly used in the isoparametric technique. The main advantage of this proposal is the determination by linearity of the coordinates of the boundary nodes $P$ in the case of item 4. Nonetheless the choice of boundary nodes ensuring our method's optimality is absolutely very wide.  \QED
\end{remark}

The fact that $W_h^g$ is a nonempty set is a trivial consequence of the three following lemmata:

\begin{lemma}
\label{TDelta}
Provided $h$ satisfies \textit{Assumption}$^{*}$, \textit{Assumption}$^{+}$ and \textit{Assumption}$^{++}$ there exist two mesh-independent constants ${\mathcal C}_{\infty}$ and ${\mathcal C}_J$ depending only on $\Gamma$ and the shape regularity of ${\mathcal T}_h$ (cf. \cite{BrennerScott}, Ch.4, Sect. 4) such that  
$\forall w \in {\mathcal P}_k(T_{\Delta})$ and $\forall T \in {\mathcal O}_h$ it holds:  
\begin{equation}
\label{LinftyTDelta}
\parallel  w \parallel_{0,\infty,T_{\Delta}} \leq {\mathcal C}_{\infty} \parallel  w \parallel_{0,\infty,\tilde{T}} \mbox{ and}
\end{equation}
\begin{equation}
\label{L2TDelta}
\parallel  w \parallel_{0,\infty,T_{\Delta}} \leq {\mathcal C}_J h_T^{-3/2} \parallel  w \parallel_{0,\tilde{T}}.
\end{equation}
\end{lemma}

\prov 
First we denote the dimension of ${\mathcal P}_k(D)$ for any bounded open   
set $D$ of $\Re^3$ by $n_k$ with $n_k = (k+3)(k+2)(k+1)/6$. \\
Let $0 < \lambda \leq 1$ be the largest possible value for the homothetic transformations $T_{\lambda}$ and $T^{'}_{\lambda}$ of $T \in {\mathcal O}_h$ centered at a vertex of $T$ not lying on $\Gamma$ and with ratios $\lambda$ and $\lambda^{-1}$, to be contained in $\tilde{T}$ and contain $T_{\Delta}$, respectively. 
Now set $\kappa :=1 - \sigma_{\mathcal T} C_{\Gamma} h_0$ and $\kappa^{'}:= 1 + \sigma_{\mathcal T} C_{\Gamma} h_0$ as two numbers depending only on $\Gamma$, where $h_0$ is the largest value of $h$ such that \textit{Assumption}$^{+}$, \textit{Assumption}$^{*}$ and \textit{Assumption}$^{++}$ hold and $\kappa$ is not less than a certain number in the interval $(0,1]$, say $1/2$, and $\sigma_{\mathcal T}$ is a shape-regularity parameter of the family of meshes in use satisfying for every ${\mathcal T}_h$, $\sigma_{\mathcal T} \geq  
\max_{T \in {\mathcal T}_h} h_T/\eta_T$, $\eta_T$ being the minimum height of $T$. From Propositions \ref{prop01} and \ref{prop03}   
together with Thales' Proportionality Theorem, it is rather easy to infer that $\kappa$ and $\kappa^{'}$ are such that the maximum diameters of tetrahedra $T_{\lambda}$ and $T^{'}_{\lambda}$ lie in the intervals $[\kappa h_T, h_T]$ and $[h_T,\kappa^{'} h_T]$, respectively. 
Since both $T_{\lambda}$ and $T^{'}_{\lambda}$ are similar 
to $T$, these tetrahedra have the same shape regularity property as any other element in ${\mathcal T}_h$, provided the maximum diameter of each member of the family of partitions in use is adjusted to take into account the thus modified maximum diameters.\\
Let $T \in {\mathcal O}_h$. Denoting by $\varphi_i$ the canonical basis function associated with the $i$-th Lagrangian node $M_i$ of $T$ extended to $T^{'}_{\lambda}$, for every $w \in {\mathcal P}_{k}(T_{\Delta})$ we can write,
\begin{equation}
\label{auxiliary6mono}
\parallel w \parallel_{0,\infty,T_{\Delta}} \leq \displaystyle \sum_{i=1}^{n_k} |w(M_i)| 
\max_{{\bf x} \in T_{\lambda}^{'}} |\varphi_i({\bf x})|.
\end{equation}
Next we resort to the master tetrahedron $\hat{T}$ with vertices $(0,0,0),(1,0,0),(0,1,0),(0,0,1)$ 
in a reference frame. ${\mathcal F}_T$ being the affine mapping from $T$ onto $\hat{T}$ let $\hat{\varphi}_i$ and $\hat{w}$ be the transformations of $\varphi_i$ and $w$ under ${\mathcal F}_T$. Let also $\hat{T}_{\lambda}$ and $\hat{T}_{\lambda}^{'}$ be the transformations of $T_{\lambda}$ and $T_{\lambda}^{'}$ under ${\mathcal F}_T$. Then it holds: 
\begin{equation}
\label{auxiliary6bis}
\parallel w \parallel_{0,\infty,T_{\Delta}} \leq \hat{C}_1 \displaystyle \sum_{i=1}^{n_k} |w(M_i)| \; \forall w \in {\mathcal P}_{k}(T_{\Delta}),
\end{equation}
where 
$$\hat{C}_1 =\displaystyle \max_{1 \leq i \leq n_k} \left[\max_{\hat{\bf x} \in \hat{T}_{\lambda}^{'}} |\hat{\varphi }_i(\hat{\bf x})| \right].$$
Owing to the equivalence of norms in the $n_k$-dimensional space ${\mathcal P}_k(\hat{T}_{\lambda})$, there exists a constant $\hat{C}_2$ depending only on $\hat{T}$, $\lambda$ and $k$ such that $\forall w \in {\mathcal P}_k(T_{\Delta})$, 
\begin{equation}
\label{auxiliary6ter}
\displaystyle \sum_{i=1}^{n_k} |w(M_i)| = \displaystyle \sum_{i=1}^{n_k} |\hat{w}({\mathcal F}_T(M_i))| \leq \hat{C}_2 \parallel \hat{w} \parallel_{0,\infty,\hat{T}_{\lambda}}.  
\end{equation} 
Combining \eqref{auxiliary6bis} and \eqref{auxiliary6ter} it easily follows that \eqref{LinftyTDelta} holds with 
${\mathcal C}_{\infty} = \hat{C}_1 \hat{C}_2$. \\
Finally we note that $volume({T}_{\lambda}) \leq \hat{\mathcal C}_{J}^2 h_T^{-3} volume(T_\lambda)$ with a constant $\hat{\mathcal C}_{J}$ independent of $T$. Then using again the equivalence of norms in ${\mathcal P}_k(\hat{T}_{\lambda})$ we infer the existence of another constants $\hat{C}_{\lambda}$ independent of $T$ for which it holds, 
\begin{equation}
\label{inverseTlambda}
\parallel \hat{w} \parallel_{0,\infty,\hat{T}_{\lambda}} \leq \hat{C}_{\lambda} 
\parallel \hat{w} \parallel_{0,\hat{T}_{\lambda}} \leq \hat{C}_{\lambda} \hat{\mathcal C}_{J} h_T^{-3/2} \parallel w \parallel_{0,T_{\lambda}} \; \forall w \in {\mathcal P}_k(T_{\Delta}).
\end{equation} 
Since $T_{\lambda} \subset \tilde{T}$, combining \eqref{auxiliary6mono}, \eqref{auxiliary6bis}, \eqref{auxiliary6ter}, \eqref{L2TDelta} must hold with ${\mathcal C}_J= \hat{C}_{\lambda} \hat{\mathcal C}_{J} {\mathcal C}_{\infty}$ independently of $\tilde{T}$ and $T_{\Delta}$. \QED \\ 

\begin{lemma}
\label{lemma2}
Provided $h$ is small enough $\forall T \in {\mathcal S}_h$, given a set of $m_k$ real values $b_{i}$, $i=1,\ldots,m_k$ with $m_k=k(k+2)(k+1)/6$, 
there exists a unique function $w_T \in {\mathcal P}_k(T)$ that takes the value of $g$ at the three vertices $S$ 
of $T$ located on $\Gamma$, at the $(k-1)(k-2)/2$ points $P$ of $\Gamma$ defined in accordance with item 4. and at the $3(k-1)$ points $Q$ of $\Gamma$ 
defined in accordance with item 5. of the above definition of $W_h^g$, and takes the value $b_i$ respectively at the $m_k$ Lagrangian nodes of $T$ not located on $\Gamma_h$. 
\end{lemma}   
 
\prov Let us first extend the vector $\vec{b}:=[b_1,b_2,\ldots,b_{m_k}]$ of $\Re^{m_k}$ into a vector of $\Re^{n_k}$ still denoted by $\vec{b}$, with $n_k:=m_k+(k+2)(k+1)/2$, by adding $n_k-m_k$ components $b_i$ which are the values of $g$ at the $(k+2)(k+1)/2$ nodes ($P$ or $Q$) of $T_{\Delta_T}$ located on $\Gamma$. If the  latter nodes were replaced by the corresponding $M \in \Gamma_h \cap T$, it is clear that the result would hold true, according to the well-known properties of Lagrange finite elements. The vector $\vec{a}$ of coefficients $a_i$ for $i=1,2,\ldots,n_k$ of the canonical basis functions $\varphi_i$ of ${\mathcal P}_k(T)$ for $1 \leq i \leq n_k$ would be precisely $b_i$ for $1 \leq i \leq n_k$. Still denoting by $M_i$ the Lagrangian nodes of $T$, $i=1,2,\ldots,n_k$, this means that the matrix $K$ whose entries are $k_{ij} := \varphi_j(M_i)$ is the identity matrix. Let $\tilde{M}_i=M_i$ if $M_i \notin \Gamma \setminus \Gamma_h$ and $\tilde{M}_i$ be the node of the type $P$ or $Q$ associated with $M_i$ otherwise. The Lemma will be proved if the $n_k \times n_k$ linear system $\tilde{K} \vec{a} = \vec{b}$ is uniquely solvable, where $\tilde{K}$ is the matrix with entries $\tilde{k}_{ij}:=\varphi_j(\tilde{M}_{i})$. Clearly we have $\tilde{K} = K + E_K$, where the entries of $E_K$ are $e_{ij}:= \varphi_j(\tilde{M}_{i}) - \varphi_j(M_{i})$. At this point we recall the constant $C^{+}_{\Gamma}$ depending only on $\Gamma$ specified in Proposition \ref{prop01} such that the length of the segment $\overline{M_i\tilde{M}_i}$ is bounded above by $C^{+}_{\Gamma} h_T^2$. From Rolle's Theorem it follows that $\forall \;i,j$, 
$|e_{ij}| \leq C^{+}_{\Gamma} h_T ^2 \parallel {\bf grad}\;\varphi_j \parallel_{0,\infty,T_{\Delta}}$.  \\
Thanks to fact that $\varphi_j \in {\mathcal P}_k(T_{\Delta})$ and to \eqref{LinftyTDelta}, $\parallel {\bf grad}\;\varphi_j \parallel_{0,\infty,T_{\Delta}} \leq {\mathcal C}_{\infty} \parallel {\bf grad}\;\varphi_j \parallel_{0,\infty,T}$. Moreover from standard arguments we know that the latter norm in turn is bounded above by a mesh-independent constant times $h_T^{-1}$. In short we have $|e_{ij}| \leq C_E h_T$ $\forall \; i,j$, where $C_E$ is a 
mesh-independent constant. Hence the matrix $\tilde{K}$ equals the identity matrix plus an $O(h_T)$ matrix $E_K$. Therefore $\tilde{K}$ is an invertible matrix, as long as $h$ is sufficiently small.    
 \QED \\

\begin{lemma}
\label{lemma3}
Provided $h$ is small enough $\forall T \in {\mathcal R}_h$, given a set of $p_k$ real values $b_{i}$, $i=1,\ldots,p_k$ with $p_k=(k+1)(k+2)(k+3)/6-(k+1)$, 
there exists a unique function $w_T \in {\mathcal P}_k(T)$ that takes the value of $g$ at the two end-points $S$ of the edge $e$ of $T$ located on $\Gamma$ and at the $k-1$ points $Q$ of $\Gamma$ defined in accordance with item 5. of the above definition of $W_h^g$, and takes the 
value $b_i$ respectively at the $p_k$ Lagrangian nodes of $T$ not located on $\Gamma_h$.
\end{lemma}
   
\prov Thanks to Proposition \ref{prop03} this lemma can be proved on the grounds of the same arguments already exploited in the proof of Lemma \ref{lemma2}. \QED \\

\section{The approximate problem}

Lemmata \ref{lemma2} and \ref{lemma3} entitle us to  
set a problem associated with the space $V_h$ and the manifold $W_h^g$, whose solution is an approximation of the solution $u$ of (\ref{Poisson}).  \\  
Before posing this problem we introduce the broken gradient operator ${\bf grad}_h$ for any function $w$ defined in $\Omega_h$ which is continuously differentiable in every $T \in {\mathcal T}_h$, given by $[{\bf grad}_h w]_{|T} \equiv {\bf grad} \; w_{|T}$ $\forall T \in {\mathcal T}_h$. \\   
Extending $f$ by zero in $\Omega_h \setminus \Omega$ and still denoting the resulting function by $f$, we wish to solve,
\begin{equation}
\label{Poisson-h}
\left\{
\begin{array}{l}
\mbox{Find } u_h \in W_h^g \mbox{ such that } a_h(u_h,v) = L_h(v) \; \forall v \in V_h \\
\\
\mbox{where } \\
\\
a_h(w,v) := \int_{\Omega_h} {\bf grad}_h\; w \cdot {\bf grad}\; v; \\
\\
L_h(v) : = \int_{\Omega_h} fv.
\end{array}
\right.
\end{equation}     

To begin with we establish the stability of \eqref{Poisson-h}. 
   
\begin{e-proposition}
\label{propo2}
Let $W_h^0$ be the space of functions corresponding to the manifold $W_h^g$ for $g \equiv 0$. Then provided $h$ is sufficiently small there exists a constant $\alpha > 0$ independent of 
$h$ such that,
\begin{equation}
\label{inf-sup}
\forall w \in W_h^0 \neq 0, \displaystyle \sup_{v \in V_h \setminus \{ 0 \}} \frac{a_h(w,v)}{\parallel {\bf grad}_h w \parallel_{0,h} \parallel {\bf grad} \; v \parallel_{0,h}} 
\geq \alpha. 
\end{equation}  
\end{e-proposition}

\prov Given $w \in W_h^0$ let $v \in V_h$ coincide with $w$ at all Lagrangian nodes of elements $T \in {\mathcal T}_h$ 
not belonging to ${\mathcal O}_h$. As for an element $T \in {\mathcal O}_h$ we set $v=w$ at the Lagrangian nodes not belonging to $\Gamma_h$, while $v=0$ at the Lagrangian nodes located on $\Gamma_h$. 
The fact that on the faces common to two elements $T^{-}$ and $T^{+}$ in ${\mathcal T}_h$, both $v_{|T^{-}}$ and $v_{|T^{+}}$ are polynomials of degree less than or equal to $k$ 
in two variables coinciding at the exact number of Lagrangian nodes required to uniquely define such a function, implies that $v$ is continuous in $\Omega_h$. Moreover for the same reason $v$ vanishes all over $\Gamma_h$. 
 \\
\noindent Let us denote by ${\mathcal M}_T$ the set of Lagrangian nodes of $T \in {\mathcal O}_h$ 
that belong to $\Gamma_h$, and are different from vertices.  
Clearly enough we have 
\begin{equation}
\label{ahwv} 
a_h(w,v) = \displaystyle \sum_{T \in {\mathcal T}_h} \int_T |{\bf grad} \; w |^2 - \displaystyle \sum_{T \in {\mathcal O}_h} \int_T {\bf grad} \; w \cdot {\bf grad} \; r_T(w), 
\end{equation}
\noindent where $r_T(w) = \sum_{M \in {\mathcal M}_T} w(M) \varphi_M$, $\varphi_M$ being the canonical basis function of the space 
${\mathcal P}_k(T)$ associated with the Lagrangian node $M$. \\
\noindent Now from standard results it holds $\parallel {\bf grad} \; \varphi_M \parallel_{0,T} \leq C_{\varphi} h_T^{1/2}$ where $C_{\varphi}$ is a mesh independent constant. Moreover, since $w(P)=0$ (resp. $w(Q)=0$), where 
$P$ (resp. $Q$) generically represent the point of $\Gamma$ corresponding to $M \in \Gamma_h$ in accordance with the definition of $W_h^0$, a simple Taylor expansion about $P$ (resp. $Q$) allows us to conclude that $|w(M)| \leq l \parallel {\bf grad} \; w \parallel_{0,\infty,T_{\Delta}}$, where $l = length(\overline{PM})$ (resp. $length(\overline{QM})$), or yet    
$|w(M)| \leq C_{\Gamma} h_T^2 \parallel {\bf grad} \; w \parallel_{0,\infty,T_{\Delta}}$. 
On the other hand from \eqref{L2TDelta} it holds $\parallel {\bf grad} \; w \parallel_{0,\infty,T_{\Delta}} \leq {\mathcal C}_J h_T^{-3/2} 
\parallel {\bf grad} \; w \parallel_{0,T}$. Plugging all those estimates into (\ref{ahwv}) we obtain:
\begin{equation}
\label{ahwvbound} 
a_h(w,v) \geq \int_{\Omega_h} |{\bf grad}_h w |^2 - C_{\varphi} {\mathcal C}_J C_{\Gamma} h
\displaystyle \sum_{T \in {\mathcal O}_h} card({\mathcal M}_T)  \parallel {\bf grad} \; w \parallel_{0,T}^2. 
\end{equation}         
Since $card({\mathcal M}_T) \leq (k+4)(k-1)/2 \; \forall T$, setting $c:= C_{\varphi} {\mathcal C}_J C_{\Gamma}[(k+4)(k-1)/2]$, it holds,
\begin{equation}
\label{ahwvbelow} 
a_h(w,v) \geq (1 - c h) \parallel {\bf grad}_h w \parallel_{0,h}^2.
\end{equation}
Now using arguments in all similar to those employed above, we easily conclude that 
\begin{equation}
\label{normbound}
\parallel {\bf grad} \; v \parallel_{0,h} \leq \parallel {\bf grad}_h w  \parallel_{0,h} + \parallel {\bf grad} \; v - {\bf grad}_h w \parallel_{0,h} \leq (1+ ch)  \parallel {\bf grad}_h w \parallel_{0,h}.
\end{equation}
Combining (\ref{ahwvbelow}) and (\ref{normbound}), provided $h \leq (2c)^{-1}$ we establish (\ref{inf-sup}) with $\alpha = 1/3$. \QED \\

Now let $u^H \in H^1(\Omega)$ be the solution of the Laplace equation $\Delta u^H = 0$ in $\Omega$ fulfilling $u^H=g$ on $\Gamma$. We may assume that 
$u^H \in H^{k+1}(\Omega)$ with $k>1$, as a trivial consequence of suitable assumptions on $g$ and $\Omega$. 
Thus we can define the interpolate $u^H_h$ of $u^H$ in $W^g_h$. Moreover the simple 
application of standard error estimates for the interpolating function (cf. \cite{BrennerScott}, Ch. 4, Sect. 4) ensures the existence of a mesh-independent constant $C$ 
such that 
\begin{equation}
\label{interperror}
\parallel {\bf grad} \; u^H - {\bf grad}_h u^H_h \parallel_{\widetilde{0,h}} \leq C h^k  | u^H |_{k+1}. 
\end{equation}
We next prove the well-posedness of \eqref{Poisson-h}. With this aim we let $u_h^0 \in W_h^0$ satisfy
\begin{equation}
\label{uh0}
a_h(u^0_h,v) = L_h^0(v) \; \forall v \in V_h \; \mbox{ where } L_h^0(v):=L_h(v)-a_h(u^H_h,v). 
\end{equation}
\begin{e-proposition}
\label{propo3}
Provided $h$ is sufficiently small, problem (\ref{Poisson-h}) has a unique solution.  
\end{e-proposition}

\prov First we note that $L^0_h$ is a continuous linear form on $V_h$, and $a_h$ is a continuous bilinear form on $W^0_h \times V_h$, the spaces 
$V_h$ and $W^0_h$ being equipped with the norms $\parallel {\bf grad}(\cdot) \parallel_{0,h}$ and $\parallel {\bf grad}_h(\cdot) \parallel_{0,h}$, 
respectively. Thus the facts that (\ref{inf-sup}) holds and $dim(V_h) = dim(W_h^0)$ imply the existence and uniqueness of $u^0_h$ according to the theory of 
non-coercive approximate linear variational problems (cf. \cite{Babuska}, \cite{Brezzi} and \cite{COAM}). 
Therefore $u_h:=u^0_h+u^H_h$ is a solution to (\ref{Poisson-h}), and its uniqueness is a direct consequence of (\ref{inf-sup}). \QED

\section{Error estimates}
We next proceed to error estimations for problem (\ref{Poisson-h}). Throughout this section 
we assume that $h$ is small enough to satisfy \textit{Assumption}$^{*}$, \textit{Assumption}$^{+}$ and \textit{Assumption}$^{++}$ and in any case $h<1$. We further assume that $\Gamma$ is at least of the piecewise $C^{k}$-class and require the minimum regularity $f \in H^{k-1}(\Omega)$ and $g \in H^{k+1/2}(\Gamma)$ for $k>1$, so that the solution $u$ of 
(\ref{Poisson}) belongs to $H^{k+1}(\Omega)$.

\subsection{Preliminaries} 
Error estimates in energy norm will be proved by comparing the solution of (\ref{uh0}) with $u^0$, where $u^0 \in H^1_0(\Omega)$ is the unique solution of the equation $-\Delta u^0 = f$ in $\Omega$. Clearly enough $u^0+u^H$ is the solution of (\ref{Poisson}) and hence $u^0$ fulfills: 
\begin{equation}
\label{u0}
a(u^0,v) = L^0(v) \; \forall v \in H^1_0(\Omega), \; \mbox{ where } L^0(v):=L(v)-a(u^H,v),
\end{equation}
with
\begin{equation}
\label{aF3D}
a(w,v) := \int_{\Omega} {\bf grad} \; w  \cdot {\bf grad}\; v \mbox{ and } L(v) : = \int_{\Omega} fv.
\end{equation} 
Henceforth we denote by $D^j w$ the $j$-th order tensor whose components are the $j$-th order partial derivatives with respect to the space variables of a function $w$ in the strong or the weak sense. Alternatively we may also write $H(w)$ instead of $D^2 w$ and ${\bf grad}\;w$ instead of $D^1 w$. \\
Many results in the sequel rely on classical inverse inequalities applying to polynomials defined in $T$ 
(see e. g. \cite{Verfuerth}) and their extensions to neighboring sets.  
Besides \eqref{L2TDelta} we shall use the following one: \\
There exists a constant ${\mathcal C}_I$ depending only on 
$k$ and the shape regularity of ${\mathcal T}_h$ (cf. \cite{BrennerScott}, Ch.4, Sect 4) such that for $1 \leq j \leq k$ it holds:  
\begin{equation}
\label{inverse}
\parallel D^j w \parallel_{0,T} \leq {\mathcal C}_I h_T^{-1} \parallel D^{j-1} w \parallel_{0,T} \; 
\forall w \in {\mathcal P}_k(T) \mbox{ and } \forall T \in {\mathcal T}_h.\\
\end{equation}

Before going into the main results we give some useful additional definitions:
\begin{itemize}
\item
$\partial (\cdot)/\partial n_T$ is the normal derivative on $\partial T$ directed outwards $T \in {\mathcal T}_h$;
;
\item
$\Gamma_T = T \cap \Gamma$ for  $T \in {\mathcal O}_h$;
\item
$\partial (\cdot)/\partial \bar{n}_T$ is the normal derivative $\partial (\cdot)/\partial n$ restricted to $\Gamma_T$ if $area(\Gamma_T)>0$;   
\item 
${\mathcal F}_h$ is the set of faces of elements in ${\mathcal O}_h$ that are not contained in $\Gamma_h$; 
\item 
$\Delta_h :=\Omega \setminus \bar{\Omega}_h$.  
\end{itemize}
\hspace{5mm} It is noteworthy that if $\Omega$ is convex the closure of $\Delta_h$ equals $\cup_{T \in {\mathcal S}_h} \Delta_T$.\\

For $T \in {\mathcal S}_h$ we further introduce the following sets and notations:
\begin{itemize}
\item
$\tilde{\Delta}_T$ is the closure of $\mathring{\Delta}_T \cap \Omega$;
\item 
$\tilde{\partial} T:= (\partial T_{\Delta} \cap \Gamma) \cup \Gamma_T$ 
($\partial T_{\Delta}$ is the boundary of $T_{\Delta}$);
\item
$\bar{\partial} T = \tilde{\partial} T \cup [\cup_{e \subset \Gamma_h \cap T} \; \delta_e^{'}]$;
\item
The normal derivative on $\bar{\partial} T \setminus \Gamma_T$ directed outwards 
$T_{\Delta} \cap \Omega$ is also denoted by $\partial (\cdot)/\partial \bar{n}_T$.  \\
\end{itemize}  

We also need the following technical lemmata.

\begin{lemma}
\label{wh}
Let $r=1/2+\epsilon$ for a certain $\epsilon$ in $(0,1)$ and $w \in H^{k+1+r}(\Omega^{'}_h)$ such that 
$w_{|\Gamma} \equiv 0$. Let $T^{'}$ be a closed set fulfilling $\tilde{T} \subseteq T^{'} \subseteq T_{\Delta}$. Given $w_h \in W_h^0$ assume that $w_h$ is extended to $\bar{\Delta}_h$ as prescribed 
in Section 3. Then there exist constants ${\mathcal C}_j$ independent of $T$ and $h$ such that for $j=1,2,\ldots,k$ it holds,
\begin{equation}
\label{Djwh}
 \parallel \! D^j(w_h-w) \! \parallel_{0,\infty,T^{'}} \leq 
{\mathcal C}_j h_T^{-j-1/2} [\parallel\! {\bf grad}(w_h-w)\! \parallel_{0,\tilde{T}} + h_T^k | w |_{k+1,\tilde{T}} + h_T^{k+r}\! \parallel w \parallel_{k+1+r,T^{'}}].
\end{equation}
\end{lemma}
      
\prov 
First of all, $I_h(w)$ being the $W^0_h$ interpolate of $w$ in $\Omega_h$   
we write $w_h-w = (w_h-I_h(w))+(I_h(w)-w)$, $I_h(w)$ being extended to $\bar{\Delta}_h$ in the same way as $w_h$. Then using \eqref{L2TDelta}we can write,
\begin{equation}
\label{wh1}
 \parallel D^j(w_h-w) \parallel_{0,\infty,T^{'}} \leq {\mathcal C}_J h_T^{-3/2} 
 \parallel D^j(w_h-I_h(w)) \parallel_{0,\tilde{T}} + \parallel D^j(I_h(w)-w) \parallel_{0,\infty,T^{'}}.
\end{equation}
Using the affine transformation ${\mathcal F}_T$ like in the proof of Lemma \ref{TDelta} and setting 
$\hat{T}^{'}:= {\mathcal F}_T(T^{'})$ we observe that  
$H^{k+1+r}(\hat{T}^{'})$ is continuously embedded in $W^{k,\infty}(\hat{T}^{'})$ 
(as much as $H^{k+1+r}(D)$ is in $W^{k,\infty}(D)$ for all open subset $D$ of $\Omega{'}$ cf. \cite{Adams}). Hence applying classical estimates for the interpolation error in fractional Sobolev norms (cf. \cite{Arcangeli}) we obtain for suitable constants $C_j$ independent of $T_{\Delta}$:
\begin{equation}
\label{wh2}
 \parallel D^j(I_h(w)-w) \parallel_{0,\infty,T^{'}} \leq C_j h_T^{k-j+\epsilon} \parallel w \parallel_{k+1+r,T^{'}} \; \mbox{ for }  j=1,2,\ldots,k.
\end{equation}
On the other hand using \eqref{inverse} we easily come up with,
\begin{equation}
\label{wh3}
\parallel D^j(w_h-I_h(w)) \parallel_{0,\tilde{T}} \leq 
[{\mathcal C}_{I} h_T]^{-j+1} [\parallel {\bf grad}(w_h-w)) \parallel_{0,\tilde{T}} + \parallel {\bf grad}(w-I_h(w)) \parallel_{0,\tilde{T}}].
\end{equation}
Moreover by standard approximation results (cf. \cite{BrennerScott}) there exists a mesh-independent constant ${\mathcal C}_L$ such that  
\begin{equation}
\label{wh4}
\parallel D^j(w-I_h(w)) \parallel_{0,\tilde{T}} \leq {\mathcal C}_L h_T^{k+1-j} | w |_{k+1,\tilde{T}} \mbox{ for } 1 \leq j \leq k.
\end{equation}
The combination of \eqref{wh1}, \eqref{wh2}, \eqref{wh3} and \eqref{wh4} with $j=1$ 
immediately yields \eqref{Djwh}.  
\QED \\

\begin{lemma}
\label{deltae}
Let $e$ be an edge of $\Gamma_h$ and also an edge of the face $F_T$ of $T \in {\mathcal S}_h$ contained in 
$\Gamma_h$ and $M_e$ be the mid-point of $e$. Denoting by   
$\vec{n}(P)$ the unit outer normal vector to $\Gamma$ at $P \in \Gamma$, assume that the plane set $\delta_e$ lies in the plane of $e$ and a point $P_e \in \Gamma$ in a perpendicular to $e$ through $M_e$ such that the inner product of $\vec{n}(P_e)$ and the unit vector in the direction of $\overrightarrow{M_eP_e}$ is bounded above by $C_{\delta} h_T$, $C_{\delta}$ being a constant independent of $T$.  
Recalling the notation $\partial v/\partial \bar{n}_T$ for the outer normal derivative on $\delta_e$ with respect to $\Delta_T$ of a function $v \in H^2(\Omega)$, there exists a constant $C_{\theta}$ independent of $T$ such that,
\begin{equation}
\label{deltatrace}
\left\{
\begin{array}{l}
|\Theta(v)| \leq C_{\theta} h \parallel v \parallel_{2} \; \forall v \in H^2(\Omega) \cap H^1_0(\Omega), \\
\mbox{where } \\
\Theta(v) :=\displaystyle \left[ \sum_{T \in {\mathcal S}_h} \sum_{e \subset F_T} 
\left\| \frac{\partial v}{\partial \bar{n}_T} \right\|_{0,\delta_e^{'}}^2 \right]^{1/2}.
\end{array}
\right.
\end{equation}
\end{lemma}

\prov
For a set $\delta^{'}_e$ whose interior is nonempty, let $\gamma_e := \Gamma \cap \delta_e^{'}$ and 
$\vec{n}_e$ be the unit normal vector on $\delta_e$ directed outwards $\Delta_T$. We first introduce an invertible mapping ${\mathcal F}_e$ from $\delta_e^{'}$ onto a (not necessarily connected) plane set $\delta_{e}^{\Gamma}$ contained in $\tilde{\partial} T$ and containing $\gamma_e$, whose Jacobian $J_{e}$ in $\delta_e^{'}$ is uniformly bounded above and below by two strictly positive constants independent of $T$. For convenience assume that the transformation of $\gamma_e$ under ${\mathcal F}_e$ is this set itself. Now $\forall M \in \delta_e^{'}$ let $P \in \delta_{e}^{\Gamma}$ be ${\mathcal F}_e(M)$. \\
\begin{figure}[h]
\label{fig5.2}
\begin{center}
\includegraphics[scale=0.5]{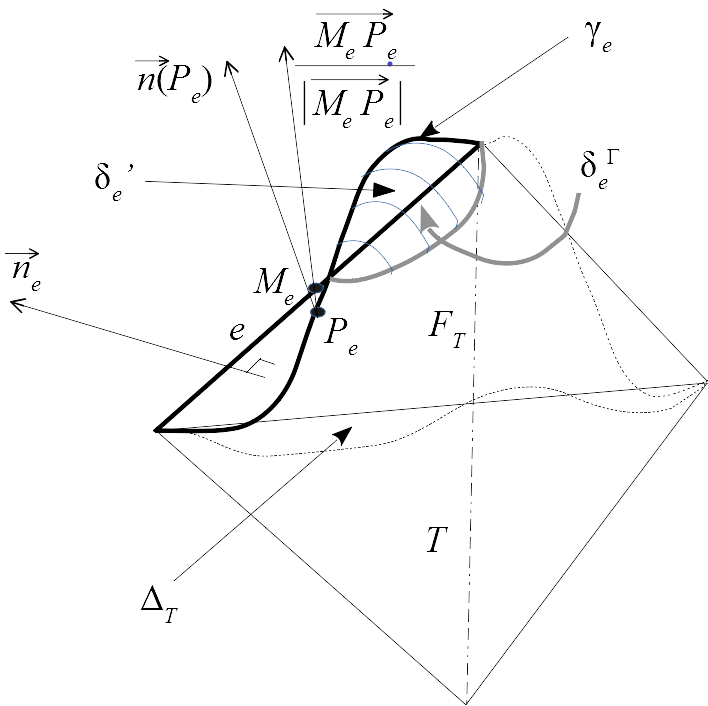}
\end{center}
\par
\caption{Sets $\delta^{'}_e$ surrounded by thick black lines and $\delta^{\Gamma}_e$ depicted in grey and pertaining data}  
\end{figure}
Referring to Figure \ref{fig5.2}, 
for every pair $(M,P) \in \delta_e^{'} \times \delta_{e}^{\Gamma}$ it is possible to construct a unique path $\eta$ leading from $P$ to $M$ entirely contained in $\tilde{\Delta}_T$ with a curvilinear abscissa $t$ such that $\vec{\tau}(t) \cdot \vec{n}_e \geq \beta >0$, $\beta$ being independent of $h$, where $\vec{\tau}$ is the unit tangent vector along $\eta$ oriented from $P$ to $M$. The paths $\eta$ are assumed to be arranged in such a manner that the Cartesian coordinates of $\delta_e$'s plane together with $t$ form a system of curvilinear coordinates in a subset of $\tilde{\Delta}_T$ containing all such paths, whose Jacobian as related to the spatial Cartesian coordinate system is bounded above and below by constants independent of $\tilde{\Delta}_T$. Moreover for every differentiable function $\omega$ a.e. in $T_{\Delta}$ it holds, 
$$\omega(M) = \omega(P) + \displaystyle \int_P^M {\bf grad}\; \omega \cdot \vec{\tau}(t) dt.$$
Therefore by the Schwarz inequality we have,
\begin{equation}
\label{Gammae}
\int_{\delta_e} w^2 \leq 2 \displaystyle \left[ \max_{P \in \delta_{e}^{\Gamma}} [J_{e}]^{-1}(P) \int_{\delta_{e}^{\Gamma}} \omega^2 + l_{e} \int_{\tilde{\Delta}_T} |{\bf grad} \; \omega |^2 \right],
\end{equation}
where $l_{e}$ is proportional to the characteristic height of $\delta_e$, i.e. $l_e$ equals a mesh-independent constant multiplied by $h_T^2$.\\
Let us apply the upper bound \eqref{Gammae} to the function $\omega= {\bf grad} \; v \cdot \vec{n}_e$. Since ${\bf grad} \; v \cdot \vec{s}=0$ on $\Gamma$, for every vector $\vec{s}$ tangent to $\Gamma$, it is clear that $\omega(P)= \vec{n}(P) \cdot \vec{n}_e [\partial v/ \partial n](P) $. However owing to the construction of $\delta_e$ and to the fact that $|\vec{n}(P_e)-\vec{n}(P)|$ is bounded by another mesh-independent constant $C_{\gamma}$ times $h_T$ for every $P \in \delta_{e}^{\Gamma}$, by a straightforward calculation we can assert that $|\vec{n}_e \cdot \vec{n}(P)|$ is bounded above by $(C_{\delta}+C_{\gamma})h_T$ for every $P \in \delta_{e}^{\Gamma}$. Plugging this result into \eqref{Gammae}, and taking into account the uniform boundedness of $[J_{e}]^{-1}$, we come up with a mesh-independent constant $C^{'}_{\theta}$ such that,
\begin{equation}
\label{deltaebound}     
\displaystyle \left\| \frac{\partial v}{\partial \bar{n}_T} \right\|_{0,\delta_e^{'}}^2 \leq \displaystyle 
C^{'}_{\theta} h_T^2 \displaystyle \left[ \left\| \frac{\partial v}{\partial n} \right\|_{0,\tilde{\partial} T}^2 + \parallel H(v) \parallel_{0,\tilde{\Delta}_T}^2 \right] \; \forall e \subset F_T.
\end{equation}
Summing up over $e \subset F_T$ and over $T \in {\mathcal S}_h$ we further obtain, 
\begin{equation}
\label{deltaebound1}
\displaystyle \sum_{T \in {\mathcal S}_h} \sum_{e \subset F_T} 
\left\| \frac{\partial v}{\partial n_T} \right\|_{0,\delta_e^{'}}^2 \leq 
3 C^{'}_{\theta} h^2 \displaystyle \left[ \left\| \frac{\partial v}{\partial n} \right\|_{0,\Gamma}^2 + |v |_{2}^2 \right].
\end{equation}
On the other hand by the Trace Theorem there exists a constant $C_t$ depending only on $\Omega$ such that  
\begin{equation}
\label{normaltrace}
\displaystyle \left\| \frac{\partial v}{\partial n} \right\|_{0,\Gamma} \leq C_t \parallel v \parallel_{2}.
\end{equation}
Plugging \eqref{normaltrace} into \eqref{deltaebound1} the result follows. \QED \\

\begin{lemma}
\label{traceomegaT}
Let $T \in {\mathcal S}_h$ and $\sigma_T$ be a portion of the boundary $\partial T_{\Delta}$ of $T_{\Delta}$ with a strictly positive area. Then for every $\omega_T \in H^1(T)$ there is a mesh-independent constant $C_{\sigma}$ such that,  
\begin{equation}
\label{intSTomegaT}
\int_{\sigma_T} |\omega_T| \leq C_{\sigma} h_T^{1/2} \displaystyle \left[ \int_{T_{\Delta}} ( \omega_T^2 
+ h_T^2 |{\bf grad} \; \omega_T|^2) \right]^{1/2}
\end{equation}
\end{lemma}

\prov
Let us resort again to the master tetrahedron $\hat{T}$ and denote the transformation of $\omega_T$ under the affine invertible mapping ${\mathcal F}_T$ from $T$ onto $\hat{T}$ by $\hat{\omega}$. Clearly enough there exists a constant $\bar{C}_{\sigma}$ independent of $T$ such that,
\begin{equation}
\label{estim2b3}
\int_{\sigma_T} |\omega_T| \leq \int_{\partial T_{\Delta}} |\omega_T| 
\leq \bar{C}_{\sigma} h_T^2 \int_{\partial \hat{T}_{\Delta}} |\hat{\omega}|. 
\end{equation}  
where $\partial \hat{T}_{\Delta}$ is the boundary of the transformation $\hat{T}_{\Delta}$ of $T_{\Delta}$ under ${\mathcal F}_T$. Next we apply the Trace Theorem to $\hat{T}_{\Delta}$. Thanks to the fact that $\Gamma$ is smooth and $h$ is sufficiently small, there exists a constant $\hat{C}_{\sigma}$ independent of $T$ such that,
\begin{equation}
\label{estim3b3} 
\displaystyle \int_{\partial \hat{T}_{\Delta}} \hat{\omega} \leq \hat{C}_{\sigma} \displaystyle \left\{ \int_{\hat{T}_{\Delta}} [ \hat{\omega}^2 + |\widehat{\bf grad}\; \hat{\omega}|^2 ]  \right\}^{1/2},   
\end{equation}  
where $\widehat{\bf grad}$ is the gradient operator for functions defined in $\hat{T}_{\Delta}$.\\
Moving back to $T_{\Delta}$ and noting that $volume(\hat{T}_{\Delta})/volume(T_{\Delta}) \leq 
C_{\Delta} h_T^{-3}$ for a certain constant $C_{\Delta}$ independent of $T$, using \eqref{estim2b3} and \eqref{estim3b3} we obtain \eqref{intSTomegaT} for a suitable $C_{\sigma}$. \QED \\

\begin{lemma}
\label{IF}
Let $T \in {\mathcal O}_h$ and $F$ be a face of $T$ belonging to ${\mathcal F}_h$. Let also $\tilde{F} = F \cap \Omega$ and $I_F$ be the 
operator $I_F: H^2(\Omega^{'}_h) + W_h^0 \rightarrow {\mathcal P}_k(\tilde{F})$ such that 
$[I_F(w)](N) = w(N)$ for all the $(k+2)(k+1)/2$ Lagrangian nodes $N$ of order $k$ on $F$, $\forall w \in 
H^2(\Omega^{'}_h) + W_h^0$. Let also $\tilde{I}_F : H^2(\Omega^{'}_h) + W_h^0 \rightarrow {\mathcal P}_k(\tilde{F})$ be the operator such that $\forall w \in H^2(\Omega^{'}_h) + W_h^0$ $[\tilde{I}_F(w)](N)=w(N)$ at all the $k(k+1)/2 + 2$ Lagrangian nodes $N$ of $T$ on $F$ of order $k$ not located in the interior of its edge $e \subset \Gamma_h$, and $[\tilde{I}_F(w)](M_i) = w(P_i)$ for $i=1,\ldots,k-1$, where the $M_i$s are the Lagrangian nodes of order $k$ of $T$ on $F$ located in the interior of $e$, $P_i \in \Gamma$ being the nodal point associated with $M_i$ in accordance with the definition of $W^0_h$.
Then if $w$ belongs to $H^{k+1+r}(\Omega_h^{'})$ ($r=1/2+\epsilon$) and 
$w_h \in W_h^0$ there exists a mesh-independent constant ${\mathcal C}_F$ such that,
\begin{equation}
\label{jump}
\parallel [I_F-\tilde{I}_F](w_h-w) \parallel_{0,\tilde{F}} \leq {\mathcal C}_F h_T^{3/2} 
[\parallel {\bf grad}(w_h-w) \parallel_{0,\tilde{T}} + h_T^{k} | w |_{k+1,\tilde{T}} + h_T^{k+r} \parallel w \parallel_{k+1+r,T_{\Delta}}].
\end{equation}     
\end{lemma}

\prov
First of all since $area(\tilde{F}) \leq area(F) \leq h_T^2/2$, we wave
\begin{equation}
\label{IF1}
\parallel [I_F - \tilde{I}_F](w_h-w) \parallel_{0,\tilde{F}} \leq \displaystyle \frac{h_T}{\sqrt{2}} 
\parallel [I_F-\tilde{I}_F](w_h-w) \parallel_{0,\infty,\tilde{F}}.
\end{equation}
Moreover recalling the canonical basis function $\varphi_i$ of ${\mathcal P}_k(T)$ associated with $M_i$ 
for $i=1,\ldots,k-1$, the construction of the operators $I_F$ and $\tilde{I}_F$ allows us to write,
\begin{equation}
\label{IF2} 
\parallel [I_F - \tilde{I}_F](w_h-w) \parallel_{0,\infty,\tilde{F}} 
\leq \displaystyle \sum_{i=1}^{k-1} |(w_h-w)(M_i)-(w_h-w)(P_i) | \parallel \varphi_i \parallel_{0,\infty,\tilde{F}} 
\end{equation} 
Since the distance between $M_i$ and $P_i$ is bounded above by $C_{\Gamma} h_T^2$, \eqref{IF2} easily yields,
\begin{equation}
\label{IF3} 
\parallel [I_F- \tilde{I}_F](w_h-w) \parallel_{0,\infty,\tilde{F}} 
\leq \tilde{C}_{\varphi} h_T^2 \parallel {\bf grad}(w_h-w) \parallel_{0,\infty,T_{\Delta}}, 
\end{equation} 
where $\tilde{C}_{\varphi}$ is a constant depending only on $k$ and $\Gamma$.\\
Now we combine \eqref{IF3} and \eqref{IF1} and recall \eqref{Djwh} with $j=1$, to establish \eqref{jump} 
with ${\mathcal C}_F = {\mathcal C}_1 \tilde{C}_{\varphi} \sqrt{2}/2$.
\QED \\

We had pointed out that the position of the plane set $\delta_e$ about $e$ is irrelevant for our method to work. It is relevant however for proving $L^2$-error estimates. With this aim 
henceforth we take for granted that the plane sets $\delta_e$ are chosen as prescribed in Lemma \ref{deltae}.
Notice that such a condition on the position of $\delta_e$ just means that $\forall e$ it is  
roughly upright with respect to $\Gamma$, which is a rather intuitive construction. 
\subsection{The case of convex domains}
At an initial stage we assume that $\Omega$ is convex.  
\begin{theorem}
\label{convconvex}
There exists a constant ${\mathcal C}(f,g)$ depending only on $f$ and $g$ such that the solution $u_h$ of (\ref{Poisson-h}) satisfies :
\begin{equation}
\label{errestconvex}
\parallel {\bf grad}_h(u - u_h) \parallel_{0,h} \leq {\mathcal C}(f,g) h^k. 
\end{equation}
\end{theorem}

\prov Owing to the convexity of $\Omega$ we have $V_h \subset H^1_0(\Omega)$. Hence the variational residual $a(u^0,v)-L^0(v)$ vanishes for every $v \in V_h$. 
On the other hand $a_h(u^0,v)=-\int_{\Omega_h} v \Delta u^0 = -\int_{\Omega} v \Delta u^0 = a(u^0,v)=L^0(v)$ if $v \in V_h$. It follows that the variational residual $a_h(u^0,v)-L^0_h(v)$ equals $L^0(v)-L^0_h(v) \; \forall v \in V_h$. 
According to \cite{COAM} we thus have: 
\begin{equation}
\label{errorbound0}
\parallel {\bf grad}_h( u^0 - u^0_h) \parallel_{0,h} \leq \displaystyle \frac{1}{\alpha} \left[ \displaystyle \inf_{w \in W_h^0} 
\parallel {\bf grad}_h (u^0 - w) \parallel_{0,h} + \displaystyle \sup_{v \in V_h \setminus \{0\}} \frac{|L^0(v)-L^0_h(v)|}{\parallel {\bf grad}\; v \parallel_{0,h}}. \right]
\end{equation}
We know that $\displaystyle \inf_{w \in W_h^0} \parallel {\bf grad}_h(u^0- w) \parallel_{0,h} \leq C h^k | u^0 |_{k+1}$. \\
Moreover $ |L^0(v)-L^0_h(v)| = |a_h(u^H_h-u^H,v)| 
\leq C h^k | u^H |_{k+1} \parallel {\bf grad} \; v \parallel_{0,h}$, according to (\ref{interperror}). \\
Summarizing, it holds
\begin{equation}
\label{errestimate0}
\parallel {\bf grad}_h(u^0 - u^0_h) \parallel_{0,h} \leq \displaystyle \frac{C}{\alpha} h^k [ | u^0 |_{k+1} + | u^H |_{k+1} ].
\end{equation}
Finally (\ref{errestconvex}) easily derives from \eqref{errestimate0} and the triangle inequality with ${\mathcal C}(f,g) = C_F C/\alpha \parallel f \parallel_{k-1} + 
C_G (1+ C/\alpha) \parallel g \parallel_{k+1/2,\Gamma}$, where $C_G$ and $C_F$ are 
constants such that $| u^0 |_{k+1} \leq C_F \parallel f \parallel_{k-1}$ and $| u^H |_{k+1} \leq C_G \parallel g \parallel_{k+1/2,\Gamma}$. \QED 
\begin{remark}
It is noticeable that the continuity of functions in $W_h^g$ is nowhere required in the above error analysis. Indeed in the generalization given in \cite{COAM} of classical error bounds such as Strang's inequalities, only the residual $a_h(u,v) - L_h(v)$ needs to be evaluated for $v \in V_h$. Thanks to the continuity of functions in $V_h$ this residual trivially vanishes. Incidentally this explains why it is not reasonable to replace $V_h$ by $W_h^0$, as one might be tempted to in order to define a symmetric approximate problem. \QED \\
\end{remark} 

If we assume that the solution of \eqref{Poisson} is a little more regular, it is possible to establish for 
problem \eqref{Poisson-h} an $O(h^{k+1})$-error estimate in the norm of $L^2(\Omega_h)$, based on \eqref{errestconvex} and a classical duality argument. We observe that for the two-dimensional analog of \eqref{Poisson-h} the validity of such an estimate was proven in \cite{arXiv}, at the price of a rather laborious analysis. In the three-dimensional case the study becomes even more complex since our method is nonconforming, in contrast to its two-dimensional version. That is why for the sake of brevity we next prove an $L^2$-error estimate in the particular case where $g \equiv 0$.
\begin{theorem}
\label{theorem1bis}
Let $k>1$ and $\Omega$ be convex. Assume that $\Omega$ is of the piecewise $C^{k+1}$-class and the solution $u$ of (\ref{Poisson}) for $g \equiv0$ belongs to $H^{k+1+r}(\Omega)$, for $r=1/2+\epsilon$ where $\epsilon>0$ can be arbitrarily small. Then the solution $u_h$ of (\ref{Poisson-h}) satisfies for a suitable constant ${\mathcal C}_0$ independent of $h$ and $u$:
\begin{equation}
\label{L2estconvex}
\parallel u - u_h \parallel_{0,h} \leq {\mathcal C}_0 h^{k+1} \parallel u \parallel_{k+1+r}.
\end{equation}
\end{theorem}

\prov 
Let $\bar{u}_h$ be the function defined in $\Omega$ such that $\bar{u}_h = u_h - u$ in $\Omega_h$,  satisfying the following condition in $\Omega \setminus \Omega_h$. Recalling the definition of the set $\Delta_T$ for $T \in {\mathcal S}_h$ illustrated in Figure 1, and the fact that the expression of $u_h$ in $T$ extends to $\Delta_T$, $\bar{u}_h$ is also given by $u_h-u$ in $\Delta_T$ $\forall T \in {\mathcal S}_h$. Notice that this also defines $\bar{u}_h$ on both sides of the plane sets $\delta_e$ depicted in Figure 1, and hence $\bar{u}_h$ is defined everywhere in $\bar{\Omega}$. \\
Now let $v \in H^1_0(\Omega)$ be the solution of 
\begin{equation}
\label{adjoint}
-  \Delta v = \bar{u}_h \; \in \Omega.
\end{equation}
Since $\Omega$ is smooth and $\bar{u}_h \in L^2(\Omega)$ we know that $v \in H^2(\Omega)$, and moreover there 
exists a constant $C_{\Omega}$ depending only on $\Omega$ such that,
\begin{equation}
\label{adjoint1}
\parallel v \parallel_{2} \leq C_{\Omega} \parallel \bar{u}_h \parallel_{0}. 
\end{equation}
Presumably $v$ does not vanish identically in $\Omega$, otherwise the analysis that follow is useless. Therefore we can write,
\begin{equation}
\label{L2est1} 
\parallel \bar{u}_h \parallel_{0} \leq C_{\Omega} \displaystyle 
\frac{- \int_{\Omega} \bar{u}_h \Delta v }{\parallel v \parallel_{2}}.
\end{equation}
Now using integration by parts we obtain,
\begin{equation}
\label{L2est2} 
\parallel \bar{u}_h \parallel_{0} \leq C_{\Omega} \displaystyle \frac{a_h(\bar{u}_h, v)
+a_{\Delta_h}(\bar{u}_h,v)-a_{\partial h}(\bar{u}_h,v)}{\parallel v \parallel_{2}}, 
\end{equation}
where $\forall w \in W_h^0 + H^1(\Omega)$ and $\forall v \in H^1(\Omega)$,
\begin{equation}
\label{aDeltah}
a_{\Delta_h}(w,v) := \int_{\Delta_h} {\bf grad}_h w \cdot {\bf grad}\; v \; 
\end{equation}
and $\forall w \in W_h^0 + H^1(\Omega)$ and $\forall v \in H^2(\Omega)$,
\begin{equation}
\label{apartialh}
\begin{array}{l}
a_{\partial h}(w,v) := \displaystyle \sum_{T \in {\mathcal R}_h} 
\int_{\partial T \setminus \tilde{F}_T} w \frac{\partial v}{\partial n_T} 
+ \displaystyle \sum_{T \in {\mathcal S}_h} 
\left[\int_{\bar{\partial} T}  w \frac{\partial v}{\partial \bar{n}_T} 
 + \int_{\partial T \setminus F_T}  w \frac{\partial v}{\partial n_T} \right], \\
\end{array}
\end{equation}
$\tilde{F}_T$ being the union of the two faces of a tetrahedron $T$ in ${\mathcal R}_h$ that do not contain its 
edge $e \subset \Gamma_h$. \\
Then we observe that $a_{\partial h}(w,v) =  - c_h(w,v) - d_h(w,v)-b_{1h}(w,v)$ where 
\begin{equation}
\label{ch}
c_h(w,v) := - \displaystyle \sum_{T \in {\mathcal R}_h} 
\int_{\partial T \setminus \tilde{F}_T} w \frac{\partial v}{\partial n_T} - \displaystyle \sum_{T \in {\mathcal S}_h} 
\int_{\partial T \setminus F_T}  w \frac{\partial v}{\partial n_T} 
\; \forall w \in W_h^0 + H^1(\Omega) \; \forall v \in H^2(\Omega), 
\end{equation}
\begin{equation}
\label{deltah}
d_h(w,v) := - \displaystyle \sum_{T \in {\mathcal S}_h}\sum_{e \subset F_T} \int_{\delta_e}  w \frac{\partial v}{\partial \bar{n}_T} \; \forall w \in W_h^0 + H^1(\Omega) \; \forall v \in H^2(\Omega),
\end{equation} 
\begin{equation}
\label{b1h}
b_{1h}(w,v) := - \displaystyle \oint_{\Gamma} w \frac{\partial v}{\partial n} \; \forall w \in W_h^0 + H^1(\Omega) \; \forall v \in H^2(\Omega).
\end{equation} 
Let $\Pi_h(v)$ be the continuous piecewise linear interpolate of $v$ in $\Omega$ at the vertices of the mesh. 
Setting $v_h=\Pi_h(v)$ in $\Omega_h$ and $v_h=0$ in $\Delta_h$ we have $v_h \in V_h$. 
Therefore it holds $a(u,v_h)=a_h(u,v_h)=L(v_h)=L_h(v_h)=a_h(u_h,v_h)$.
Now we split $a_{\Delta_h}(\bar{u}_h,v)$ into the sum $a_{\Delta_h}(\bar{u}_h,v-\Pi_h(v))
+a_{\Delta_h}(\bar{u}_h,\Pi_h(v))$ and apply First Green's identity in $\Delta_T$ for $T \in {\mathcal S}_h$. 
Since $\Pi_h(v)_{|\Gamma_h} \equiv 0$, we come up with 
$a_{\Delta_h}(\bar{u}_h,\Pi_h(v)) = b_{2h}(\bar{u}_h,\Pi_h(v))+b_{3h}(\bar{u}_h,\Pi_h(v))$, where
\begin{equation}
\label{b2h}
b_{2h}(w,z):= - \displaystyle \sum_{T \in {\mathcal S}_h}  \int_{\Delta_T} z \Delta w  \mbox{ for } w \in W_h^0 + H^2(\Omega) \mbox{ and } z \in H^{1}(\Omega),  
\end{equation}
and
\begin{equation}
\label{b3h}
b_{3h}(w,z) := \displaystyle \sum_{T \in {\mathcal S}_h} \int_{\bar{\partial} T} \frac{\partial w}{\partial \bar{n}_T} z \mbox{ for } w \in W_h^0 + H^2(\Omega) \mbox{ and } z \in H^1(\Omega).
\end{equation}
Further setting $e_h(v) = v -\Pi_h(v)$ together with,
\begin{equation}
\label{b4h} 
b_{4h}(w,z) := a_{\Delta_h}(w,z) \mbox{ for } w \in W_h^0 + H^1(\Omega) \mbox{ and } z \in H^1(\Omega),  
\end{equation}  
it follows that,
\begin{equation}
\label{L2est3} 
\begin{array}{l}
\parallel \bar{u}_h \parallel_{0} \leq C_{\Omega} \displaystyle \left[ \frac{a_h(\bar{u}_h,e_h(v))}
{\parallel v \parallel_{2}}\right. \\
\\
\left. \displaystyle \frac{b_{1h}(\bar{u}_h,v)+b_{2h}(\bar{u}_h,\Pi_h(v))+b_{3h}(\bar{u}_h,\Pi_h(v))+
b_{4h}(\bar{u}_h,e_h(v))+c_h(\bar{u}_h,v)+d_h(\bar{u}_h,v)}{\parallel v \parallel_{2}} \right]. 
\end{array}
\end{equation} 
From classical results, for a mesh-independent constant $C_{V}$ it holds
\begin{equation}
\label{interP1}
\parallel {\bf grad} \; e_h(v) \parallel_{0,h} \leq \parallel {\bf grad} \; e_h(v) \parallel_{0} \leq C_{V} h 
| v |_{2}.
\end{equation} 
Therefore, combining \eqref{errestconvex}, \eqref{interP1} and \eqref{L2est3}, and setting  
$\tilde{\mathcal C}_0=C_{\Omega} C_{V} {\mathcal C}(f,0)$ we have,
\begin{equation}
\label{L2est4}
\begin{array}{l} 
\parallel \bar{u}_h \parallel_{0} \leq \tilde{\mathcal C}_0 h^{k+1} | u |_{k+1} + C_{\Omega} \times\\
\\
\displaystyle \frac{b_{1h}(\bar{u}_h,v)+b_{2h}(\bar{u}_h,\Pi_h(v))+b_{3h}(\bar{u}_h,\Pi_h(v))+
b_{4h}(\bar{u}_h,e_h(v))+c_h(\bar{u}_h,v)+d_h(\bar{u}_h,v)}{\parallel v \parallel_{2}}.
\end{array}
\end{equation}
In \eqref{L2est4} the functionals $c_h$ and $d_h$ account for the discontinuity of functions in $W_h^0$, which does not occur in the two-dimensional case. The functionals $b_{ih}$ for $i=1,2,3,4$ in turn are residual-like terms that also appear in the two-dimensional case, though in a significantly simpler form. That is why it is necessary to carry out a thorough study of the latter too, which we do next.\\  
\indent As for $b_{1h}$ we first note that according to \eqref{normaltrace} we have,
\begin{equation}
\label{L2estimate1}
b_{1h}(\bar{u}_h,v) \leq C_t \parallel \bar{u}_h \parallel_{0,\Gamma} \parallel v \parallel_{2}.
\end{equation}
Moreover we have,   
\begin{equation}
\label{L2estimate2}
 \parallel \bar{u}_h \parallel_{0,\Gamma} = 
\displaystyle \left[\sum_{T \in {\mathcal S}_h} \int_{\tilde{\partial} T} \bar{u}_h^2 \right]^{1/2}.
\end{equation}
Let $\tilde{F}_T$ be the domain on the plane of $F_T$ delimited by the projections of the curves $\gamma_e \subset \delta_e$ for the three edges $e$ of $F_T$ (notice that $\tilde{F}_T$ is nothing but the projection of $\tilde{\partial} T$ onto the plane of $F_T$). The construction of $\delta_e$ together with both the regularity and the assumed degree of refinement of the mesh, allow us to assert that $\tilde{F}_T$ is contained in $F_T^{'}$. Therefore, recalling the local orthogonal frame $(O_H;x,y)$ of 
the plane of $F_T$, choosing for instance its origin to be a vertex of $T$ in $\Gamma$, $\tilde{\partial} T$ can be uniquely parametrized by 
the function $f_T$, in such a way that the spacial coordinates of any $P \in \tilde{\partial} T$  
are $(x,y,f_T(x,y))$ in the direct orthogonal spatial frame $(O_H;x,y,z)$ such that 
$z=0$ for points in $\tilde{F}_T$.
Let $\breve{u}_h$ be the function of $(x,y)$ defined in $\tilde{F}_T$ by $\breve{u}_h(x,y) = \bar{u}_h(x,y,f_T(x,y))$. \\
Notice that by construction, for all meshes and tetrahedra in ${\mathcal S}_h$ under consideration the ratio between the diameter of $\tilde{F}_T$ and the maximum diameter of the circles inscribed in this set is bounded above by a mesh-independent constant.  
Then, since $\breve{u}_h$ vanishes at $(k+2)(k+1)/2$ different points of $\tilde{F}_T$, from well-known results in interpolation theory in two-dimensional domains satisfying such a uniformity conditon (cf. \cite{BrennerScott}, Ch. 4, sect.4), there exists a mesh-independent constant $C_L$such that,
\begin{equation}
\label{auxiliary0}
\displaystyle \left[\int_{\tilde{F}_T} |\breve{u}_h(x)|^2 dxdy \right]^{1/2} \leq C_L h_T^{k+1} \displaystyle \left[\int_{\tilde{F}_T} \left|D_{x,y}^{k+1}\breve{u}_h(x,y)| \right|^2 dxdy \right]^{1/2}, 
\end{equation}
recalling that $D_{x,y}^j w$ is the $j$-th order tensor, whose components are the $j$-th order partial derivatives 
of a function $w$ with respect to $x$ and $y$.\\
Next we observe that, owing to Proposition \ref{prop02} there exist mesh-independent constants 
$c^j_{\Gamma}$ such that $\forall T \in {\mathcal S}_h$,
\begin{equation}
\label{auxiliary1}
\begin{array}{l}
\displaystyle \max_{M \in \tilde{F}_T} |D_{x,y}^j f_T(M)| \leq c^{j}_{\Gamma} h_T^{2-j}, \; j=1,\ldots,k+1.  
\end{array}
\end{equation}
On the other hand taking into account that the derivatives of $u_h$ of order greater than $k$ vanish in $T_{\Delta}$, straightforward calculations using the chain rule yield for suitable mesh-independent constants $c_{i}$, $i=0,1,\ldots,k$,
\begin{equation}
\label{auxiliary2}
| D_{x,y}^{k+1}\breve{u}_h(M)| \leq  c_{0} | D^{k+1}u | + \displaystyle 
\sum_{i=1}^{k} c_{i} h_T^{1-i} | D^{k+1-i} \bar{u}_h | \; \forall M \in \tilde{F}_T.  
\end{equation}
Notice that all the partial derivatives appearing on the right hand side of \eqref{auxiliary2} are to be understood at a (variable) point $P \in \tilde{\partial} T$ associated with $M \in \tilde{F}_T$. \\
Furthermore, owing to the bounds of the first partial derivatives of $f_T$, the surface element $dS$ on $\tilde{\partial} T$ equals $\xi(M) dxdy$ where $|\xi(M)| \leq \tilde{C}$ 
$\forall M \in \tilde{F}_T$ and conversely $dxdy = \zeta(P) dS$ with $|\zeta(P)| \leq \underline{C}$ $\forall P \in 
\tilde{\partial}T$, $\tilde{C}$ and $\underline{C}$ being independent of $T$.\\ 
On the other hand we observe that the union of all $\tilde{\partial}_T$ is nothing but $\Gamma$. Thus after straightforward calculations,  from \eqref{auxiliary0} and \eqref{auxiliary2} we come up with a mesh-independent constant $\tilde{C}_1$ such that,
\begin{equation}
\label{auxiliary4}
\parallel \bar{u}_h \parallel_{0,\Gamma}^2 
\leq \tilde{C}_1^2  \displaystyle \left[ h^{2(k+1)} \int_{\Gamma}  |D^{k+1}u|^2  + \sum_{T \in {\bf S}_h} h_T^{2(k+1)} 
\displaystyle \int_{\tilde{\partial} T} \displaystyle \sum_{j=1}^k h_T^{2(1-j)}|D^{k+1-j}\bar{u}_h|^2  \right].\\
\end{equation}
Now from the Trace Theorem \cite{Adams} we know that there exists a constant $C_{r}$ such that,
\begin{equation}
\label{auxiliary5}
\int_{\Gamma}  |D^{k+1}u|^2  \leq 
C_{r}^2 \parallel u \parallel_{k+1+r}^2
\end{equation}
On the other hand we clearly have.
\begin{equation}
\label{area}
area(\tilde{\partial} T) \leq \tilde{C} h_T^2/2.
\end{equation}
Hence using the curved element $T_{\Delta}$ associated with
$T$, we can write:
\begin{equation}
\label{auxiliary6}
\begin{array}{l}
\displaystyle \int_{\tilde{\partial} T}  
\displaystyle \sum_{j=1}^k h_T^{2(1-j)}|D^{k+1-j}\bar{u}_h|^2 
\leq \displaystyle \frac{\tilde{C}}{2} h_T^2 \displaystyle \sum_{j=1}^k h_T^{2(1-j)} 
\parallel D^{k+1-j}\bar{u}_h \parallel_{0,\infty,T_{\Delta}}^2.
\end{array}
\end{equation} 
Now setting ${\mathcal C}_{1} = {\mathcal C}_J$, according to Lemma \ref{wh} we have for $1 \leq j \leq k$ :
\begin{equation}
\label{skin}
\begin{array}{l}
\parallel D^{k+1-j}\bar{u}_h \parallel_{0,\infty,T_{\Delta}} \leq 
{\mathcal C}_{k+1-j} h_T^{j-k} \times \\
 h_T^{-3/2} \left(\parallel {\bf grad}\; \bar{u}_h \parallel_{0,T} 
 + h_T^k | u |_{k+1,T} + h_T^{k+r} \parallel u \parallel_{k+1+r,T} \right).
\end{array}
\end{equation}
Plugging \eqref{skin} into \eqref{auxiliary6} we come up with 
\begin{equation}
\label{auxiliary7}
\begin{array}{l}
\displaystyle \sum_{T \in {\mathcal S}_h} h_T^{2(k+1)} 
\displaystyle \int_{\tilde{\partial} T}  
 \sum_{j=1}^k h_T^{2(1-j)}|D^{k+1-j}\bar{u}_h|^2    
\leq {\bf C}_{k} \times \\
\displaystyle \sum_{T \in {\mathcal S}_h} \left( h_T^{3} \parallel {\bf grad}\; \bar{u}_h \parallel_{0,T}^2 
+ h_T^{2k+3} | u |_{k+1,T}^2 + h_T^{2k+3+2r} \parallel u \parallel_{k+1+r,T}^2 \right), 
\end{array}
\end{equation}
where ${\bf C}_{k}$ is another mesh-independent constant.\\
Now recalling \eqref{errestconvex} we can write:
\begin{equation}
\label{auxiliary8}
\displaystyle \sum_{T \in {\mathcal S}_h} h_T^{3} \parallel {\bf grad} \; \bar{u}_h \parallel_{0,T}^2 
\leq [{\mathcal C}(f,0)]^2 h^{2k+3} \; | u |_{k+1}^2.
\end{equation} 
Combining \eqref{auxiliary8} and \eqref{auxiliary7}   
and taking into account \eqref{auxiliary4}, \eqref{auxiliary5} and the fact that $h <1$, we easily obtain,
\begin{equation}
\label{auxiliary9}
\parallel \bar{u}_h \parallel_{0,\Gamma} 
\leq \bar{C}_1 h^{k+1} \displaystyle \left[ h^{1/2}| u |_{k+1} + \parallel u \parallel_{k+1+r} \right],
\end{equation} 
for a suitable mesh-independent constant $\bar{C}_1$.\\
It follows from \eqref{L2estimate1} and \eqref{auxiliary9} that for $C_{b1} = \bar{C}_1 C_t $ it holds:
\begin{equation}
\label{estimateb1}
b_{1h}(\bar{u}_h,v) \leq C_{b1} h^{k+1}[ h^{1/2} | u |_{k+1} + \parallel u \parallel_{k+1+r}] \parallel v \parallel_{2}.
\end{equation} 
\indent Now we turn our attention to $b_{2h}$.\\
First of all observing that ${\bf grad}\;\Pi_h(v)$ is constant in $T_{\Delta}$ for $T \in {\mathcal S}_h$ and $\Pi_h(v) =0$ on $\Gamma_h$, by Rolle's Theorem
\begin{equation}
\label{estim1b2}
|\Pi_h(v)(P)|\leq C_{\Gamma} h_T^2 \parallel {\bf grad} \; \Pi_h(v) \parallel_{0,\infty,T} \; \forall P \in \bar{\partial} T \mbox{ and } \forall T \in {\mathcal S}_h.
\end{equation}
On the other hand since for $h$ sufficiently small $area(F_T)  \leq area(\tilde{\partial}T)$ it holds 
$volume(\Delta_T) \leq \tilde{C} C_{\Gamma} h_T^4/2$. Thus \eqref{estim1b2} yields
\begin{equation}
\label{estim2b2}
b_{2h}(\bar{u}_h,\Pi_h{v}) \leq \displaystyle \frac{\tilde{C} C_{\Gamma}^2}{2}  \sum_{T \in {\mathcal S}_h}  h_T^6 \sqrt{3} \parallel H(\bar{u}_h) \parallel_{0,\infty,T_{\Delta}} 
\parallel {\bf grad} \; \Pi_h(v) \parallel_{0,\infty,T}.
\end{equation}
Using \eqref{LinftyTDelta} we further obtain:
\begin{equation}
\label{estim3b2}
b_{2h}(\bar{u}_h,\Pi_h{v}) \leq \displaystyle \tilde{C} C_{\Gamma}^2 {\mathcal C}_J  \sum_{T \in {\mathcal S}_h}  h_T^{9/2} 
\parallel H(\bar{u}_h) \parallel_{0,\infty,T_{\Delta}} 
\parallel {\bf grad} \; \Pi_h(v) \parallel_{0,T}.
\end{equation}
Next applying \eqref{Djwh} with $j=2$ we rewrite \eqref{estim3b2} as, 
\begin{equation}
\label{estim4b2} 
\begin{array}{l}
b_{2h}(\bar{u}_h,\Pi_h{v}) \leq \displaystyle \frac{\tilde{C} C_{\Gamma}^2 {\mathcal C}_J}{2}  \sum_{T \in {\mathcal S}_h}  h_T^{9/2} 
{\mathcal C}_2 h_T^{-5/2} \times \\
 \left( \parallel {\bf grad}\; \bar{u}_h \parallel_{0,T} + h_T^k | u |_{k+1,T} + h_T^{k+r} \parallel u \parallel_{k+1+r,T_{\Delta}} \right) 
\parallel {\bf grad} \; \Pi_h(v) \parallel_{0,T}.
\end{array}
\end{equation}
Now from standard interpolation results (cf. \cite{Ciarlet}) we know that for a mesh-independent constant $C_{\Pi}$ it holds, 
\begin{equation}
\label{estim5b2}
\parallel {\bf grad} \; \Pi_h(v) \parallel_{0,h} \leq C_{\Pi} \parallel v \parallel_{2}; \\
\end{equation} 
Thus, using the Cauchy-Schwarz inequality and recalling \eqref{errestconvex}, from \eqref{estim4b2} and \eqref{estim5b2} we easily infer the existence of a mesh-independent constant $C_{b2}$ such that,
\begin{equation}
\label{estimateb2} 
\begin{array}{l}
b_{2h}(\bar{u}_h,\Pi_h{v}) \leq C_{b2} h^{k+1}  
( h | u |_{k+1} + h^{1+r}  \parallel u \parallel_{k+1+r} ) \parallel v \parallel_{2}.
\end{array}
\end{equation}
\indent Next we estimate $b_{3h}$.\\
Recalling \eqref{b3h} and the fact that $\parallel {\bf grad}\;\Pi_h(v) \parallel_{0,\infty,T_{\Delta}} = \parallel {\bf grad}\;\Pi_h(v) \parallel_{0,\infty,T}$, we first take $\omega_T:=|{\bf grad}\;\bar{u}_{h_{|T_{\Delta}}}|$ and $\sigma_T = \bar{\partial} T$ for every $T \in {\mathcal S}_h$. Then since $\Pi_h(v) = 0$ on $F_T$, we first have,
\begin{equation}
\label{estim1b3}
b_{3h}(\bar{u}_h,\Pi_h(v)) \leq \displaystyle \sum_{T \in {\mathcal S}_h} 
\int_{\bar{\partial} T} \omega_T \Pi_h(v) \leq C_{\Gamma}  \displaystyle \sum_{T \in {\mathcal S}_h} h_T^2 \parallel {\bf grad} \; \Pi_h(v) \parallel_{0,\infty,T} \int_{\sigma_T} \omega_T.
\end{equation} 
Now using Lemma \ref{traceomegaT} together with \eqref{L2TDelta} and setting $\breve{C}_{3}:=C_{\Gamma} C_{\sigma} \tilde{C} {\mathcal C}_J$,
we easily conclude that,
\begin{equation}
\label{estim4b3}
b_{3h}(\bar{u}_h,\Pi_h(v)) \leq \breve{C}_{3} \displaystyle \sum_{T \in {\mathcal S}_h} 
h_T \parallel {\bf grad} \; \Pi_h(v) \parallel_{0,T} \left[ \int_{T_{\Delta}} [ \omega_T^2 + h_T^2|{\bf grad}\;\omega_T|^2 ] \right]^{1/2}. 
\end{equation}
Replacing $\omega_T$ by its expression in terms of ${\bf grad}\;\bar{u}_h$ and using the Cauchy-Schwarz inequality together with the fact that $volume(T_{\Delta}) \leq C_{\Delta}^{'} h_T^3$ by a straightforward geometric argument we obtain for $\bar{C}_3 =\breve{C}_{3}^2 C_{\Delta}^{'}$:   
\begin{equation}
\label{estim5b3}
| b_{3h}(\bar{u}_h,\Pi_h(v))|^2 \! \leq \! \bar{C}_3  \! \parallel \! {\bf grad} \; \Pi_h(v) \! \parallel_{0,h}^2 \! \displaystyle  \! \sum_{T \in {\mathcal S}_h}\! h_T^{5}\! \left( \parallel \! {\bf grad} \; \bar{u}_h \! \parallel_{0,\infty,T_{\Delta}}^2 \!\! + h_T^2 \! \parallel \! H(\bar{u}_h) \! \parallel_{0,\infty,T_{\Delta}}^2  \right).
\end{equation}
Now using \eqref{Djwh} with $j=1$ and $j=2$ together with \eqref{errestconvex}, elementary calculations lead to another mesh-independent constant $\bar{\mathcal C}_2$ such that, 
\begin{equation}
\label{estim7b3}
\displaystyle  \sum_{T \in {\mathcal S}_h}\! h_T^5 \! \left( \parallel \! {\bf grad}\;\bar{u}_{h} \! \parallel_{0,\infty,T_{\Delta}}^2 \! + h_T^2 \parallel \! H(\bar{u}_{h}) \! \parallel_{0,\infty,T_{\Delta}}^2 \right)  \leq \bar{\mathcal C}_2 h^{2k+2} \displaystyle \left( | u |_{k+1}^2 + h^{2r} \parallel u \parallel_{k+1+r}^2 \right).
\end{equation}
Finally plugging \eqref{estim7b3} into \eqref{estim5b3}, recalling \eqref{estim5b2} and setting $C_{b3}=[\bar{\mathcal C}_2 \bar{C}_{3}]^{1/2} C_{\Pi}$ we come up with, 
\begin{equation}
\label{estimateb3} 
b_{3h}(\bar{u}_h,\Pi_h(v)) \leq C_{b3} h^{k+1} (| u |_{k+1} + h^{r}  \parallel u \parallel_{k+1+r} )\parallel v \parallel_{2}. 
\end{equation}
\indent We pursue the proof with the estimation of $b_{4h}$.\\
To begin with we have,
\begin{equation}
\label{estim1b4}
b_{4h}(\bar{u}_h,v-\Pi_h(v)) \leq \displaystyle \sum_{T \in {\mathcal S}_h}  
\parallel {\bf grad} \; \bar{u}_h \parallel_{0,\Delta_T} \parallel {\bf grad}( v - \Pi_h(v)) 
\parallel_{0,T_{\Delta}}, 
\end{equation}
Furthermore we trivially have,
\begin{equation}
\label{estim2b4}
b_{4h}(\bar{u}_h,v-\Pi_h(v)) \leq \displaystyle \sum_{T \in {\mathcal S}_h} [volume(\Delta_T)]^{1/2}
\parallel {\bf grad}( v - \Pi_h(v))\parallel_{0,T_{\Delta}} \parallel {\bf grad} \; \bar{u}_h \parallel_{0,\infty,T_{\Delta}}. 
\end{equation}
Then using \eqref{Djwh} with $j=1$, from \eqref{estim2b4} we obtain
\begin{equation}
\label{estim3b4}
\begin{array}{l}
b_{4h}(\bar{u}_h,v-\Pi_h(v)) \leq {\mathcal C}_1 \displaystyle \sum_{T \in {\mathcal S}_h} [volume(\Delta_T)]^{1/2} \parallel {\bf grad}( v - \Pi_h(v))\parallel_{0,T_{\Delta}} \times \\
 h_T^{-3/2} \left(\parallel {\bf grad} \; \bar{u}_h \parallel_{0,T} + h_T^{k} | u |_{k+1,T} + h_T^{k+r} \parallel u \parallel_{k+1+r,T_{\Delta}} \right).
\end{array}
\end{equation}
Since $volume(\Delta_T) \leq C_{\Gamma} h_T^4$, from \eqref{estim3b4} we further obtain,
\begin{equation}
\label{estim4b4}
\begin{array}{l}
b_{4h}(\bar{u}_h,v-\Pi_h(v)) \leq  \displaystyle \sum_{T \in {\mathcal S}_h} 
{\mathcal C}_{1} C_{\Gamma}^{1/2} h_T^{2} \parallel {\bf grad}( v - \Pi_h(v)) \parallel_{0,T_{\Delta}} \times \\
 h_T^{-3/2} (\parallel {\bf grad} \; \bar{u}_h \parallel_{0,T} + h_T^{k} | u |_{k+1,T} + h_T^{k+r} \parallel u \parallel_{k+1+r,T} ).
\end{array}
\end{equation}
Now plugging \eqref{interP1} into \eqref{estim4b4} and applying the Cauchy-Schwarz inequality together with  \eqref{errestconvex}, we infer the existence of a mesh-independent constant $C_{b4}$ such that,
\begin{equation}
\label{estimateb4}
b_{4h}(\bar{u}_h,v-\Pi_h(v)) \leq C_{b4} h^{k+1} \displaystyle \left( h^{1/2} | u |_{k+1} 
+ h^{1/2+r} \parallel u \parallel_{k+1+r} \right)
| v |_{2}. 
\end{equation}
\indent Now we switch to the estimates of $c_h$ and $d_h$.\\
\indent As for $c_h$, we first observe that by the Trace Theorem the normal derivative of $v \in H^2(\Omega)$ across the interfaces of elements in ${\mathcal O}_h$ has no jumps. Thus roughly speaking the estimation of $c_h$ reduces to estimating the jumps of $\bar{u}_h$ on such interfaces. With this aim 
we resort to the operator $\tilde{I}_F$ defined in Lemma \ref{IF} where $F \in {\mathcal F}_h$.    
Notice that by construction $\tilde{I}_F(w)$ coincides on both sides of such an 
$F$ $\forall w \in W_h^0$, and clearly this property also holds for $\tilde{I}_F(u)$. Therefore we can write:
\begin{equation}
\label{estim1ch}
c_h(\bar{u}_h,v) = - \displaystyle \sum_{T \in {\mathcal R}_h} \sum_{F \in \partial T \setminus \tilde{F}_T}
\int_{F} (\bar{u}_h- \tilde{I}_F(\bar{u}_h))\frac{\partial v}{\partial n_T} - \displaystyle \sum_{T \in {\mathcal S}_h} \sum_{F \in \partial T \setminus F_T}
\int_{F} (\bar{u}_h- \tilde{I}_F(\bar{u}_h)) \frac{\partial v}{\partial n_T}. 
\end{equation}
Furthermore, since both $u$ and $I_F(u)$ coincide on both sides of any $F \in {\mathcal F}_h$, the former can be replaced by  the latter in \eqref{estim1ch}, or yet $\bar{u}_h$ can be replaced by $I_F(\bar{u}_h)$ therein. \\
Now we resort to Lemma \ref{traceomegaT} with $\sigma_T = \tilde{F}=F$ and $\omega_T = |{\bf grad} \; v |$ and to 
Lemma \ref{IF}. In doing so, after applying the Cauchy-Schwarz inequality to \eqref{estim1ch}, we easily obtain,
\begin{equation}
\label{estim2ch}
| c_h(\bar{u}_h,v) |^2 \leq [3{\mathcal C}_F C_{\sigma} ]^2 \! \displaystyle \sum_{T \in {\mathcal O}_h} \! h_T^4 \!
(\parallel {\bf grad} \; \bar{u}_h \parallel_{0,T} \! + h_T^{k} | u |_{k+1,T} + h_T^{k+r} \! \parallel u \parallel_{k+1+r,T_{\Delta}} )^2 \! \parallel v \parallel_2^2. 
\end{equation}
Finally using \eqref{errestconvex}, from \eqref{estim2ch} we come up with a mesh-independent constant $C_c$ such that,
\begin{equation}
\label{estimatec}
c_{h}(\bar{u}_h,v) \leq C_{c} h^{k+1} [h| u |_{k+1} + h^{1+r} \parallel u \parallel_{k+1+r}] \parallel v \parallel_{2}. 
\end{equation}
\indent In order to estimate $d_h$ we resort to Lemma \ref{deltae}. Indeed again by the Cauchy-Schwarz inequality and \eqref{deltatrace} we have
\begin{equation}
\label{estim1dh}
\left\{
\begin{array}{l}
d_h(\bar{u}_h,v) \leq C_{\theta} h \Lambda(\bar{u}_h) \parallel v \parallel_{2} \\
\mbox{where}\\
\Lambda(w) := \displaystyle \left[ \sum_{T \in {\mathcal S}_h} \sum_{e \subset F_T} 
\| w \|_{0,\delta_e}^2 \right]^{1/2}.
\end{array}
\right.
\end{equation}
Let us estimate $\Lambda(\bar{u}_h)$. \\
Since $\bar{u}_h=0$ at the end-points of $e$ and $area(\delta_e) \leq C_{\Gamma} h_T^3$ by Rolle's Theorem we clearly have $\| \bar{u}_h \|_{0,\delta_e}^2 \leq $ $C_{\Gamma} h_T^5 
\| {\bf grad} \; \bar{u}_h \|_{0,\infty,T_{\Delta}}^2$. Then by the same 
tricks already employed several times in this proof, in particular the use of \eqref{Djwh} with $j=1$ 
and \eqref{errestconvex}, without any difficulty the following estimate holds:
\begin{equation}
\label{estim2dh}
\Lambda(\bar{u}_h) \leq C_{\Lambda} h \left( h^{k} | u |_{k+1} + h^{k+r} \| u \|_{k+1+r} \right),
\end{equation} 
where $C_{\Lambda}$ is a mesh-independent constant.\\
Finally combining \eqref{estim1dh} and \eqref{estim2dh} and setting $C_d=C_{\theta} C_{\Lambda}$ we obtain,
\begin{equation}
\label{estimated}
d_h(\bar{u}_h,v) \leq C_{d} h^{k+1} \left( h| u |_{k+1} + h^{1+r} \| u \|_{k+1+r} \right)  
\| v \|_{2}.
\end{equation}
Plugging \eqref{estimateb1}, \eqref{estimateb2}, \eqref{estimateb3}, \eqref{estimateb4}, \eqref{estimatec} 
and \eqref{estimated} into \eqref{L2est4}, owing to the fact that $h <1$, we immediately obtain \eqref{L2estconvex} with 
${\mathcal C}_0=\tilde{\mathcal C}_0 + 2C(\Omega) (C_{b1}+C_{b2}+C_{b3}+C_{b4}+C_{c}+C_{d})$.  \QED 

\subsection{The case of non-convex domains} 
The case of a non-convex $\Omega$ is more delicate because the residual $a_h(u,v) - L_h(v)$ is not even defined for $v \in V_h$. 
Let us then consider a smooth domain $\Omega^{'}$ close to $\Omega$ which strictly contains $\Omega^{'}_h$ for all $h$ sufficiently small. More precisely, denoting by $\Gamma^{'}$ the boundary of $\Omega^{'}$ we assume that $|area(\Gamma^{'})-area(\Gamma)| \leq \varepsilon$ for $\varepsilon$ conveniently small. Henceforth we consider that $f$ was also extended to $\Omega^{'} \setminus \Omega$. We denote the extended $f$ by $f^{'}$, which is arbitrarily chosen, except for the requirement that $f^{'} \in H^{k-1}(\Omega^{'} )$. There are different ways to achieve such a regularity and in this respect the author refers for instance to 
\cite{Lions} or \cite{Necas}. \\
\indent Then instead of (\ref{Poisson-h}) we solve:
\begin{equation}
\label{uhprime} 
\left\{
\begin{array}{l}
\mbox{Find } u_h \in W_h^g \mbox{ such that } \\
a_h(u_h,v)=L^{'}_h(v) := \int_{\Omega_h} f^{'} v$ $\forall v \in V_h. 
\end{array}
\right.
\end{equation}
Akin to problem (\ref{Poisson-h}) and thanks to (\ref{inf-sup}), problem (\ref{uhprime}) has a unique solution.  
This fact allows us to claim the following preliminary result:  
\begin{theorem}
\label{convergence}
Assume that for $f^{'} \in H^{k-1}(\Omega^{'})$ there exists a function $u^{'}$ defined in $\Omega^{'}$ having the following properties:
\begin{itemize}
\item 
$-\Delta u^{'} = f^{'}$ in $\Omega^{'}$;
\item
$u^{'}_{|\Omega} = u$;
\item
$u^{'} = g$ a.e. on $\Gamma$; 
\item
$u^{'} \in H^{k+1}(\Omega^{'})$.
\end{itemize}
Then for $k>1$ and a suitable constant ${\mathcal C}^{'}$ independent of $h$ it holds:
\begin{equation}
\label{errorestimate}
\parallel {\bf grad}_h(u - u_h) \parallel_{\widetilde{0,h}} \leq {\mathcal C}^{'} | u^{'} |_{k+1,\Omega^{'}} h^k.
\end{equation} 
\end{theorem}

\prov Here, instead of adapting the distance inequalities in (\cite{COAM}) to this specific situation, we employ a more straightforward argument. 
First we recall (\ref{inf-sup}) to note that $\forall w \in W^g_h$ we have:
\begin{equation}
\label{firstbound}
\parallel {\bf grad}_h(u_h - w) \parallel_{0,h} \leq \displaystyle \frac{1}{\alpha} \displaystyle \sup_{v \in V_h \setminus \{0\}} 
\frac{|a_h(u_h,v)-a_h(w,v)|}{\parallel {\bf grad}\; v \parallel_{0,h}}.
\end{equation}
Since $a_h(u^{'},v)=-\int_{\Omega_h} v \Delta u^{'}= L^{'}_h(v)=a_h(u_h,v) \; \forall v \in V_h$ we can further write for every $w \in W^g_h$:
\begin{equation}
\label{secondbound}
\parallel {\bf grad}_h(u_h - w) \parallel_{0,h} \leq \displaystyle \frac{1}{\alpha} \displaystyle \sup_{v \in V_h \setminus \{0\}} 
\frac{|a_h(u^{'}-w,v)|}{\parallel {\bf grad}\; v \parallel_{0,h}} \leq \displaystyle \frac{1}{\alpha} \parallel {\bf grad}_h(u^{'} - w) \parallel_{0,h}. 
\end{equation}
From the triangle inequality this further yields:
\begin{equation}
\label{thirdbound}
\parallel {\bf grad}_h(u_h - u^{'}) \parallel_{0,h} \leq \displaystyle \left[ 1 + \frac{1}{\alpha} \right] \parallel {\bf grad}_h(u^{'} - w) \parallel_{0,h}. 
\end{equation}
Choosing $w$ to be the $W^g_h$-interpolate of $u^{'}$ in $\Omega_h$, and using standard interpolation results (cf. \cite{BrennerScott}), from 
(\ref{thirdbound}) we establish (\ref{errorestimate}). \QED \\

In principle the knowledge of a regular extension $f^{'}$ of the right hand side datum $f$ associated with a regular extension $u^{'}$ of $u$ is necessary to solve problem (\ref{uhprime}). However in most practical cases, neither such an extension of $f$, nor $u^{'}$ 
satisfying the assumptions of Theorem \ref{convergence} associated with a given regular extension $f^{'}$ of $f$ 
is known. 
Nevertheless using some results available in the literature it is possible to identify cases where such an extension $u^{'}$ does exist. Let us consider for instance a simply connected domain $\Omega$ of the $C^{\infty}$-class and a datum $f$ 
infinitely differentiable in $\bar{\Omega}$. Taking an extension $f^{'} \in C^{\infty}(\Omega^{'}) \cap H^{k-1}(\Omega^{'})$ of $f$ to an enlarged domain 
$\Omega^{'}$ also of the $C^{\infty}$-class, we first solve $-\Delta u_0 = f^{'}$ in $\Omega^{'}$ and $u_0 = 0$ on $\Gamma^{'}$. According to well-known results (cf. \cite{LionsMagenes}) $u_0 \in C^{\infty}({\Omega}^{'})$ and hence the trace $g_0$ of $u_0$ on $\Gamma$ belongs to $C^{\infty}(\Gamma)$. Next we denote by $u_H$ the harmonic function in $\Omega$ such that $u_H = g_0$ on $\Gamma$. Let $r_0$ be the radius of the largest (open) ball $B$ contained in $\Omega$ and $O=(x_0,y_0,z_0)$ be its center. Assuming thet $f^{'}$ is not too wild, so that the Taylor series of $u_H(x,y,z_0)$ and 
$[\partial u_H/\partial z](x,y,z_0)$ centered at $O$ converge in a disk of the plane $z=z_0$ centered at $O$ with radius equal to 
$r_0 \sqrt{2}+ \delta$ for a certain $\delta > 0$, according to \cite{Coffman} there exists a harmonic extension of $u^{'}_H$ to the ball $B_0^{'}$ centered at $O$ with radius $r_0 +\delta \sqrt{2}$. Clearly in this case, as long as $\delta$ is large enough for $B^{'}$ to contain $\Omega^{'}$, we can define 
$u_0^{'}:=u_0-u_H^{'}$ as a function in $H^{k+1}(\Omega^{'})$ that vanishes on $\Gamma$. Now further assuming that $g \in C^{\infty}(\Gamma)$ we can also define an extension of the harmonic function $u^H$ whose value is $g$ on $\Gamma$ into $u^{H'} \in H^{k+1}(\Omega^{'})$ in the very same manner as $u_H$. The extension $u^{'}$ 
of $u$ to $\Omega^{'}$ given by $u^{'} := u^{H'} + u^{'}_0$ satisfies the required properties.\\ 
In the general case however, a convenient way to bypass the uncertain existence of an extension $u^{'}$ satisfying the assumptions of Theorem \ref{convergence}, is to resort to numerical integration on the right hand side. Under certain conditions rather easily satisfied, this leads to the definition of an alternative approximate problem, in which only values of $f$ (in $\Omega$) come into play. This trick is inspired by the one of Ciarlet and Raviart in their work on the isoparametric finite element method (cf. \cite{CiarletRaviart} and \cite{Ciarlet}). To be more specific, these celebrated authors employ the following argument, assuming that $h$ is small enough: if a numerical integration formula is used, which has no integration points different from vertices on the faces of a tetrahedron, then only values of $f$ (in $\Omega$) will be needed to compute the corresponding approximation of $L^{'}_h(v)$. This means that the knowledge of $u^{'}$, and thus of $f^{'}$, will not be necessary for implementation purposes. Moreover, provided the accuracy of the numerical integration formula is compatible with method's order, the resulting modification of (\ref{uhprime}) will be a method of order $k$ in the norm  
$\parallel \cdot \parallel_{\widetilde{0,h}}$ of ${\bf grad} \; u - {\bf grad}_h u_h$. \\
Nevertheless it is possible to get rid of the above argument based on numerical integration in the most important cases in practice, namely, those of quadratic and cubic Lagrange finite elements. Let us see how this works. \\
First of all we consider that $f$ is extended by zero in $\Delta_{\Omega}:=\Omega^{'} \setminus \bar{\Omega}$, and resort to the extension $u^{'}$ of $u$ to the same set constructed in accordance to Stein et al. \cite{Stein}. This extension does not satisfy $\Delta u^{'}=0$ in $\Delta_{\Omega}$ but the function denoted in the same way such that $u^{'}_{|\Omega} = u$ does belong to $H^{k+1}(\Omega^{'})$. Since $k > 1$ this means in particular that the traces of the functions $u$ and $u^{'}$ coincide on 
$\Gamma$ and that $\partial u /\partial n = - \partial u^{'}/\partial n^{'} = 0$ a.e. on $\Gamma$ where the normal derivatives on the right hand side of this relation is the outer normal derivative with respect to $\Delta_{\Omega}$ (the trace of the Laplacian of both functions also coincide  on $\Gamma$ but this is not relevant for our purposes). Based on this extension of $u$ to $\Omega_h$ for all such polyhedra of interest, we next prove the following results for the approximate problem (\ref{Poisson-h}), without assuming that $\Omega$ is convex, and still denoting by $f$ the function identical to the right hand side datum of (\ref{Poisson}) in $\Omega$, that vanishes identically in $\Delta_{\Omega}$.

\begin{theorem} 
\label{P2}
If $k=2$ there exists a mesh independent constant $C_2$ such that the unique solution $u_h$ to (\ref{Poisson-h}) satisfies:
\begin{equation}
\label{estimateP2} 
\begin{array}{l}
\parallel {\bf grad}_h(u - u_h) \parallel_{\widetilde{0,h}} \leq C_2 h^{2} G(u^{'}) \\
\\
\mbox{with } G(u^{'}):=| u^{'} |_{3,\Omega^{'}} + h^{1/2} \parallel \Delta u^{'} \parallel_{0,\Omega^{'}}, 
\end{array}
\end{equation}
$u^{'} \in H^3(\Omega^{'})$ being the regular extension of $u$ to $\Omega^{'}$ constructed in accordance to Stein et al. \cite{Stein}.
\end{theorem}

\prov
First we recall (\ref{firstbound}), from which we obtain:
\begin{equation}
\label{fourthbound}
\parallel {\bf grad}_h(u_h - w) \parallel_{0,h} \leq \displaystyle \frac{1}{\alpha} \displaystyle \sup_{v \in V_h \setminus \{0\}} 
\frac{|a_h(u^{'},v)-L_h(v)| + |a_h(u^{'}-w,v)|}{\parallel {\bf grad}\; v \parallel_{0,h}}.
\end{equation}
Thanks to the following facts the first term in the numerator of (\ref{fourthbound}) can be dealt with in the following manner: Since $u^{'} \in H^3(\Omega^{'})$ 
we can apply First Green's identity to $a_h(u^{'},v)$ thereby getting rid of integrals on portions of $\Gamma$; next defining $\Delta^{'}_T:=T \setminus \Omega$, we note that 
$\Delta u + f=0$ in every $T \in {\mathcal T}_h \setminus {\mathcal O}_h$; this is also 
true of elements $T$ not belonging to the subset ${\mathcal Q}_h$ of ${\mathcal O}_h$ consisting of elements $T$ such that $\Delta^{'}_T$ is not restricted to a set of vertices of $\Omega_h$; finally we recall that $\Delta u^{'} + f$ vanishes identically in $\tilde{T}$ and observe that the interior of $\Delta^{'}_T$ is not empty $\forall T \in {\mathcal Q}_h$. In short we can write:  
\begin{equation}
\label{ahFhprime}
|a_h(u^{'},v)-L_h(v)| = \displaystyle \sum_{T \in {\mathcal Q}_h} \int_{\Delta^{'}_T} -\Delta u^{'} v \;  
\leq \displaystyle \sum_{T \in {\mathcal Q}_h} \parallel \Delta u^{'} \parallel_{0,\Delta^{'}_T} \parallel v \parallel_{0,\Delta^{'}_T}.
\end{equation}   
Let us first consider the case where $T \in {\mathcal S}_h \cap {\mathcal Q}_h$. We recall that for the mesh-independent constant $C_{\Gamma}$ it holds
\begin{equation}
\label{Rolle} 
|v({\bf x})| \leq C_{\Gamma} h_T^2 \parallel {\bf grad} \; v \parallel_{0,\infty,\Delta^{'}_T}, \; 
\forall {\bf x} \in \Delta^{'}_T. 
\end{equation}
On the other hand from \eqref{L2TDelta} we infer that
$\parallel {\bf grad} \; v \parallel_{0,\infty,\Delta^{'}_T} \leq {\mathcal C}_J h_T^{-3/2} \parallel {\bf grad} \; v \parallel_{0,T}$. Then noticing that $volume(\Delta^{'}_T)$ is bounded by a constant depending only on $\Omega$ multiplied by $h_T^4$ and using \eqref{Rolle}, we obtain for a certain mesh-independent constant $C_Q$:
\begin{equation}
\label{fifthbound}
\parallel \Delta u^{'} \parallel_{0,\Delta^{'}_T} \parallel v \parallel_{0,\Delta^{'}_T} \leq C_Q h_T^{5/2} \parallel \Delta u^{'} \parallel_{0,\Delta^{'}_T} \parallel 
{\bf grad} \; v \parallel_{0,T} \; \forall T \in {\mathcal Q}_h \cap {\mathcal S}_h.
\end{equation}
Now we consider the elements $T$ in the set ${\mathcal Q}_h \cap {\mathcal R}_h$. Since in this case the measure of $\Delta^{'}_T$ is bounded above by a constant depending only on $\Omega$ multiplied by $h_T^5$, we obtain for such elements a bound similar to (\ref{fifthbound}) with $h_T^3$ instead of $h_T^{5/2}$. 
Since $h_T << 1$ by assumption we can assert that (\ref{fifthbound}) also holds for elements in this set.\\
Now plugging (\ref{fifthbound}) into (\ref{ahFhprime}) and applying the Cauchy-Schwarz inequality, we easily come up with,
\begin{equation}
\label{sixthbound}
|a_h(u^{'},v)-L_h(v)| \leq C_Q h^{5/2} \displaystyle \parallel \Delta u^{'} \parallel_{0,\Omega^{'}} \parallel {\bf grad} \; v \parallel_{0,h}.
\end{equation}
Finally combining (\ref{sixthbound}) and (\ref{fourthbound}) and using the triangle inequality we easily establish the validity of error estimate (\ref{estimateP2}). \QED \\ 

\begin{theorem} 
\label{P3}
If $k=3$ there exists a mesh independent constant $C_3$ such that the unique solution $u_h$ to (\ref{Poisson-h}) satisfies:
\begin{equation}
\label{estimateP3} 
\parallel {\bf grad}_h(u - u_h) \parallel_{\widetilde{0,h}} \leq C_3 h^{3}[ | u^{'} |_{4,\Omega^{'}} + h^{1/2} \parallel \Delta u^{'} \parallel_{0,\infty,\Omega^{'}}]
\end{equation}
where $u^{'} \in H^4(\Omega^{'})$ is the regular extension of $u$ to $\Omega^{'}$ constructed in accordance to Stein et al. \cite{Stein}.
\end{theorem}

\prov 
First of all we point out that, according to the Sobolev Embedding Theorem \cite{Adams}, $\Delta u^{'} \in L^{\infty}(\Omega^{'})$, since $u^{'} \in H^4(\Omega^{'})$ by assumption. \\
Following the same steps as in the proof of Theorem \ref{P2} up to equation (\ref{ahFhprime}), the latter becomes for a certain mesh-independent constant $C_R$,
\begin{equation}
\label{ahFh3}
|a_h(u^{'},v)-L_h(v)| \leq C_R \displaystyle \sum_{T \in {\mathcal Q}_h} h_T^{4} \parallel \Delta u^{'} \parallel_{0,\infty,\Omega^{'}} \parallel v \parallel_
{0,\infty,\Delta^{'}_T}. 
\end{equation}
Using the same arguments leading to (\ref{fifthbound}) this yields in turn, for a constant $C_S$ equal to $C_{\Gamma} C_R {\mathcal C}_J$:
\begin{equation}
\label{ahFh4}
|a_h(u^{'},v)-L_h(v)| \leq C_S \displaystyle \sum_{T \in {\mathcal Q}_h} h_T^{9/2} \parallel \Delta u^{'} \parallel_{0,\infty,\Omega^{'}} 
\parallel {\bf grad} \; v \parallel_{0,T}. 
\end{equation}
Further appying the Cauchy-Schwarz inequality to the right hand side of (\ref{ahFh4}) we easily obtain:
\begin{equation}
\label{ahFh5}
|a_h(u^{'},v)-L_h(v)| \leq C_S h^{7/2} \parallel \Delta u^{'} \parallel_{0,\infty,\Omega^{'}} \displaystyle \left[ \sum_{T \in {\mathcal Q}_h} h_T^{2} \right]^{1/2} 
\parallel {\bf grad} \; v  \parallel_{0,h}. 
\end{equation}
From the fact that the family of meshes in use is regular we know that 
\begin{equation}
\label{skinvolume}
\displaystyle \left[ \sum_{T \in {\mathcal Q}_h} h_T^{2} \right]^{1/2} \leq C_{\Gamma}^{'} \mbox{ independently of }h. 
\end{equation} 
Plugging \eqref{skinvolume} into (\ref{ahFh5}) and the resulting 
relation into (\ref{fourthbound}) we immediately establish error estimate (\ref{estimateP3}). \QED \\

A simple and useful consequence of Theorems \ref{P2} and \ref{P3} is the following,
\begin{corollary}
The solution $u_h$ of \eqref{Poisson-h} satisfies,
\begin{equation}
\label{estimatePk} 
\parallel {\bf grad}_h(u^{'} - u_h) \parallel_{0,h} \leq C_k h^{k} G(u^{'})
\end{equation}
where 
\begin{equation}
\label{Gk}
\left\{
\begin{array}{l}
G(u^{'}) =  | u^{'} |_{3,\Omega^{'}} + h^{1/2} \parallel \Delta u^{'} \parallel_{0,\Omega^{'}} \mbox{ for } k=2\\
\mbox{and} \\
G(u^{'})=| u^{'} |_{4,\Omega^{'}} + h^{1/2} \parallel \Delta u^{'} \parallel_{0,\infty,\Omega^{'}}\mbox{ for } k=3. 
\end{array}
\right.
\end{equation}
\end{corollary}

\prov Estimate \eqref{estimatePk} trivially results from \eqref{fourthbound} if we 
add and subtract $u^{'}$ inside the norm on the left hand side and apply the triangle inequality. \QED \\

Akin to Theorem \ref{theorem1bis}, it is possible to establish error estimates in the $L^2$-norm in the case of a non-convex $\Omega$, by requiring some more regularity from the solution $u$ of \eqref{Poisson}. However,  unless the assumptions of Theorem \ref{convergence} hold, optimality is not attained for $k>2$. This is because of the absence of $u$ from the nonempty domain $\Delta^{'}_h : = \Omega_h \setminus \Omega$, whose volume is an invariant $O(h^2)$ whatever $k$. Roughly speaking, integrals in $\Delta^{'}_h$ of expressions in terms of the approximate solution $u_h$ dominate the error, in such a way that those terms cannot be reduced to less than an $O(h^{7/2})$, even under additional regularity assumptions. \\
\indent Most steps in the proof of the following result rely on arguments essentially identical to those already exploited to prove Theorem \ref{theorem1bis}. Therefore we will focus on aspects specific to the non-convex case. \\
The proof of error estimates in the $L^2$-norm is rather long. Thus for the sake of brevity, and without loss of essential results we confine ourselves here again to the case of homogeneous boundary conditions. In addition to this, in order to avoid technicalities even more intricate than those already involved in the proofs of our $L^2$-error estimates, we shall make the additional assumption on the mesh specified in the theorem that follows. 
However we observe that besides being reasonable, such an assumption is by no means essential for the underlying result to hold.\\
Furthermore, although this is by no means mandatory, the proof of the following theorem is significantly simplified if we assume that $u^{'} \in H^{3+r}(\Omega^{'})$.

\begin{theorem}
\label{L2P2} Let $k=2$ and $g \equiv 0$. Assume that the mesh is such that every pair of elements in 
${\mathcal R}_h$ has no common face. 
Further assume that $\Omega$ is of the piecewise $C^3$-class and $u \in H^{3+r}(\Omega)$ for 
$r=1/2+\epsilon$, $\epsilon>0$ being arbitrarily small and consider the 
extension $u^{'}$ of $u$ to $\Omega^{'}$ in $H^{3}(\Omega^{'})$ 
constructed in accordance to Stein et al. \cite{Stein}. Then the following error estimate holds:
\begin{equation}
\label{L2estP2} 
\parallel u - u_h \parallel_{\widetilde{0,h}} \leq C^{'}_0 h^{3} [ G(u^{'})
+ \parallel u^{'} \parallel_{3+r,\Omega^{'}} ],   
\end{equation} 
where $C^{'}_0$ is a mesh-independent constant and $G(u^{'})$ is given by \eqref{Gk}. \\
\end{theorem}  

\prov Let $\bar{u}_h$ be the function defined in $\Omega$ by $\bar{u}_h:=u_h-u$.  
$v \in H^1_0(\Omega)$ being the function satisfying \eqref{adjoint}-\eqref{adjoint1}, we have: 
\begin{equation}
\label{L2estP2bis}
\parallel \bar{u}_h \parallel_{\widetilde{0,h}} \leq \parallel \bar{u}_h \parallel_{0} \leq 
C_{\Omega} \displaystyle \frac{- \int_{\Omega} \bar{u}_h \Delta v}{\parallel v \parallel_{2}}.
\end{equation}
First of all we observe that $\Gamma^{'}_h = \displaystyle \cup_{T \in {\mathcal Q}_h} \Gamma_{T}$ and moreover in the case under study $area(\Gamma^{'}_h) >0$. \\
Now using integration by parts we obtain,  
\begin{equation}
\label{L2estP2ter}
\parallel \bar{u}_h \parallel_{\widetilde{0,h}} \leq 
C_{\Omega} \displaystyle \frac{\displaystyle \tilde{a}_h(\bar{u}_h,v)+a_{\Delta_h}(\bar{u}_h,v)
- \tilde{a}_{\partial h}(\bar{u}_h,v)}{\parallel v \parallel_{2}},
\end{equation}
where the bilinear form $a_{\Delta_h}$ is defined in \eqref{aDeltah} and for 
$w \in W_h^0 + H^1(\Omega)$, $z \in H^1(\Omega)$ and $v\in H^2(\Omega)$,   
\begin{equation}
\label{tildeah}
\tilde{a}_h(w,z) 
:= \int_{\tilde{\Omega}_h}  {\bf grad}_h w \cdot {\bf grad} \;z,  
\end{equation}
and
\begin{equation}
\label{tildeapartialh}
\tilde{a}_{\partial h}(w,v) := \displaystyle \sum_{T \in {\mathcal R}_h} 
\int_{[(\partial T \cap \Omega)\setminus \tilde{F}_T] \cup \Gamma_T} w \frac{\partial v}{\partial n_T} + \displaystyle \sum_{T \in {\mathcal S}_h} 
\left[\int_{\bar{\partial} T \cap \bar{\Omega}}  w \frac{\partial v}{\partial \bar{n}_T} 
\! + \! \int_{(\partial T \cap \Omega) \setminus F_T}  w \frac{\partial v}{\partial n_T} \right].
\end{equation}
Similarly to \eqref{apartialh} we write 
\begin{equation}
\label{splitapartial}
\tilde{a}_{\partial h}(w,v) =  - c^{'}_h(w,v) - d^{'}_h(w,v) - b_{1h}(w,v),
\end{equation} 
where $b_{1h}(w,v)$ is defined by \eqref{b1h} and $\forall w \in W_h^0 + H^1(\Omega)$ and $\forall v \in H^2(\Omega)$, 
\begin{equation}
\label{tildegammah}
c^{'}_h(w,v) := - \displaystyle \sum_{T \in {\mathcal R}_h} 
\int_{(\partial T \cap \Omega) \setminus \tilde{F}_T} w \frac{\partial v}{\partial n_T} - \displaystyle \sum_{T \in {\mathcal S}_h} 
\int_{(\partial T \cap \Omega) \setminus F_T}  w \frac{\partial v}{\partial n_T} 
\; 
\end{equation}
\begin{equation}
\label{tildedeltah}
d^{'}_h(w,v) := - \displaystyle \sum_{T \in {\mathcal S}_h}\sum_{e \subset F_T} \int_{\delta_e^{'}}  w \frac{\partial v}{\partial \bar{n}_T}.
\end{equation}  
Notice that, in contrast to the convex case (the closure of) $\Delta_h$ is not the union of the sets $\Delta_T$ as $T$ sweeps ${\mathcal S}_h$, but rather $\bar{\Delta}_h := \displaystyle \cup_{T \in {\mathcal S}_h} \tilde{\Delta}_T$, bearing in mind that the interior of $\tilde{\Delta}_T$ can obviously be an empty set for certain tetrahedra in ${\mathcal S}_h$. \\
\indent On the other hand, denoting by $\partial \cdot/\partial \tilde{n}_h$ the outer normal derivative on $\tilde{\Gamma}_h$, 
$\forall v_h \in V_h$ we have,
\begin{equation}
\label{L2estP2qua}
a_h(u_h,v_h)= -\int_{\tilde{\Omega}_h} v_h \Delta u =  
- \displaystyle \oint_{\tilde{\Gamma}_h} \frac{\partial u}{\partial \tilde{n}_h} v_h + \tilde{a}_h(u,v_h).
\end{equation}
But since any function in $V_h$ vanishes identically on $\Gamma_h$, recalling the definition of ${\mathcal Q}_h$ in the proof of Theorem \ref{P2} together with the set $\Delta^{'}_T = T \setminus \tilde{T}$, $\forall T \in {\mathcal Q}_h$ we necessarily have,
\begin{equation}
\label{L2estP2quaqua}
\tilde{a}_h(u_h,v_h) + \displaystyle \sum_{T \in {\mathcal Q}_h} \int_{\Delta^{'}_T} {\bf grad} \; u_h \cdot 
{\bf grad} \; v_h = a_h(u_h,v_h)=   
- \displaystyle \int_{\Gamma^{'}_h} \frac{\partial u}{\partial n}v_h + \tilde{a}_h(u,v_h).
\end{equation}
In doing so we define, 
\begin{equation}
\label{b5h}
b_{5h}(w,z) := \displaystyle \sum_{T \in {\mathcal Q}_h} \int_{\Delta^{'}_T} z \Delta w \;
\forall w \in W_h^0 \mbox{ and } \forall z \in V_h.
\end{equation}
together with,
\begin{equation}
\label{b6h}
b_{6h}(w,z) := \displaystyle \sum_{T \in {\mathcal Q}_h} 
\int_{\Gamma_T} \frac{\partial w}{\partial n}z \; \forall  
w \in W_h^0 + H^2(\Omega) \mbox{ and } z \in V_h.
\end{equation}
Then applying integration by parts in $\Delta^{'}_T$ it easily follows from \eqref{L2estP2quaqua} that 
\begin{equation}
\label{L2estP2qui}
-\tilde{a}_h(\bar{u}_h,v_h) + b_{5h}(u_h,v_h) + b_{6h}(\bar{u}_h,v_h)=0 \; \forall v_h \in V_h.
\end{equation}
Taking $v_h = \Pi_h(v)$, using again the error function $e_h(v)=v-\Pi_h(v)$ and plugging \eqref{L2estP2qui} into \eqref{L2estP2ter} we come up with,
\begin{equation}
\label{L2estP2sex}
\parallel \bar{u}_h \parallel_{\widetilde{0,h}} \leq 
C_{\Omega} \displaystyle \frac{ 
b_{5h}(u_h,\Pi_h(v))\!+\!b_{6h}(\bar{u}_h,\Pi_h(v))\!+\!\tilde{a}_h(\bar{u}_h,e_h(v))\!-\!\tilde{a}_{\partial h}(\bar{u}_h,v)\!+\!a_{\Delta_h}(\bar{u}_h,v)}{\parallel v \parallel_{2}},
\end{equation}
where $\tilde{a}_{\partial h}$ is split as indicated in \eqref{splitapartial}. \\
Next we redefine $b_{2h}$ given by \eqref{b2h} in order to take into account the sets   
$\tilde{\Delta}_T$ rather than $\Delta_T$ as follows:  
\begin{equation}
\label{b2hprime}
b^{'}_{2h}(w,z):=\displaystyle \sum_{T \in {\mathcal S}_h} - \int_{\tilde{\Delta}_T} z \Delta w \mbox{ for } w \in W_h^0 + H^2(\Omega) \mbox{ and } z \in H^{1}(\Omega),  
\end{equation}
Now $\partial \tilde{\Delta}_T$ being the boundary of $\tilde{\Delta}_T$ let $\partial^{'} \tilde{\Delta}_T:=
\partial \tilde{\Delta}_T \setminus \Gamma_h$. Then from the fact that $\Pi_h(v)$ vanishes identically on $\Gamma_h$, using integration by parts and recalling that the notation $\partial \cdot / \partial \bar{n}_T$ is used to represent the normal derivative on $\partial^{'} \tilde{\Delta}_T$ directed outwards $\tilde{\Delta}_T$, we obtain:
\begin{equation}
\label{L2estP2sept}
a_{\Delta_h}(\bar{u}_h,v)=a_{\Delta_h}(\bar{u}_h,e_h(v)) +  
 \displaystyle \sum_{T \in {\mathcal S}_h} \int_{\partial^{'} \tilde{\Delta}_T} \displaystyle 
 \frac{\partial \bar{u}_h}{\partial \bar{n}_T} \Pi_h(v) + b^{'}_{2h}(\bar{u}_h,\Pi_h(v)).
\end{equation} 
Next akin to $b_{2h}$ we adjust the definition \eqref{b3h} of $b_{3h}$ for $w \in W_h^0 + H^2(\Omega)$ and 
$z \in H^1(\Omega)$ into 
\begin{equation}
\label{b3hprime}
b_{3h}^{'}(w,z) := \displaystyle \sum_{T \in {\mathcal R}_h} \int_{\Gamma_T} 
\frac{\partial w}{\partial \bar{n}_T} z 
+ \displaystyle \sum_{T \in {\mathcal S}_h} \int_{\partial^{'} \tilde{\Delta}_T 
\cup \Gamma_T} \frac{\partial w}{\partial \bar{n}_T} z 
\end{equation}
recalling also that $\partial \cdot /\partial \bar{n}_T$ coincides with $\partial \cdot /\partial n$ on 
$\Gamma_T$ $\forall T \in {\mathcal O}_h$. \\
Actually from \eqref{b3hprime} and since $\Gamma_T = \emptyset$ if $T \notin {\mathcal Q}_h$, we conclude that
\begin{equation}
\label{b36} 
\displaystyle \sum_{T \in {\mathcal S}_h} \int_{\bar{\partial}\tilde{\Delta}_T} \displaystyle 
 \frac{\partial \bar{u}_h}{\partial \bar{n}_T} \Pi_h(v) + b_{6h}(\bar{u}_h,\Pi_h(v)) = b_{3h}^{'}(\bar{u}_h,\Pi_h(v)).
\end{equation}
Thus recalling the definition \eqref{b4h} of $b_{4h}$, \eqref{L2estP2sept} \eqref{b36} combined with 
\eqref{L2estP2sex} yields:
\begin{equation}
\label{L2estP2oct}
\left\{
\begin{array}{l}
\parallel \bar{u}_h \parallel_{\widetilde{0,h}} \leq 
C_{\Omega} \displaystyle \frac{\tilde{a}_h(\bar{u}_h,e_h(v))+b_{1h}(\bar{u}_h,v)+b_{5h}(u_h,v_h)+{\mathcal L}(\bar{u}_h,v)}{\parallel v \parallel_{2}},\\
\mbox{where }\\
{\mathcal L}(\bar{u}_h,v):= b^{'}_{2h}(\bar{u}_h,\Pi_h(v)) + b_{3h}^{'}(\bar{u}_h,\Pi_h(v))+ 
b_{4h}(\bar{u}_h,e_h(v))+ c^{'}_h(\bar{u}_h,v) + d^{'}_h(\bar{u}_h,v).
\end{array}
\right.
\end{equation}
The estimation of $\tilde{a}_h(\bar{u}_h,e_h(v))$ is a trivial variant of the one in Theorem 
\ref{theorem1bis}, that is,
\begin{equation}
\label{estildeah}
\tilde{a}_h(\bar{u}_h,e_h(v)) \leq C_2 \tilde{C}_V h^{3} G(u^{'}) |v|_{2},
\end{equation}
where $\tilde{C}_V$ is an interpolation error constant such that 
\begin{equation}
\label{intildeP1}
\parallel {\bf grad} (v-\Pi_h(v))\parallel_{\widetilde{0,h}} \leq \tilde{C}_V h | v |_{2}.
\end{equation}
\indent The bilinear form $b_{4h}$ can be estimated like in Theorem \ref{theorem1bis} with 
minor modifications. The estimates of the bilinear forms $b_{2h}^{'}$, $b_{3h}^{'}$, $c^{'}_h$ and $d^{'}_h$ also follow the main lines of those of $b_{2h}$, $b_{3h}$, $c_h$ and $d_h$, respectively given in Theorem \ref{theorem1bis}, taking $k=2$. Through the use of \eqref{estimatePk} instead of \eqref{errestconvex} we come up with final results of the same qualitative nature. Actually if we replace $u$ with $u^{'}$ in the different intermediate steps that come into play they become practically the same. Keeping in mind that $u^{'}$ is 
claimed to be in $H^{3+r}(\Omega^{'})$, one can figure this out by applying Lemma \ref{wh} with $w=u^{'}$ and replacing here and there $| u |_{3}$ by $G(u^{'})$ and $\parallel u \parallel_{3+r}$ by $\parallel u^{'} \parallel_{3+r,\Omega^{'}}$. As a consequence, all that is left to do is to estimate $b_{1h}(\bar{u}_h,v)$ and $b_{5h}(u_h,\Pi_h(v))$. \\

As for $b_{1h}$, to begin with we define $\Gamma_{{\mathcal R}} := \displaystyle \cup_{T \in {\mathcal R}_h \cap {\mathcal Q}_h} \Gamma_T$ and $\Gamma_{{\mathcal S}} := \Gamma \setminus \Gamma_{{\mathcal R}}$. Then we split $b_{1h}$ in the following fashion: 
\begin{equation}
\label{splitb1}
\left\{
\begin{array}{l}
b_{1h} = b^S_{1h}+b^R_{1h}, \mbox{ where}\\
b^R_{1h}(w,v) := \int_{\Gamma_{\mathcal R}} w \frac{\partial v}{\partial n},\\
b^S_{1h}(w,v) := \int_{\Gamma_{\mathcal S}} w \frac{\partial v}{\partial n}.
\end{array}
\right.
\end{equation}
$b^S_{1h}$ involves the sum of integrals on $\tilde{\partial} T$ for $T \in {\mathcal S}_h$ only, which can be estimated like in Theorem \ref{theorem1bis}. However in order to ensure sufficient differentiability, this is at the price of the enlargement of $\Gamma_{\mathcal S} \cap T_{\Delta}$ into $\Gamma \cap {\mathcal H}_T$ where ${\mathcal H}_T$ is the trihedral formed by the faces of $T \in {\mathcal S}_h$ in ${\mathcal F}_h$, and the natural extension of $u_h$ to ${\mathcal H}_T \cap \Omega$. The final result is qualitatively the same with a constant $C^S_{b1}$ similar to $C_{b1}$, namely,
\begin{equation}
\label{estimateb1S}
b_{1h}^S(\bar{u}_h,v) \leq C_{b1}^S h^3 [ h^{1/2} | u |_{3} + \parallel u \parallel_{3+r}] \parallel v \parallel_{2}.
\end{equation} 
In order to estimate $b_{1h}^R$ we first note that obviously enough it holds:
\begin{equation}
\label{estimb1R}
b^R_{1h}(\bar{u}_h,v) \leq C_t \displaystyle \left[ \sum_{T \in {\mathcal R}_h \cap {\mathcal Q}_h} 
\parallel \bar{u}_h \parallel_{0,\Gamma_T}^2 \right]^{1/2} \parallel v \parallel_{2}.
\end{equation}
According to the constructions previously advocated, if $T \in {\mathcal R}_h \cap {\mathcal Q}_h$ a subset of $\delta_e$ with a non-zero measure lies in the interior of $T$, where $e$ is the edge of $T$ contained in $\Gamma_h$. Moreover the underlying portion of the limiting curve $\gamma_e$ of $\delta_e$ is contained in $\Gamma_T$. Let $s$ be the curvilinear abscissa along $\gamma_e$ and $t$ be the abscissa along the intersections $\psi(s)$ of $\Gamma_T$ with the planes orthogonal to $\delta_e$ at the successive points along $\gamma_e \cap T$, in such a way that $\Gamma_T$ can be uniquely parametrized by $(s,t)$. Let $A(s)$ and $B(s)$ be the end-points of $\psi(s)$. In doing so, for a constant $C_q$ depending only on $\Gamma$, we trivially have, 
\begin{equation}
\label{estimb1R1}
\parallel \bar{u}_h \parallel_{0,\Gamma_T}^2 \leq C_q \displaystyle \int_{\gamma_e \cap T} \left[ \left|\int_{A(s)}^{B(s)} \bar{u}_h^2 dt \right| \right] ds.
\end{equation}
Now we observe that the length of $\psi(s)$ is bounded above by $C_{\psi} h_T^2$ $\forall s$, where $C_{\psi}$ is a constant independent of $T$. It follows that $|\bar{u}_h(Q) - \bar{u}_h(P)| \leq C_{\psi} h_T^2 \parallel {\bf grad} \; \bar{u}_h \parallel_{0,\infty,\tilde{T}}$ for $Q \in \gamma_e \cap T$ with abscissa $s$ and $\forall P \in \psi(s)$. Hence 
after straightforward calculations we can write,
\begin{equation}
\label{estimb1R2}
\parallel \bar{u}_h \parallel_{0,\Gamma_T}^2 \leq 2 C_q C_{\psi} h_T^2 \displaystyle \int_{\gamma_e} \left[\bar{u}_h^2 +  C_{\psi}^2 h_T^4 \parallel {\bf grad} \; \bar{u}_h \parallel_{0,\infty,\tilde{T}}^2 \right] ds.
\end{equation}
Now since $\bar{u}_h$ vanishes at three distinct points of $\gamma_e$ we can use a result in \cite{arXiv} according to which 
\begin{equation}
\label{estimb1R3}
\displaystyle \int_{\gamma_e} \bar{u}_h^2 ds \leq C_o h_T^4 \displaystyle 
\int_{\gamma_e} \left[ |H(\bar{u}_h)|^2 + | {\bf grad} \; \bar{u}_h |^2 \right],
\end{equation} 
where $C_o$ is a mesh-independent constant. \\
Noting that the length of $\gamma_e$ is an $O(h_T)$ from \eqref{estimb1R2} and \eqref{estimb1R3} we further obtain for another constant $C_n$ independent of $T$,
\begin{equation}
\label{estimb1R4}
\parallel \bar{u}_h \parallel_{0,\Gamma_T}^2 \leq C_n h_T^7 [ \parallel H(\bar{u}_h) \parallel_{0,\infty,\tilde{T}}^2 + \parallel {\bf grad} \; \bar{u}_h \parallel_{0,\infty,\tilde{T}}^2 ].
\end{equation} 
Now we resort to \eqref{Djwh} taking $T^{'} = \tilde{T}$ and $j=1,2$. In view of this \eqref{estimb1R4} easily leads to a mesh-independent constant $C_{\epsilon}$ for which it holds
\begin{equation}
\label{estimb1R5}
\parallel \! H(\bar{u}_h) \! \parallel_{0,\infty,\tilde{T}}^2 + \parallel \!~{\bf grad} \; \bar{u}_h \! \parallel_{0,\infty,\tilde{T}}^2  \leq C_{\epsilon} ( h_T^{-5} \parallel \! {\bf grad}\; \bar{u}_h \! \parallel_{0,\tilde{T}}^2 + h_T^{-1} | u |^2_{3,\tilde{T}} + h_T^{2\epsilon} \parallel \! u \! \parallel^2_{3+r,\tilde{T}} ).
\end{equation}
Collecting \eqref{estimb1R5} and \eqref{estimb1R4} and taking into account \eqref{estimateP2} and \eqref{estimb1R} we establish the existence of a mesh-independent constant $C_{b1}^R$ such that,
\begin{equation}
\label{estimateb1R}
b_{1h}^R(\bar{u}_h,v) \leq C_{b1}^R h^3 [ G(u^{'}) + h^{r} \| u \|_{3+r}] \parallel v \parallel_{2}.
\end{equation}
\indent In order to estimate $b_{5h}(u_h,v)$ we proceed as follows.\\
First we apply \eqref{Rolle} to obtain $|[\Pi_h(v)]({\bf x})| \leq C_{\Gamma} h_T^2 \parallel {\bf grad}\; \Pi_h(v) \parallel_{0,\infty,T}$ $\forall {\bf x} \in \Delta^{'}_T$ and $\forall T \in {\mathcal Q}_h$. Thus taking into account that the volume of $\Delta^{'}_T$ is bounded above by $C_{\Gamma} h_T^4$ for $T \in {\mathcal S}_h$ and by $C_{\Gamma}^2 h_T^5$ for $T \in {\mathcal R}_h$, by a straightforward argument we can write for a suitable constant $\bar{C}_{\Gamma}$ depending only on $\Gamma$:
\begin{equation}
\label{step1}
b_{5h}(u_h,\Pi_h(v)) \leq \displaystyle \sum_{T \in {\mathcal Q}_h} \bar{C}_{\Gamma} h_T^6 
\parallel \Delta u_h \parallel_{0,\infty,T} \parallel {\bf grad}\; \Pi_h(v) \parallel_{0,\infty,T}.
\end{equation}
Since all the components of $[{\bf grad}\;\Pi_h(v)]_{|T}$ and $[H(u_h)]_{|T}$ are in ${\mathcal P}_0(T)$ and those of $[{\bf grad}\;u_{h}]_{|T}$ are in ${\mathcal P}_1(T)$, in all the norms involving $\Pi_h(v)$ and $u_h$ appearing in \eqref{step1} $T$ can be replaced by $\tilde{T}$. Thus by \eqref{L2TDelta} and the Schwarz inequality we successively establish, 
\begin{equation}
\label{step2}
b_{5h}(u_h,\Pi_h(v)) \leq \bar{C}_{\Gamma} {\mathcal C}_J^2 \displaystyle \sum_{T \in {\mathcal Q}_h} 
 h_T^3 \parallel \Delta u_h \parallel_{0,\tilde{T}} \parallel {\bf grad}\; \Pi_h(v) \parallel_{0,\tilde{T}}, 
\end{equation}   
\begin{equation}
\label{step3}
b_{5h}(u_h,v_h) \leq \bar{C}_{\Gamma} {\mathcal C}_J^2 \displaystyle \left\{ \sum_{T \in {\mathcal Q}_h} h^6_T
\parallel \Delta u_h \parallel_{0,\tilde{T}}^2 \right\}^{1/2} \parallel {\bf grad}\; \Pi_h(v) \parallel_{\widetilde{0,h}}.
\end{equation}
Incidentally by the Sobolev Embedding Theorem $H^{3+r}(\Omega)$ is continuously embedded in $W^{2,\infty}(\Omega)$, which means that there exists a constant $C_e$ depending only on $\Omega$ such that,
\begin{equation}
\label{embedding}
\parallel u \parallel_{2,\infty} \leq C_e \parallel u \parallel_{3+r}.
\end{equation} 
Therefore replacing $u_h$ by $\bar{u}_h+u$ we can write $\forall T \in {\mathcal Q}_h$:
\begin{equation}
\label{step4}
\parallel \Delta u_h \parallel_{0,\tilde{T}}^2 
\leq 6 volume(T) \left( \parallel H(\bar{u}_h) \parallel_{0,\infty,\tilde{T}}^2 + 
C_e \parallel u \parallel_{3+r}^2 \right).
\end{equation}
Now we resort again to \eqref{Djwh} with $T^{'} = \tilde{T}$, $j=2$ and $w = u$. After straightforward calculations we can write for a mesh-independent constant $\bar{C}_5$,
\begin{equation}
\label{step4bis}
\parallel \Delta u_h \parallel_{0,\tilde{T}}^2 \displaystyle 
\leq \bar{C}_5 \left( h_T^{-2} \parallel {\bf grad} \; \bar{u}_h \parallel_{0,\tilde{T}}^2 + 
h_T^{2} | u |_{3,\tilde{T}}^2 + h_T^{2+2r} \parallel u \parallel_{3+r,\tilde{T}}^2 
+ h_T^{3} \parallel u \parallel_{3+r}^2 \right).
\end{equation}
Thus plugging \eqref{step4bis} into \eqref{step3} and using Theorem \ref{P2} together with \eqref{skinvolume} we easily obtain for another mesh-independent constant $\breve{C}_5$:
\begin{equation}
\label{step5}
b_{5h}(u_h,v_h) \leq  \breve{C}_5 h^3 \displaystyle \left[ h G(u^{'}) + (h^{1+r}+C_{\Gamma}^{'} h^{1/2}) {\color{blue} 
\parallel u \parallel_{3+r} } \right] {\color{blue} \parallel {\bf grad}\; \Pi_h(v) \parallel_{\widetilde{0,h}}}.
\end{equation}
On the other hand setting $\tilde{C}_{\Pi}=\sqrt{1+\tilde{C}_V^2 diam(\Omega)^2}$, \eqref{intildeP1} easily yields, 
\begin{equation}
\label{interpol}
\parallel {\bf grad}\; \Pi_h(v) \parallel_{\widetilde{0,h}} \leq \tilde{C}_{\Pi} \parallel v \parallel_{2}.
\end{equation} 
. Hence there exists a mesh-independent constant 
$C_{b5}$ such that,
\begin{equation}
\label{estimateb5}
b_{5h}(u_h,v_h) \leq  C_{b5} h^{3} \displaystyle \left\{ h G(u^{'}) 
+ h^{1/2} \parallel u \parallel_{3+r} \right\} \parallel v  \parallel_{2}.
\end{equation} 
Finally plugging into \eqref{L2estP2sex} the upper bounds \eqref{estildeah} and \eqref{estimateb4} together with \eqref{estimateb2}, \eqref{estimateb3}, \eqref{estimatec} and \eqref{estimated} with $b^{'}_{2h}$, $b^{'}_{3h}$, $c^{'}_h$ and $d^{'}_h$ instead of $b_{2h}$, $b_{3h}$, $c_h$ and $d_h$, and replacing $|u|_{k+1}$ by $G(u^{'})$ with $k=2$ and $\parallel u \parallel_{3+r}$ by $\parallel u^{'} \parallel_{3+r,\Omega^{'}}$ on the right hand side of all those inequalities, estimates 
\eqref{estimateb1S}, \eqref{estimateb1R} combined with \eqref{splitb1} together with \eqref{estimateb5} complete the proof. \QED \\ 


\section{Numerical experiments}

In this section we assess the accuracy of the method studied in Sections 3, 4, 5 - referred to hereafter as the \textit{new method} -, by solving equation (\ref{Poisson}) in some relevant test-cases, 
taking $k=2$. Comparisons with the isoparametric technique and the approach consisting of shifting boundary conditions from the true boundary to the boundary of the approximating polyhedron are also carried out. Hereafter the latter technique will be called the \textit{polyhedral approach}.
In all the examples numerical integration of the right hand side term was performed with the $15$-point Gauss quadrature formula given in \cite{Zienkiewicz}, with fourteen-digit accurate coefficients.  

\subsection{Consistency check}

In order to dissipate any skepticism about the performance of our method, we first solved the model problem with a constant right hand side equal to $2(a^{-2}+b^{-2}+1)$ in the ellipsoid centered 
at the origin given by the inequality $p(x,y,z) \leq 1$ where $p(x,y,z)=(x/a)^2+(y/b)^2+z^2$. Taking $g \equiv 0$, the exact solution is the quadratic function $1-p$, and 
thus the new method is expected to reproduce it up to machine precision for any mesh. i.e., except for round-off errors. Incidentally we observe that the isoparametric version of the finite element method does not enjoy the same property. Hence from this pont of view it is not a consistent method, for it can only reproduce exactly linear functions (up to machine precision).  \\
Here we used a mesh consisting of $3072$ tetrahedra resulting from the transformation of a standard uniform $6 \times 8 \times 8 \times 8$ mesh of a unit cube $\Omega_0$ into tetrahedra having one edge coincident with a diagonal parallel to the line $x=y=z$ of a cube with edge length equal to $1/8$, resulting from a first subdivision of $\Omega_0$ into $8^3$ equal cubes. The final tetrahedral mesh of the ellipsoid octant corresponding to positive values of $x,y,z$, contains the same number of elements and is generated by mapping the unit cube into the latter domain through the transformation of Cartesian coordinates into spherical coordinates using a procedure described in \cite{RBC}. \\
It turns out that the error in the broken $H^1$-semi-norm $\parallel {\bf grad}(\cdot) \parallel_{0,h}$ resulting from computations with $a=0.6$ and $b=0.8$, equals approximately $0.10599965 \times 10^{-13}$, for an exact value of ca. $1.0214597$. From these computations done in double precision  the numerical solution can be considered to be exact, taking into account the precision of the numerical integration coefficients. At this point we emphasize that our nonconforming method is fully algebraic consistent in the sense of \cite{Veeser}, in contrast to the isoparametic technique. Indeed it enjoys the property of reproducing exactly conforming piecewise polynomial solutions which are locally of degree $\leq k$.\\
It is noteworthy that the absolute error measured in the same way for the polyhedral approach is about $0.01663104$, while it equals ca. $0.01001501$ if the isoparametric technique is employed with the same degree of mesh refinement, as seen in Subsection 6.4. This means relative errors of about $1.6$ percent and $1.0$ percent, respectively. One might object that this is not so bad for a rather coarse mesh. However substantial gains with the new method over the polyhedral or the isoparametric approach will be manifest in the examples that follow. 
           
\subsection{Test-problems in a convex domain} 
We next validate error estimates (\ref{errestconvex}) and \eqref{L2estconvex} by assessing method's accuracy in $\Omega_h$. With this aim we solved two test-problems with known exact solution. Corresponding results are reported below.\\

\noindent \underline{Test-problem 1:} Here $\Omega$ is the unit sphere centered at the origin. We take the exact solution $u = \rho^2-\rho^4$ where $\rho^2 
= x^2 + y^2 +z^2$, which means that $g \equiv 0$ and $f=-6+20 \rho^2$. Owing to symmetry we consider only the octant sub-domain given by $x>0$, $y>0$ and $z>0$ by prescribing Neumann boundary conditions on $x=0$, $y=0$ and $z=0$.  
We computed with quasi-uniform meshes defined by a single integer parameter $J$, constructed by the procedure proposed in \cite{RBC} and described in main lines at the beginning of this section. Roughly speaking the mesh of the computational sub-domain is the spherical-coordinate counterpart of the uniform partition of the unit cube $(0,1) \times (0,1) \times (0,1)$ into $J^3$ identical cubic cells. Each element of the final mesh is the transformation of a tetrahedron out of six resulting from the subdivision of each cubic cell; the latter have as an edge the cell's diagonal parallel to the line $x=y=z$. Since both the mesh and the solution are symmetric with respect to the three Cartesian axes computations were effectively performed only for a third of the chosen octant sub-domain. \\
\indent In Table 1 we display the absolute errors in the norms $\parallel {\bf grad}(\cdot) \parallel_{0,h}$ and $\parallel \cdot \parallel_{0,h}$ 
for increasing values of $J$, namely, $J=4,8,12,16,20$. Since the true value of $h$ equals $\kappa/J$ for a suitable constant $\kappa$, as a reference we set $h = 1/J$ to simplify things.      
\noindent As one infers from Table 1, the approximations obtained with the new method perfectly conform to the theoretical estimates (\ref{errestconvex}) and \eqref{L2estconvex}. Indeed as $J$ increases the errors in the 
broken $H^1$-semi-norm decrease roughly like $h^2$ as predicted. The error in the $L^2$-norm in turn tends to decrease as an $O(h^3)$.
 In Table 2 we display the same kind of results obtained with the polyhedral approach. As one can observe the error in the broken $H^1$-semi-norm decreases roughly 
like $h^{1.5}$, as predicted by the mathematical theory of the finite element method, while the errors in the $L^2$-norm seem to behave like an $O(h^2)$.\\

\begin{table*}[h!]
{\small 
\centering
\begin{tabular}{ccccccc} &\\ [-.3cm]  
$h$ & $\longrightarrow$ & 1/4 & $1/8$ & $1/12$ & $1/16$ & $1/20$ 
\tabularnewline & \\ [-.3cm] \hline &\\ [-.3cm]
$\parallel {\bf grad}_h(u-u_h) \parallel_{0,h}$ & $\longrightarrow$ & 0.187649 E-1 & 0.499091 E-2 & 0.225836 E-2 & 0.128114 E-2 & 0.823972 E-3   
\tabularnewline &\\ [-.3cm] \hline &\\ [-.3cm] 
$\parallel u-u_h \parallel_{0,h}$ & $\longrightarrow$ & 0.653073 E-3 & 0.845686 E-4 & 0.253348 E-4 & 0.107516 E-4 & 0.552583 E-5   
\tabularnewline &\\ [-.3cm] \hline &\\ [-.3cm]
\end{tabular}
\caption{Errors with the new method measured in two different manners for Test-problem 1.} 
}
\label{table5}
\end{table*}

\begin{table*}[h!]
{\small 
\centering
\begin{tabular}{ccccccc} &\\ [-.3cm]  
$h$ & $\longrightarrow$ & 1/4 & $1/8$ & $1/12$ & $1/16$ & $1/20$ 
\tabularnewline & \\ [-.3cm] \hline &\\ [-.3cm]
$\parallel {\bf grad}_h(u-u_h) \parallel_{0,h}$ & $\longrightarrow$ & 0.257134 E-1 & 0.917910 E-2 & 0.50152682 E-2 & 0.326410 E-2 & 0.233854 E-2  
\tabularnewline &\\ [-.3cm] \hline &\\ [-.3cm] 
$\parallel u-u_h \parallel_{0,h}$ & $\longrightarrow$ & 0.454733 E-2 & 0.113568E-2 & 0.502166 E-3 & 0.281468 E-3 & 0.179698 E-3   
\tabularnewline &\\ [-.3cm] \hline &\\ [-.3cm]
\end{tabular}
\caption{Errors for the polyhedral approach measured in two different manners for Test-problem 1.} 
}
\label{table5s}
\end{table*}

\noindent \underline{Test-problem 2:} In order to make sure that the previous example has no particularity due to the simple form of the domain, we now consider  
$\Omega$ to be the ellipsoid centered at the origin with semi-axes $a$, $b$ and $1$.  
We take $g \equiv 0$ and $f = -\Delta u$ for the exact solution $u = [1-(x/a)^2-(y/b)^2-z^2][1-(x/b)^2-(y/a)^2-z^2]$.  
In view of the symmetry with respect to the planes $x=0$, $y=0$ and $z=0$, computations are restricted to the octant sub-domain given by $x>0$, $y>0$ and $z>0$, by prescribing Neumann boundary conditions on $x=0$, $y=0$ and $z=0$. We computed with quasi-uniform meshes defined by a single integer parameter $J$, constructed in a way in all analogous to the procedure described in Test-problem 1, i.e. the one proposed in \cite{RBC} for spheroidal domains.
Like in the case of the ellipsoid considered at the beginning of this section, this means that the mesh of the computational sub-domain is a spherical-coordinate counterpart of the $6(J \times J \times J)$ uniform mesh of the unit cube $(0,1) \times (0,1) \times (0,1)$. \\
\indent Taking again $a=0.6$ and $b=0.8$, we display in Table 3 the errors in the norms $\parallel {\bf grad}(\cdot) \parallel_{0,h}$ and $\parallel \cdot \parallel_{0,h}$, 
for increasing values of $J$, namely, $J=2,4,8,12$, for the new method and the polyhedral approach, respectively. For simplicity we quite abusively set again $h = 1/J$.      
\noindent As one infers from Table 3, akin to Test-problem 1, the approximations obtained with the new method are 
also in full agreement with the theoretical estimates \eqref{errestconvex} and \eqref{L2estconvex}. Indeed as $J$ increases the errors in the $L^2$-norm of error function's broken gradient decrease roughly as $(1/J)^2$, as predicted. Moreover here again, the error in the $L^2$-norm behaves roughly like an $O(h^3)$. On the other hand Table 4 certifies again the losses in order for the 
polyhedral approach, close to those observed for Test-problem 1.             
 
\begin{table*}[h!]
{\small 
\centering
\begin{tabular}{cccccc} &\\ [-.3cm]  
$h$ & $\longrightarrow$ & $1/2$ & $1/4$ & $1/8$ & $1/12$  
\tabularnewline & \\ [-.3cm] \hline &\\ [-.3cm]
$\parallel {\bf grad}_h(u-u_h) \parallel_{0,h}$ & $\longrightarrow$ & 0.117716 E+0 & 0.353096 E-1 & 0.943753 E-2 & 0.427408 E-2    
\tabularnewline &\\ [-.3cm] \hline &\\ [-.3cm] 
$\parallel u-u_h \parallel_{0,h}$ & $\longrightarrow$ & 0.705684 E-2 & 0.956478 E-3 & 0.122026 E-3 & 0.364375 E-4   
\tabularnewline &\\ [-.3cm] \hline &\\ [-.3cm]
\end{tabular}
\caption{Errors with the new method measured in two different manners for Test-problem 2.} 
}
\label{table6}
\end{table*}

\begin{table*}[h!]
{\small 
\centering
\begin{tabular}{cccccc} &\\ [-.3cm]  
$h$ & $\longrightarrow$ & $1/2$ & $1/4$ & $1/8$ & $1/12$  
\tabularnewline & \\ [-.3cm] \hline &\\ [-.3cm]
$\parallel {\bf grad}_h(u-u_h) \parallel_{0,h}$ & $\longrightarrow$ & 0.124723 E+0 & 0.368763 E-1 & 0.104133 E-1 & 0.501084 E-2    
\tabularnewline &\\ [-.3cm] \hline &\\ [-.3cm] 
$\parallel u-u_h \parallel_{0,h}$ & $\longrightarrow$ & 0.807272 E-2 & 0.163738 E-2 & 0.365620 E-3 & 0.157317 E-3   
\tabularnewline &\\ [-.3cm] \hline &\\ [-.3cm]
\end{tabular}
\caption{ Errors for the polyhedral approach measured in two different manners for Test-problem 2.} 
}
\label{table6s}
\end{table*}

\subsection{Test-problem in a non-convex domain}

\noindent \underline{Test-problem 3:} The aim of the following test-problem is to assess the behavior of the new method when $\Omega$ is not convex, taking now a non-polynomial exact solution. 
More precisely (\ref{Poisson}) is solved in the torus $\Omega$ with minor radius $r_m$ and major radius $r_M$. This means that the torus' inner radius $r_i$ equals $r_M-r_m$ and its outer radius $r_e$ equals $r_M+r_m$. Hence $\Gamma$ is given by the equation $(r_M-\sqrt{x^2+y^2})^2 + z^2 = r_m^2$. We only consider problems with symmetry about the $z$-axis, and with respect to the plane $z=0$. For this reason we may work with a computational domain given by 
$\{(x,y,z) \in \Omega \; | \; z \geq 0; \; 0 \leq \theta \leq \pi/4 \mbox{ with } \theta = atan(y/x)\}$. A family of meshes of this domain depending on a single even integer parameter $I$ containing $6 I^3$ tetrahedra is generated by the following procedure. First we generate a partition of the cube $(0,1) \times (0,1) \times (0,1)$ into $I^3/2$ equal rectangular boxes by subdividing the edges parallel to the $x$-axis, the $y$-axis and the $z$-axis into $2I$, $I/2$ and $I/2$ equal segments, respectively. Then each box is subdivided into six tetrahedra having an edge parallel to the line $4x=y=z$. This mesh with $3I^3$ tetrahedra is transformed into the mesh of the quarter cylinder $\{(x,y,z) \; | \; 0 \leq x \leq 1, \; y \geq 0, \; z \geq 0, \; y^2 + z^2 \leq 1 \}$, following the transformation of the mesh consisting of $I^2/2$ equal right triangles formed by the faces of the mesh elements contained in the unit cube's section given by $x=j/(2I)$, for $j=0,1,\ldots,2I$. 
The latter transformation is based on the mapping of the Cartesian coordinates $(y,z)$ into the 
polar coordinates $(r,\varphi)$ with $r=\sqrt{y^2+z^2}$, using a procedure of the same nature as the one described in \cite{RBC} (cf. Figure 4). Then the resulting mesh of the quarter cylinder is transformed into the mesh with $6I^3$ thetrahedrons of the half cylinder $\{(x,y,z) \; | \; 0 \leq x \leq 1, \; -1 \leq y \leq 1, \; z \geq 0, \; y^2 + z^2 \leq 1 \}$ by symmetry with respect to the plane $y=0$. Finally this mesh is transformed into the computational mesh (of an eighth of half-torus) by first 
mapping the Cartesian coordinates $(x,y)$ into polar coordinates $(\rho,\theta)$, with $\rho = r_M + y r_m$ and $\theta = x \pi/4$, and then the latter coordinates into new Cartesian coordinates $(x,y)$ using the relations $x = \rho cos \theta$ and $y = \rho sin \theta$. Notice that the faces of the final tetrahedral mesh on the sections of the torus given by $\theta = j \pi/(8I)$, for $j=0,1,\ldots,2I$, form a triangular mesh of a disk with radius equal to $r_m$, having the pattern  illustrated in Figure 4 for a quarter disk, taking $I=4$, $\theta=0$ and $r_m=1$ (cf. \cite{RBC}).                     
  
\begin{figure}[h]
\label{fig4}
\begin{center}
\includegraphics[scale=0.43]{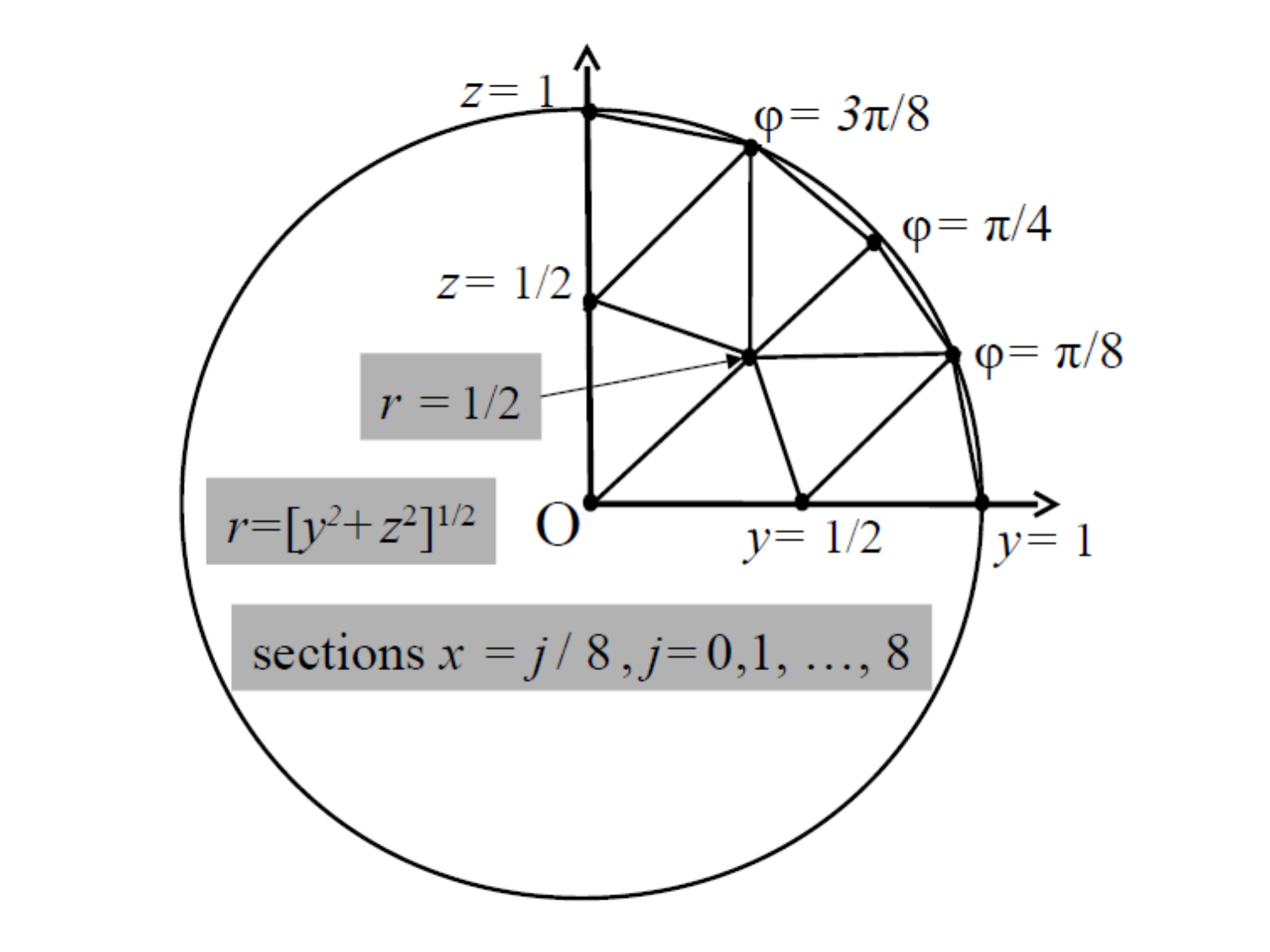}
\end{center}
\par
\caption{Trace of the intermediate mesh of 1/4 cylinder on sections $x=j/(2I)$, $0 \leq j \leq 2I$, for $I=4$}
\end{figure} 

Recalling that here $\rho=\sqrt{x^2+y^2}$, we take $r_M=5/6$, $r_m=1/6$ and $f^{'}=6-5/(3\rho)$. For $g \equiv 0$ the exact solution is given by 
$u = 1/36-z^2-(5/6-\rho)^2$. Obviously enough we take the same expression for $u^{'}$.\\ 
\indent In Table 5 we display the errors in the norm $\parallel {\bf grad}(\cdot) \parallel_{0,h}$ and in the norm of $L^2(\Omega_h)$,  
for increasing values of $I$, namely $I=2^m$ for $m=1,2,3,4$. 
Now we take as a reference $h=\pi/(8I)$. \\ 
\noindent As one can observe from Table 5, here again the quality of the approximations obtained with the new method is in very good agreement with the theoretical result 
(\ref{errorestimate}), for as $I$ increases the errors in the broken $H^1$-semi-norm decrease roughly as $1/I^2$  as predicted. On the other hand here again the errors in the $L^2$-norm are in agreement with \eqref{L2P2} for they decrease roughly like $1/I^3$. Table 6 in turn shows a qualitative erosion 
of the solution errors obtained by means of the polyhedral approach similar to the case of convex domains.       
           
\begin{table*}[h!]
{\small 
\centering
\begin{tabular}{cccccc} &\\ [-.3cm]  
$h$ & $\longrightarrow$ & $\pi/32$ & $\pi/64$ & $\pi/128$ & $\pi/256$ 
\tabularnewline & \\ [-.3cm] \hline &\\ [-.3cm]
$\parallel {\bf grad}_h(u^{'}-u_h) \parallel_{0,h}$ & $\longrightarrow$ & 0.786085 E-3 & 0.205622 E-3 &  0.522963 E-4 & 0.131844 E-4    
\tabularnewline &\\ [-.3cm] \hline &\\ [-.3cm] 
$\parallel u^{'}-u_h \parallel_{0,h}$ & $\longrightarrow$ & 0.133794 E-4 & 0.171222 E-5 & 0.214555 E-6  & 0.269187 E-7  
\tabularnewline &\\ [-.3cm] \hline &\\ [-.3cm]
\end{tabular}
\caption{Errors with the new method measured in two different manners for Test-problem 3.}
}
\label{table7}
\end{table*}

\begin{table*}[h!]
{\small 
\centering
\begin{tabular}{cccccc} &\\ [-.3cm]  
$h$ & $\longrightarrow$ & $\pi/32$ & $\pi/64$ & $\pi/128$ & $\pi/256$  
\tabularnewline & \\ [-.3cm] \hline &\\ [-.3cm]
$\parallel {\bf grad}_h(u^{'}-u_h) \parallel_{0,h}$ & $\longrightarrow$ &  0.829181 E-2 & 0.327176 E-2 & 0.119077 E-2 & 0.425739 E-3   
\tabularnewline &\\ [-.3cm] \hline &\\ [-.3cm] 
$\parallel u^{'}-u_h \parallel_{0,h}$ & $\longrightarrow$ & 0.579150 E-3 & 0.143425 E-3 & 0.343823 E-4 & 0.834136 E-5    
\tabularnewline &\\ [-.3cm] \hline &\\ [-.3cm]
\end{tabular}
\caption{ Errors for the polyhedral approach measured in two different manners for Test-problem 3.}
}
\label{table7s}
\end{table*}

\subsection{Comparison with the isoparametric technique}

The results in Subsections 5.1, 5.2 and 5.3 validate the finite-element methodology studied in this article in the three-dimensional case. A priori it is an advantageous alternative in many respects to more classical techniques such as the isoparametric version of the finite element method. This is because its most outstanding features are not only universality but also simplicity, and eventually accuracy and CPU time too, although the two latter aspects were not our point from the beginning. 
Nevertheless we have compared our technique with the isoparametric one in terms of accuracy, by solving with 
both methods for $k=2$ the Poisson equation in the same domain as in Test-problem 2 and for the same exact solution. \\
Here again, owing to symmetry, we considered only the octant domain given by $x>0$, $y>0$ and $z>0$ by prescribing Neumann boundary conditions on $x=0$, $y=0$ and $z=0$.\\
\indent We supply in Table 7 the $L^2$-norms of the gradient of the error function and of this function itself, and maximum error at the nodes of the mesh, that is, a pseudo-$L^{\infty}$-seminorm that we denote by $\parallel \cdot \parallel_{0,\infty,h}$. The isoparametric solution is denoted by $\tilde{u}_h$. On the other hand the subscript $0,\tilde{h}$ replaces $0,h$ in the $L^2$-norms for the isoparametric case, in order to signify that the integrations take place in a curved domain approximating $\Omega$ instead of $\Omega_h$.    
We took again $a=0.6$, $b=0.8$ and computed with the same kind of meshes defined by a single integer parameter $J$ as for Test-problem 2. \\From Table 7 one can observe that both methods are of the same order as expected. However the new method was more accurate than  isoparametric elements all the way, especially in terms of nodal values. \\
On the other hand we report that both methods are roughly equivalent in terms of CPU time.
 
\begin{table*}[h!]
{\small 
\centering
\begin{tabular}{ccccccc} &\\ [-.3cm]  
$h$ & $\longrightarrow$ & $1/2$ & $1/4$ & $1/8$ & $1/12$ & $1/16$ 
\tabularnewline & \\ [-.3cm] \hline &\\ [-.3cm]
$\parallel {\bf grad}(u-u_h) \parallel_{0,h}$ & $\longrightarrow$ & 0.117716 E+0 & 0.353096 E-1 & 0.943754 E-2 & 0.427408 E-2 & 0.242528 E-2   
\tabularnewline &\\ [-.3cm] \hline &\\ [-.3cm] 
$\parallel {\bf grad}(u-\tilde{u}_h) \parallel_{0,\tilde{h}}$ & $\longrightarrow$ &  0.139311 E+0 & 0.390893 E-1 & 0.100150 E-1 &  0.445839 E-2 & 0.250611 E-2 
\tabularnewline &\\ [-.3cm] \hline &\\ [-.3cm]
\tabularnewline &\\ [-.3cm] \hline &\\ [-.3cm]
$\parallel u-u_h \parallel_{0,h}$ & $\longrightarrow$ & 0.705684 E-2 & 0.956478 E-3 & 0.122026 E-3 & 0.364375 E-4 & 0.154448 E-4 
\tabularnewline &\\ [-.3cm] \hline &\\ [-.3cm] 
$\parallel u-\tilde{u}_h \parallel_{0,\tilde{h}}$ & $\longrightarrow$ &  0.752197 E-2 & 0.105638 E-2  & 0.131730 E-3 &  0.386297E-4 & 0.161845 E-4
\tabularnewline &\\ [-.3cm] \hline &\\ [-.3cm]
\tabularnewline &\\ [-.3cm] \hline &\\ [-.3cm]
$\parallel u-u_h \parallel_{0,\infty,h}$ & $\longrightarrow$ & 0.360639 E-1 & 0.693934 E-2 & 0.106156 E-2 & 
0.331707 E-3 & 0.143288 E-3  
\tabularnewline & \\ [-.3cm] \hline &\\ [-.3cm]
$\parallel u-\tilde{u}_h \parallel_{0,\infty,h}$ & $\longrightarrow$ & 0.409800 E-1 & 0.791483 E-2  & 0.123837 
E-2 & 0.389720 E-3 & 0.168969 E-3  
\tabularnewline &\\ [-.3cm] \hline &\\ [-.3cm]
\end{tabular} 
\caption{Errors with the new and the isoparametric approach for Test-problem 1 and $k=2$.} 
}
\label{table8}
\end{table*} 

\section{A nonconforming method with mean-value degrees of freedom}

Our technique to handle Dirichlet conditions prescribed on curved boundaries has a wide scope of applicability. The aim of this section is to illustrate this assertion once more, in the case of a nonconforming method with degrees of freedom other than function nodal values.\\
Incidentally for many well-known nonconforming finite element methods the construction of an isoparametric counterpart brings no improvement. This does not prevent suitable parametric elements from being successfully employed in this case. However to the best of author's knowledge studies in this direction are incipient. This fact motivates us to show in this section that our technique for handling curvilinear boundaries can be optimally extended in a straightforward manner to finite element methods, which are nonconforming even in the case of polytopes.\\ 
The method to be studied here is based on the same type of piecewise quadratic interpolation as the one introduced in \cite{JJAM}, in order to optimally represent the velocity in the framework of the stable solution of incompressible viscous flow problems. Actually the corresponding velocity representation enriched by the quartic bubble-functions of the tetrahedra combined with a discontinuous piecewise linear pressure in each tetrahedron is a sort of nonconforming three-dimensional analog of the popular conforming Crouzeix-Raviart mixed finite element \cite{CrouzeixRaviart} for solving the same kind of flow problems in two-dimension space. Here we use such a nonconforming approach to solve the model problem \eqref{Poisson}. With this aim we confine ourselves to the case of homogeneous boundary conditions for the sake of simplicity, though without any loss of essential aspects.\\

To begin with we recall the space $V_h^{*}$ of test-functions defined in $\Omega_h$, associated with the method under consideration. \\
Generically denoting by $F$ and $e$ a face and an edge of a tetrahedron 
$T \in {\mathcal T}_h$ respectively, by $A_e$ and $B_e$ the end-points of $e$ and by $M_e$ the mid-point of 
$e$, any function $v \in V_h^{*}$ restricted to every $T$ is a polynomial of degree less than or equal to two, defined upon the following set of degrees degrees of freedom:
\begin{itemize}
\item The four values $\mu_F(v)$ of $v$ at the centroids of $F$;
\item The six mean values $\nu_e(v)$ along $e$, where $\nu_e(v)=0.4 v(M_e) + 0.3[v(A_e)+v(B_e)]$.
\end{itemize} 
$\forall v \in V_h^{*}$ and $\forall F$ and $e$, we require that both $\mu_F(v)$ and $\nu_e(v)$ coincide for all tetrahedra of the mesh sharing the face $F$ or the edge $e$; moreover we require that both $\mu_F(v)$ and 
$\nu_e(v)$ vanish whenever $F$ or $e$ is contained in $\Gamma_h$. Clearly enough these requirements are not sufficient to ensure the continuity in $\Omega_h$ of a function in $V_h^{*}$, and hence this space is not a subspace of $H^1_0(\Omega_h)$. \\
The set of local canonical quadratic basis functions in a tetrahedron $T \in {\mathcal T}_h$ associated with the above degrees of freedom can be found in \cite{JJAM}. It is noteworthy that the gradients of all of them are an $O(h_T^{-1})$. This is a key property for the proof of Lemma \ref{lemma4} hereafter. \\

Similarly to the case of the standard Lagrangian piecewise quadratic elements, we define the trial-function space 
$W^{*}_h$ in the same way as $V_h^{*}$, except for the fact that the 
degrees of freedom associated with faces $F$ and edges $e$ contained in $\Gamma_h$ are modified as follows:
For a given function $w \in W^{*}_h$, $\mu_F(w)$ is replaced by $\tilde{\mu}_F(w)$ defined to be the value of $w$ 
at the point $P$ lying in the nearest intersection with $\Gamma$ of the perpendicular to $F$ passing through the centroid of $F$; referring to Figure 3, $\nu_e(w)$ is replaced by $\tilde{\nu}_e(w):=0.4 w(Q_e) + 0.3[w(A_e) + 
w(B_e)]$, where $Q_e$ is the nearest intersection with $\Gamma$ of the perpendicular to $e$ in $\delta_e$ passing through $M_e$. $\forall w \in W^{*}_h$ we require that both $\tilde{\mu}_F(w)$ and $\tilde{\nu}_e(w)$ vanish for every face $F$ or edge $e$ contained in $\Gamma_h$. \\
Similarly to Lemmata \ref{lemma2} and \ref{lemma3}, we have 
\begin{lemma}
\label{lemma4} 
Provided $h$ is small enough, $\forall T \in {\mathcal S}_h$ (resp. $\forall T \in {\mathcal R}_h$), given a set of $m$ real values generically denoted by $b$, with $m=6$ (resp. $m=9$), there exists a unique function $w_T \in {\mathcal P}_2(T)$ such that $\tilde{\mu}_F(w_T)=0$ and $\tilde{\nu}_e(w_T)=0$ if $F$ and $e$ are a face or an edge of $T$ contained in $\Gamma_h$, and such that $\mu_F(w_T)$ and $\nu_e(w_T)$ take the assigned value 
$b$, if neither $F$ nor $e$ is a face or an edge of $T$ contained in $\Gamma_h$.   
\end{lemma}   

\prov
The proof of this result goes very much like the one of Lemma \ref{lemma2}. This is essentially because the 
absolute value of the difference between both $\mu_F(w_T)$ and $\tilde{\mu}_F(w_T)$, and $\nu_e(w_T)$ and $\tilde{\nu}_e(w_T)$ is bounded above by $C_{\Gamma} h_T^2 
\parallel {\bf grad} \; w_T \parallel_{0,\infty,T}$, for every face $F$ or edge $e$ contained in $\Gamma_h$. \QED \\

Lemma \ref{lemma4} allows us to assert that $W_h^{*}$ is indeed a nonempty function space, whose dimension equals the one of $V_h^{*}$. Moreover if $u^{*}$ is a function in $H^2(\Omega) \cap H^1_0(\Omega)$, we can define $I^{*}_h(u^{*}) \in W_h^{*}$ to be the function given by $\mu_F(I^{*}_h(u^{*}))=\mu_F(u^{*})$ and $\nu_e(I^{*}_h(u^{*}))=\nu_e(u^{*})$ for all the faces $F$ and edges $e$ of tetrahedra in ${\mathcal T}_h$ not contained in $\Gamma_h$. Akin to the operator $I_h$, from standard interpolation results it is easy to see that $I^{*}_h$ enjoys the following property: \\
There exists a mesh-independent constant $C_P$ such that $\forall u^{*} \in H^3(\Omega) \cap H^1_0(\Omega)$ it holds, 
\begin{equation}
\label{interperror*}
 \parallel {\bf grad}_h (u^{*} - I^{*}_h(u^{*})) \parallel_{\widetilde{0,h}} \leq C_P h^2 | u^{*} |_{3}.
\end{equation} 
It also follows that the following problem can be considered to approximate \eqref{Poisson}:
\begin{equation}
\label{Poisson-h*}
\left\{
\begin{array}{l}
\mbox{Find } u^{*}_h \in W^{*}_h \mbox{ such that } a^{*}_h(u^{*}_h,v) = L_h(v) \; \forall v \in  V_h^{*}, \\
\\
\mbox{where}\\
\\
a_h^{*}(w,v) := \int_{\Omega_h} {\bf grad}_h w \cdot {\bf grad}_h v, \mbox{ for } w \in W^{*}_h + H^1(\Omega_h), \; v \in V^{*}_h.
\end{array}
\right.
\end{equation}
and $L_h$ was defined in \eqref{Poisson-h} with $f \equiv 0$ in $\Omega_h \setminus \Omega$. \\

The well-posedness of problem \eqref{Poisson-h*} is a direct consequence of the following propositions analogous to Propositions \ref{propo2} and \ref{propo3}: 

\begin{e-proposition}
\label{propo4}
If $h$ is sufficiently small there exists a constant $\alpha^{*} > 0$ independent of 
$h$ such that,
\begin{equation}
\label{inf-sup*}
\forall w \in W_h^{*} \neq 0, \displaystyle \sup_{v \in V_h^{*} \setminus \{ 0 \}} \frac{a^{*}_h(w,v)}{\parallel {\bf grad}_h w \parallel_{0,h} \parallel {\bf grad} \; v \parallel_{0,h}} 
\geq \alpha^{*}. 
\end{equation}  
\end{e-proposition}

\prov In order to prove this result, for a given $w \in W_h^{*}$ we construct $v \in V^{*}_h$ in such a way that 
$\mu_F(v) = \mu_F(w)$ and $\nu_e(v) = \nu_e(w)$ for all faces $F$ and edges $e$ of the mesh not contained in 
$\Gamma_h$. Then proceeding exactly like in the proof of Proposition \ref{propo2},  
\eqref{inf-sup*} is thus established. \QED \\

\begin{e-proposition}
\label{propo5}
Provided $h$ is sufficiently small, problem (\ref{Poisson-h*}) has a unique solution.  
\end{e-proposition}

\prov Clearly enough it holds 
$$a^{*}_h(w,v) \leq \parallel {\bf grad}_h w \parallel_{0,h} \parallel {\bf grad}_h v \parallel_{0,h} \; 
\forall (w,v) \in W^{*}_h \times V^{*}_h.$$ 
On the other hand according to \cite{JJAM} $\parallel {\bf grad}_h (\cdot) \parallel_{0,h}$ is a norm of $V^{*}_h$, certainly equivalent to the norm of $L^2(\Omega_h)$. Therefore $L_h$ is a continuous linear form over $V^{*}_h$. Hence resorting to the well-known theory of weakly coercive linear variational problems (cf. \cite{Babuska}, \cite{Brezzi} and \cite{COAM}), the result directly follows from Proposition \ref{propo4}. \QED \\

Next we establish error estimates for problem (\ref{Poisson-h*}). Here again we distinguish the convex case from the non-convex case. \\
First we have:

\begin{theorem}
\label{convexestim}
Assume that $f \in H^{1}(\Omega)$ and $g \equiv 0$. As long as $h$ is sufficiently small, if $\Omega$ is a convex domain smooth enough for the solution $u$ of (\ref{Poisson}) to belong to $H^{3}(\Omega)$, there exists a constant $C^{*}(f)$ depending only on $f$ such that the solution $u^{*}_h$ of (\ref{Poisson-h*}) satisfies :
\begin{equation}
\label{errestconvex*}
\parallel {\bf grad}_h(u - u^{*}_h) \parallel_{0,h} \leq C^{*}(f) h^2. 
\end{equation}
\end{theorem}

\prov According to \cite{COAM}, using Proposition \ref{propo4} we can write:
\begin{equation}
\label{bound*}
 \parallel {\bf grad}_h(u - u^{*}_h) \parallel_{0,h}\leq \displaystyle \frac{1}{\alpha^{*}} \left[ 
\parallel {\bf grad}_h(u - I^{*}_h(u)) \parallel_{0,h} + \sup_{v \in V^{*}_h \setminus \{0\}} 
\frac{|a^{*}_h(u,v)-L_h(v)|}{\parallel {\bf grad}_h v \parallel_{0,h}} \right].
\end{equation}

\prov Taking into account \eqref{interperror*}, all we have to do is to estimate the sup term on
the right hand side of \eqref{bound*}. But this is a matter that was already addressed in \cite{JJAM}. More precisely the required estimate is a consequence of the fact that the $L^2$-projection of the trace on a face $F$ of the mesh of any function $v \in V^{*}_h$ onto the space ${\mathcal P}_1(F)$, is a linear combination 
of the values $\mu_F(v)$ and $\nu_e(v)$, where $e$ here generically represents the edges of $F$. 
Actually this property implies the existence of a mesh-independent constant $C_R$ such that,
\begin{equation}
\label{estimaresidual}
|a^{*}_h(u,v)-L_h(v)| \leq C_R h^2 | u |_{3} \parallel {\bf grad}_h v \parallel_{0,h}.
\end{equation}
Then \eqref{errestconvex*} directly follows from \eqref{bound*}, \eqref{interperror*} and \eqref{estimaresidual}. \QED \\ 

\begin{theorem} 
\label{concavestim}
Assume that $u \in H^3(\Omega)$. Provided $h$ is sufficiently small, there exists a mesh-independent constant 
$\tilde{C}^{*}$ such that the unique solution $u^{*}_h$ to (\ref{Poisson-h*}) satisfies:
\begin{equation}
\label{errestconcave*} 
\begin{array}{l}
\parallel {\bf grad}_h(u - u_h^{*}) \parallel_{\widetilde{0,h}} \leq \tilde{C}^{*} h^{2} \parallel u^{'} \parallel_{3,\Omega^{'}}, \\
\end{array}
\end{equation}
$u^{'} \in H^3(\Omega^{'})$ being the regular extension of $u$ to $\Omega^{'}$ constructed in accordance to Stein et al. \cite{Stein}.
\end{theorem}

\prov
First of all combining \eqref{Poisson-h*} with 
Proposition \ref{propo4} we can write:
\begin{equation}
\label{concave1}
\parallel {\bf grad}_h(u^{*}_h - I^{*}_h(u)) \parallel_{0,h} \leq \displaystyle \frac{1}{\alpha^{*}} \displaystyle \sup_{v \in V^{*}_h \setminus \{0\}} 
\frac{|a^{*}_h(u^{'},v)-L_h(v)| + |a^{*}_h(u^{'}-I^{*}_h(u),v)|}{\parallel {\bf grad}\; v \parallel_{0,h}}.
\end{equation}
The first term in the numerator of (\ref{concave1}) can be estimated in the following manner. \\
Following the same steps as in Theorem \ref{P2}, and recalling the subset ${\mathcal Q}_h$ of ${\mathcal O}_h$ together with the subset $\Delta^{'}_T$ of $T \in {\mathcal Q}_h$ defined therein, we apply First Green's identity to $a^{*}_h(u^{'},v)$. Noticing that $v$ is not continuous across the inter-element boundaries, and recalling the notation $\partial T$ for the boundary of $T \in {\mathcal T}_h$ and $\partial (\cdot)/\partial n_T$ for the normal derivative on $\partial T$ oriented outwards $T$ we obtain: 
\begin{equation}
\label{concave2}
\left\{
\begin{array}{l}
|a_h^{*}(u^{'},v)-L_h(v)| = c_h^{*}(u^{'},v) + d_h^{*}(u^{'},v) \\
\\
\mbox{where} \\
c_h^{*}(u^{'},v) = \displaystyle \sum_{T \in {\mathcal T}_h} \int_{\partial T} v \frac{\partial u^{'}}{\partial n_T} \\
\mbox{and}\\
d_h^{*}(u^{'},v) = -\displaystyle \sum_{T \in {\mathcal Q}_h} \int_{\Delta^{'}_T} \Delta u^{'} v .
\end{array}
\right.
\end{equation}
$c_h^{*}(u^{'},v)$ can be estimated by means of standard arguments for nonconforming finite elements. More specifically in the case under study (cf. \cite{JJAM}) an estimate of the same nature as \eqref{estimaresidual} applies to $c_h^{*}$, i.e., 
\begin{equation}
\label{concave3} c_h^{*}(u^{'},v) \leq C_R h^2 | u^{'} |_{3,\Omega_h} \parallel {\bf grad}_h v \parallel_{0,h}.
\end{equation}
As for bilinear form $d_h^{*}$ first we observe that, 
\begin{equation}
\label{concave4}
d_h^{*}(u{'},v) \leq \displaystyle \sum_{T \in {\mathcal Q}_h} [volume(\Delta^{'}_T)]^{1/2} \parallel \Delta u^{'} \parallel_{0,\Delta^{'}_T} \parallel v \parallel_{0,\infty,\Delta^{'}_T}.
\end{equation}
Since $\mu_F(v) = 0$ for all face $F$ contained in $\Omega_h$, there exists a mesh-independent constant $C_{\Gamma}^{*}$ such that 
\begin{equation}
\label{concave5}
\parallel v \parallel_{0,\infty,\Delta^{'}_T} \leq 
\parallel v \parallel_{0,\infty,T} \leq C_{\Gamma}^{*} h_T \parallel {\bf grad} \; v \parallel_{0,\infty,T}.
\end{equation} 
Using \eqref{L2TDelta} like in Theorem \ref{P2}, from \eqref{concave5} we thus have 
\begin{equation}
\label{concave6}
\parallel v \parallel_{0,\infty,\Delta^{'}_T} \leq C^{*}_{\Gamma} {\mathcal C}_J h_T^{-1/2} \parallel {\bf grad} \; v \parallel_{0,T}.
\end{equation} 
Noticing that $volume(\Delta^{'}_T)$ is bounded by $h_T^4$ multiplied by a constant $C^{*}_{\Omega}$ depending only on $\Omega$, for both $T \in {\mathcal S}_h \cap {\mathcal Q}_h$ and $T \in {\mathcal R}_h \cap {\mathcal Q}_h$, from straightforward calculations it follows that, 
\begin{equation}
\label{concave7}
\parallel \Delta u^{'} \parallel_{0,\Delta^{'}_T} \leq [C^{*}_{\Omega}]^{1/4} h_T 
\displaystyle \left[ \int_{\Delta^{'}_T} (\Delta u^{'})^4  \right]^{1/4} \; \forall T \in {\mathcal Q}_h.
\end{equation}
Then combining \eqref{concave4}, \eqref{concave5}, \eqref{concave6} and \eqref{concave7}, applying the Cauchy-Schwarz inequality to the summation over $T$,  
and setting $C^{*}_S:= [C^{*}_{\Omega}]^{3/4} C^{*}_{\Gamma} {\mathcal C}_J$ 
we come up with, 
\begin{equation}
\label{concave8}
d_h^{*}(u^{'},v) \leq C_S^{*} h^{2} \displaystyle \left\{ \sum_{T \in {\mathcal Q}_h} h_T \left[ 
\int_{\Delta^{'}_T} (\Delta u^{'})^4  \right]^{1/2} \right\}^{1/2} \parallel {\bf grad}_h v \parallel_{0,h}.
\end{equation}
Applying againn the Cauchy-Schwarz inequality to the summation on the right hand side of \eqref{concave8} 
we readily obtain,
\begin{equation}
\label{concave9}
d_h^{*}(u^{'},v) \leq C_S^{*} h^{2} \displaystyle \left( \sum_{T \in {\mathcal Q}_h} h_T^2  \right)^{1/4} \displaystyle \left[
\sum_{T \in {\mathcal Q}_h} \int_{\Delta^{'}_T} (\Delta u^{'})^4 \right]^{1/4} \parallel {\bf grad}_h v \parallel_{0,h},   
\end{equation}
or yet using \eqref{skinvolume},
\begin{equation}
\label{concave10}
d_h^{*}(u^{'},v) \leq C_S^{*} [C_{\Gamma}^{'}]^{1/2} h^{2} \parallel \Delta u^{'} \parallel_{0,4,\Omega_h}
\parallel {\bf grad}_h v \parallel_{0,h}.   
\end{equation}
Since $H^1(\Omega^{'})$ is continuously embedded in $L^4(\Omega^{'})$ (cf. \cite{Adams}), from \eqref{concave10} we infer the existence of a mesh-independent constant $C^{*}_R$ such that 
\begin{equation}
\label{concave0}
d_h^{*}(u^{'},v) \leq C_R^{*} h^{2} 
\parallel \Delta u^{'} \parallel_{1,\Omega^{'}} \parallel {\bf grad}_h v \parallel_{0,h}, 
\end{equation} 
Now we plug \eqref{concave3} and \eqref{concave0} into \eqref{concave2}, and then the resulting inequality into 
\eqref{concave1}. Finally using the trivial variant of \eqref{interperror*} according to which 
\begin{equation}
\label{interperrorprime}
\parallel {\bf grad}_h (u^{'} - I^{*}_h(u^{'})) \parallel_{0,h} \leq C_P^{'} h^2 | u^{'} |_{3,\Omega^{'}} 
\end{equation}
for a suitable constant $C_P^{'}$ together with the triangle inequality, the result follows. \QED \\

\begin{remark} It is not sure that optimal error esimates in the $L^2$-norm analogous to those given in Theorem \ref{L2P2} apply to the nonconforming finite element studied in this section. One of the reasons for such a shortcoming would be the fact that some properties exploited to estimate the bilinear form $b_{1h}$ no longer hold in the present case. Notice that $L^2$-error estimates applying to this element have not even been established for polyhedral domains. Hence the study of the case of polyhedra is the first step to take in order to carry out an $L^2$-error analysis for curved domains. \QED 
\end{remark} 

\section{Final comments}         
To conclude we make some comments on the methodology studied in this work. \\ 
\begin{enumerate}
\item First of all a word on method's generality. The technique illustrated here in the framework of the solution of the Poisson equation with Dirichlet boundary conditions in curved domains with standard or non standard Lagrange finite elements provides a simple and reliable manner to overcome technical difficulties brought about by more complex problems.
Moreover the principles it is based upon trivially extend to situations of greater complexity than the one of Lagrangian finite elements, in contrast to the isoparametric technique for example.   
For example, Hermite finite element methods to solve second- or fourth-order problems in curved domains with normal-derivative degrees of freedom can also be dealt with very easily by means of our method. This was shown for instance in \cite{AMIS} and in \cite{CFM2017}. 
\item As the reader may have noticed, in case $\Omega$ is a polyhedron the method studied in this paper coincides with the standard Galerkin FEM, as long as the boundary nodes are chosen in the same manner, for instance, the Lagrangian nodes located on $\Gamma$.
\item As for equations with homogeneous Neumann boundary conditions $\partial u /\partial n = 0$ on $\Gamma$ (as long as $f$ satisfies the underlying scalar condition) our method practically coincides with the standard Lagrange finite element method. Indeed, the fact that the degrees of freedom on $\Gamma_h$ are shifted to $\Gamma$ is not supposed to bring about any improvement. However it is well-known that even for the standard method there is order erosion for $k \geq 2$, unless in the variational formulation the domain of integration is taken closer to $\Omega$ than $\Omega_h$. For more details the author refers to \cite{BE2}. Besides this, if inhomogeneous Neumann boundary conditions are prescribed, optimality can only be recovered if the linear form 
$L_h$ is modified, in such a way that boundary integrals for elements $T \in {\mathcal S}_h$ are shifted to a curved boundary approximation sufficiently close to $\Gamma$. But definitively, these are issues that have nothing to do with our method, which is basically aimed at resolving those related to the prescription of degrees of freedom in the case of Dirichlet boundary conditions. 
\item As we should observe our method leads to linear systems of equations with a non-symmetric matrix, even when the original problem is symmetric. Moreover in order to compute the element matrix and right hand side vector for an element in ${\mathcal O}_h$, the inverse of an $n_k \times n_k$ matrix has to be computed. However this represents a rather small extra effort, in view of the significant progress already accomplished in Computational Linear Algebra. 
\item The assumption made throughout the paper that meshes be sufficiently fine (also made by celebrated finite-element authors in the same context) 
is of academic interest only. This assertion is supported by several computations with meshes for which 
$h$ was equal to a half diameter of the domain. Even in such extreme cases the new method behaved pretty well and produced coherent results with respect to successively refined meshes.
\item The use of our method to handle curvilinear boundaries is not restricted to smooth ones. For instance it can also be applied to the case of boundaries of the $C^0$-class consisting of a set of curved faces. Notice that in this case it is advisable to adjust the mesh in such a way that the intersection of adjacent smooth boundary portions are approximated by a polygonal line formed by edges of elements in the mesh, but this procedure is not compulsory. The main constraint is the one of any higher order method: to take the best advantage of the theoretical order the method provides with, the exact solution should be sufficiently smooth. However in general the required regularity will not hold for this type of domains. 
\item As already pointed out in the Introduction, from the author's point of view, an outstanding merit of the new method relies on the use polynomial algebra. This simplifies things significantly especially in the case of complex non linear problems, as compared to methods based on rational functions such as the isoparametric technique. Indeed in the latter case a judicious choice of quadrature formulae to compute element matrices is a must. In contrast exact integration can always be used for this purpose to implement the new method. 
\item Another clear advantage of our method is related to mesh generation, since only straight-edged elements are used. For this reason mesh data structures are as simple as in the case of polyhedral domains. Moreover for some geometries isoparametric elements tend to have locally negative jacobians, which may spoil simulation accuracy with this technique. Clearly enough this situation is completely avoided if our method is employed.
\item A comment is in order on the combination of the new method with widespread techniques to improve accuracy, such as the adaptive finite element method and  
the $h-p$ method. First of all increasing polynomial degree ($p$) is achieved without touching a fixed background mesh, as one can infer from method's description. 
Notice that this is clearly not the case of isoparametric finite elements. Furthermore mesh ($h$) refinement by bisection in the presence of a curvilinear boundary is not more complicated here than for any other method. Of course the same procedure can be exploited for mesh adaptivity. However it is well-known that nested finite-element subspaces are not generated by bisection for Dirichlet boundary conditions prescribed on a curvilinear boundary. This makes application of multigrid methods more tricky, though perfectly possible.
\end{enumerate}     

\indent 
As a conclusion we must say that this is a rather lengthy article 
owing to its basically mathematical content. In spite of this we reported a numerical validation of the methodology it deals with, and showed it to be competitive in terms of both accuracy and CPU time, as compared to existing tetrahedron-based finite-element techniques to handle curvilinear boundaries. But more than this, it provides a simple possibility to tackle the problem in cases where alternatives are either too complicated or simply unknown. It would be interesting to perform as well comparative numerical studies with other higher order methods which do not really belong to the finite-element family, such as finite volumes, discontinuous Galerkin and isogeometric analysis. Fair comparisons with methods not based on tetrahedral cells such as classical spectral elements and, why not, finite differences, could also be accomplished with a straight-edged hexahedron-based version of our method. In short, several perspectives are open for future work. 

\noindent \underline{Acknowledgment:}
The author gratefully acknowledges the financial support provided by CNPq through grant 307996/2008-5. Many thanks are also due to Enrique  Zuazua for helpful discussions.\QED

\section*{APPENDIX - On a mesh-independent upper bound for the Hessian of $f_T$} 
\begin{center}
\end{center}

This Appendix is aimed at presenting convincing arguments leading to the conclusion that the euclidean norm of the Hessian 
${\mathcal H}(f_T)$ of the function $f_T$ defined in Subsection 2.4 is bounded above independently of the element $T$ and the mesh under consideration, as long as the latter is sufficiently fine. We recall that $f_T$ expresses the boundary $\Gamma$ in terms of the coordinates $x$ and $y$ sweeping the planes of the faces of $\Gamma_h$ generically denoted by $F_T$. \\ 
The argument is essentially the same as in problem's two-dimensional counterpart, although in the three-dimensional case additional 
complicated technicalities come into play. Keeping this in mind, we shall first develop the argument in detail as applied to a two-dimensional domain. Then afterwards we just show how it easily extends to the three-dimensional case.\\

Let $\Omega$ be a curved two-dimensional domain of the piecewise $C^2$-class. Here also we consider a regular family of triangular finite-element meshes ${\mathcal T}_h$ fitting $\Omega$ in such a way that all the vertices of the polygon $\Omega_h := \cup_{T \in {\mathcal T}_h} T$ lie on $\Gamma$. We denote by $h_T$ the maximum edge length of $T \in {\mathcal T}_h$. We assume that no mesh in the family has a triangle with more than one edge contained in $\Gamma_h$. 
Let $e_T \subset \Gamma_h$ be the edge of a triangle $T$ having two vertices on $\Gamma$. Denoting by $g_T$ the length of $e_T$,  
let $x$ be the abscissa along $e_T$ in the interval $[0,g_T]$, whose orientation plays no role. Akin to the three-dimensional case we assume that the mesh is fine enough for the portion of $\Gamma$ comprised between the two vertices of $T$ lying on $\Gamma$ to be uniquely 
represented by a function $f_T$ of $x$. More precisely, any point $P$ of such a portion of $\Gamma$ has coordinates $(x,f_T(x))$ in the cartesian coordinate system of the plane $(O,x,y)$, whose origin $O \in \Gamma$ is the afore chosen vertex of $T$.\\
Next we prove that there is a constant $C_{{\mathcal H}}$ independent of both $T$ and the usual mesh parameter (size) $h$, such that,
$|f_T^{''}(x)| \leq C_{{\mathcal H}}$ $\forall x \in (0,g_T)$ and $\forall T \in {\mathcal T}_h$. With this aim we first recall that the  curvature $\kappa$ of $\Gamma$ at $P$ can be locally expressed in terms of $f_T$, in such a way that (see e.g. \cite{Goldman}):

\[ |\kappa(P)| = \displaystyle \frac{|f_T^{''}(x)|}{[1+|f_T^{'}(x)|^2]^{3/2}} \; \forall x \in [0,g_T] \]          

Let ${\mathcal C}_{max} := \max_{P \in \Gamma} |\kappa(P)|$ and ${\mathcal H}_{max}:=\max_{x \in [0,g_T]} |f_T^{''}(x)|$. Since $f_T(0)=f_T(g_T)=0$, there is necessarily an abscissa $x_0 \in [0,g_T]$ at which $f_T^{'}$ vanishes, and hence we can write $|f^{'}_T(x)| = |\int_{x_0}^x f^{''}_T(s) ds|$ for $x \in [0,g_T]$. Then by straightforward calculations we have,

\[ |f_T^{''}(x)|^2 \leq {\mathcal C}_{max}^2 (1+g_T^2 {\mathcal H}_{max}^2)^{3} \; \forall x \in [0,g_T].\]

Now we assume that $g_T \leq \alpha / {\mathcal C}_{max}$, where $\alpha$ is less than or equal to $2 \sqrt{3}/9$.  
This means that the upper bound for ${\mathcal H}_{max}$ we are searching for satisfies,
\[ \displaystyle \frac{{\mathcal H}_{max}^2}{{\mathcal C}_{max}^2} \leq 1 + 3 \alpha^2 \displaystyle \frac{{\mathcal H}_{max}^2}{{\mathcal C}_{max}^2}+3 \alpha^4 \displaystyle \frac{{\mathcal H}_{max}^4}{{\mathcal C}_{max}^4}+\alpha^6 \displaystyle \frac{{\mathcal H}_{max}^6}{{\mathcal C}_{max}^6}.\]

For convenience we set $t:={\mathcal H}_{max}^2/{\mathcal C}_{max}^2$ and $\beta = \alpha^2 (\leq 4/27)$. Next we check whether there exists $t_1 > 0$ such that $0 \leq \varphi(t):= t(1-3\beta-3 \beta^2 t-\beta^3 t^2) \leq 1$ for every $t$ in $[0,t_1]$. Straightforward calculations show that if $\beta \leq 4/27$, the function $\varphi(t)$ is non negative for $0 \leq t \leq t_{max} := (-3 + \sqrt{4 \beta^{-1}-3})/(2 \beta)$ with $t_{max} > 0$.  
Moreover in the interval $[0,t_{max}]$ $\varphi$ attains a minimum at both $t=0$ and $t=t_{max}$ and only a local maximum greater than one at the point $t_0=(-3+\sqrt{\beta^{-1}})/(3 \beta)$. It follows that there exists a point $t_1 \in (0, t_{0})$ depending only on $\beta$ such that $\varphi(t_1) = 1$ and hence $\max_{x \in [0,g_T]} |f_T^{''}(x)| \leq \sqrt{t_1} {\mathcal C}_{max}$. Notice that this upper bound is uniform and holds for all $T$ having an edge on $\Gamma_h$. \\

Next we turn our attention to the case addressed in this work, to which we apply an argument of the same nature as above. More precisely, by these means we conclude that, provided the mesh is sufficiently fine, the euclidean norm of the Hessian of $f_T$ is bounded above by a constant expressed in terms of the principal curvatures of $\Gamma$. However, clearly enough, the calculations are 
considerably more complicated than in the two-dimensional case. That is why we confine ourselves here to sketching the argument. \\
In order to avoid non essential complications, we consider below only faces $F_T$ of the acute type. The case where $F_T$ has obtuse angles can be handled similarly, though at the price of some rather cumbersome modifications.\\  
First of all we refer to formulas (2) and (3) of \cite{Albin}, adapted to the case of a surface with an explicit equation 
$z-\phi(x,y)=0$. Its Gaussian curvature $\kappa_G$ and mean curvature $\kappa_M$ are expressed as follows:

\begin{equation*}  
\left\{
\begin{array}{l}
\kappa_G = \displaystyle \frac{\Delta \phi + \phi_x^2 \phi_{yy} + \phi_y^2 \phi_{xx} - 2 \phi_x \phi_y \phi_{xy}}{2(1+|{\bf grad} \; \phi|^2)^{3/2}}\\
\kappa_m = \displaystyle \frac{\phi_{xx} \phi_{yy}-\phi_{xy}^2 + \phi_x^2 \phi_{yy} + \phi_y^2 \phi_{xx} + \phi_{xx} \phi_{yy} - \phi_{xy}^2}{
(1+|{\bf grad} \; \phi|^2)^2},
\end{array}
\right.
\end{equation*}
where the first and second partial derivatives of $\phi$ with respect to $x$ and $y$ are represented by corresponding subscripts. 
Here we will be dealing with the case of portions of $\Gamma$ for which $\phi=f_T$.\\
Denoting by $s_e$ the abscissa along an edge $e$ of $F_T$, since $f_T$ vanishes at the end-points of $e$, for each $e$ there is a point 
$M_{e,0} \in e$ such that $[\partial f_T/\partial s_e](M_{e,0})=0$. Therefore we can write:
\[ [\partial f_T/\partial s_e](E)= \int_{M_{e,0}}^E \partial^2 f_T/\partial s_e^2 ds_e \; \forall E \in e. \]
This means that $$\|[\partial f_T/\partial s_e](E)| \leq g_T \max_{S \in e}|[\partial^2 f_T/\partial s_e^2](S)| \; \forall E \in e.$$
On the other hand, $\forall M \in F_T^{'}$ we can write $$[\partial f_T/\partial s_e](M) = [\partial f_T/\partial s_e](E_M) + 
\int_{E_M}^M \partial^2 f_T/(\partial \nu_e \partial s_e) d\nu_e,$$ where $\nu_e$ is an abscissa with origin in $e$ along the perpendicular to $e$ through $M$ oriented in an arbitrary way, and $E_M$ is the point of $e$ for which $\nu_e =0$. Taking the previous inequality into account, this implies in turn that $\forall M \in F_T^{'}$, 
$$[\partial f_T/\partial s_e](M) \leq g_T \{\max_{N \in F_T^{'}}|[\partial^2 f_T/\partial s_e^2](N)| + 
\max_{N \in F_T^{'}} |[\partial^2 f_T/(\partial \nu_e \partial s_e)](N)| \}.$$
\\
Since $e$ can be any edge of $F_T$, we readily conclude that there exists a constant ${\mathcal C}_{\Theta}$ depending only on the smallest angle of all faces $F_T$ such that, 
\[ | [{\bf grad}\; f_T](M)| \leq C_{\Theta} g_T \displaystyle \max_{N \in F^{'}_T} |[{\mathcal H}(f_T)](N)| \; \forall M \in F_T^{'}.\]
Now we observe that for every $P \in \Gamma$ covered by the domain $F^{'}_T$ we have, 
$$|4 [\kappa_G(P)]^2 - 2 \kappa_m(P)| = |[{\mathcal H}(f_T)](M)|^2 + {\mathcal D}_T(M)| \{1+|[{\bf grad}\; f_T](M)|^2\}^{-3},$$ 
where $M$ is the point in $F^{'}_T$ corresponding to $P \in \Gamma$, and ${\mathcal D}_T$ denotes a remainder consisting of the sum of products of two or four first order partial derivatives of $f_T$ multiplied by one or the product of two second order partial derivatives of $f_T$. Taking into account the denominators in the expressions of $\kappa_G$ and $\kappa_m$, it is not difficult to figure out that ${\mathcal H}_{max}:= \max_{M \in F^{'}_T}|[{\mathcal H}(f_T)](M)|^2$ 
satisfies, 
\[ {\mathcal H}_{max}^2 \leq {\mathcal C}_{max}^2 [1+({\mathcal C}_{\Theta} g_T {\mathcal H}_{max})^2]^3 + \max_{M \in F^{'}_T}|{\mathcal D}_T(M)|,\] 
with ${\mathcal C}_{max}:= \max_{P \in \Gamma} \sqrt{|4 \kappa_G^2(P) - 2 \kappa_m(P)|}$. \\
Next using the Young's inequality $ab \leq (a^2+b^2)/2$ ($a$ being a first order derivative of $f_T$ and $b$ being the product of a first order derivative and a second order derivative of $f_T$), we note that there exists a mesh-independent constant $C_{\mathcal D}$ such that 
$$\displaystyle \max_{M \in F^{'}_T}|{\mathcal D}_T(M)| \leq C_{\mathcal D} 
\{\displaystyle \max_{M \in F^{'}_T}(|[{\bf grad}\; f_T](M)|^2+|[{\bf grad}\; f_T](M)|^4
+|[{\bf grad}\; f_T](M)|^6 \} (1+{\mathcal H}_{max}^2).$$ 
Then, recalling that $|[{\bf grad}\; f_T](M)| \leq {\mathcal C}_{\Theta} g_T {\mathcal H}_{max}$ $\forall M \in F^{'}_T$, setting 
$t={\mathcal H}_{max}^2/{\mathcal C}_{max}^2$, after trivial adjustments, we come up with 
$$\phi(t):= t[1- p(t)] \leq 1,$$ 
where $p(t)$ is the polynomial $a_1+a_2 t +a_3 t^2 + a_4 t^3$ with strictly positive coefficients $a_i$, $i=1,2,3,4$ expressed in terms of strictly positive powers of $g_T$.\\ 
Now let $\alpha$ be such that $g_T \leq h_T \leq \alpha /{\mathcal C}_{max}$ by assumption. Then the $a_i$s can be made conveniently small by taking $\alpha$ small enough. Finally, except for the degree of $p(t)$ (i.e. $3$ instead of $2$), the remainder of the argument is the same as in the two-dimensional case, thanks to the smallness of the $a_i$s. Otherwise stated, there exists a number $t_1>0$ independent of the meshes in use such that,
\[  \displaystyle \max_{M \in F^{'}_T} |[{\mathcal H}(f_T)](M)| \leq \sqrt{t_1} {\mathcal C}_{max} \; \forall T \in {\mathcal S}_h, \]
where ${\mathcal C}_{max} = \displaystyle \max_{P \in \Gamma} |4 [\kappa_G(P)]^2 - 2 \kappa_m(P)|$, 
$\kappa_G(P)$ and $\kappa_m(P)$ being the Gaussian curvature and the mean curvature of $\Gamma$ at $P$, respectively. \rule{2mm}{2mm}

 \end{document}